\documentclass[11pt,leqno]{article}
\usepackage{amssymb}
\usepackage[mathscr]{eucal}
\usepackage{amsmath,amssymb,latexsym,theorem,bbm}
\usepackage{color,url}
\usepackage{appendix}
\usepackage{accents}

\newcommand{\CC}{\mathbb{C}}

\newcommand{\NN}{\mathbb{N}}

\newcommand{\RR}{\mathbb{R}}

\newcommand{\ZZ}{\mathbb{Z}}

\newcommand{\tA}{\tilde{A}}
\newcommand{\ttA}{\Tilde{\Tilde{A}}}

\newcommand{\tb}{\widetilde{b}}
\newcommand{\tB}{\widetilde{B}}

\newcommand{\bM}{{\boldsymbol{M}}}

\newcommand{\bQ}{{\boldsymbol{Q}}}

\newcommand{\bv}{{\boldsymbol{v}}}

\newcommand{\bx}{{\boldsymbol{x}}}

\newcommand{\bZ}{{\boldsymbol{Z}}}

\newcommand{\bfeta}{{\boldsymbol{\eta}}}

\newcommand{\bpsi}{{\boldsymbol{\psi}}}
\newcommand{\btpsi}{{\widetilde{\bpsi}}}
\newcommand{\bzero}{{\boldsymbol{0}}}

\newcommand{\cA}{{\mathcal A}}
\newcommand{\cB}{{\mathcal B}}

\newcommand{\cD}{{\mathcal D}}

\newcommand{\cF}{{\mathcal F}}

\newcommand{\cL}{{\mathcal L}}

\newcommand{\cN}{{\mathcal N}}
\newcommand{\cP}{{\mathcal P}}

\newcommand{\dd}{\mathrm{d}}
\newcommand{\ee}{\mathrm{e}}
\newcommand{\ii}{\mathrm{i}}

\newcommand{\cont}{\mathrm{cont}}

\DeclareMathOperator*{\argmax}{arg\,max}

\newcommand{\EE}{\operatorname{\mathbb{E}}}
\newcommand{\PP}{\operatorname{\mathbb{P}}}

\newcommand{\hb}{\widehat{b}}

\newcommand{\tbeta}{\widetilde{\beta}}
\newcommand{\ttheta}{\widetilde{\theta}}
\newcommand{\tmu}{\widetilde{\mu}}

\newcommand{\tF}{\widetilde{F}}

\newcommand{\tR}{\widetilde{R}}

\newcommand{\tv}{\widetilde{v}}

\newcommand{\tY}{\widetilde{Y}}

\newcommand{\teta}{\widetilde{\eta}}
\newcommand{\hsigma}{\widehat{\sigma}}

\newcommand{\vare}{\varepsilon}

\renewcommand{\mid}{\,|\,}

\renewcommand{\leq}{\leqslant}
\renewcommand{\geq}{\geqslant}

\newcommand{\stoch}{\stackrel{\PP}{\longrightarrow}}
\newcommand{\distr}{\stackrel{\cD}{\longrightarrow}}
\newcommand{\distre}{\stackrel{\cD}{=}}

\newcommand{\as}{\stackrel{{\mathrm{a.s.}}}{\longrightarrow}}
\newcommand{\ase}{\stackrel{{\mathrm{a.s.}}}{=}}

\newcommand{\bbone}{\mathbbm{1}}
\newcommand{\ns}{{\lfloor ns\rfloor}}
\newcommand{\nks}{{\lfloor n_ks\rfloor}}

\newcommand{\proofend}{\hfill\mbox{$\Box$}}

\numberwithin{equation}{section}

\theoremstyle{change} \theorembodyfont{\em}
\newtheorem{Lem}{Lemma.}[section]
\newtheorem{Thm}[Lem]{Theorem.}
\newtheorem{Pro}[Lem]{Proposition.}

\newtheorem{Def}[Lem]{Definition.}

\theorembodyfont{\rm}
\newtheorem{Rem}[Lem]{Remark.}
\newtheorem{Ex}[Lem]{Example.}

\begin{document}

\begin{center}
 {\bfseries\Large
  Asymptotic properties of maximum likelihood estimator \\[2mm]
   for the growth rate of a stable CIR process \\[2mm]
   based on continuous time observations} \\[7mm]
{\sc\large
 M\'aty\'as $\text{Barczy}^{*,\diamond}$, \ Mohamed $\text{Ben Alaya}^{**}$,}\\[2mm]
 {\sc\large Ahmed $\text{Kebaier}^{***}$ \ and \ Gyula $\text{Pap}^{****}$}
\end{center}

\vskip0.2cm

\noindent
 * MTA-SZTE Analysis and Stochastics Research Group,
   Bolyai Institute, University of Szeged,
   Aradi v\'ertan\'uk tere 1, H--6720 Szeged, Hungary.

\noindent
 ** Laboratoire De Math\'ematiques Rapha\"el Salem, UMR 6085, Universit\'e De Rouen,
    Avenue de L'Universit\'e Technop\^ole du Madrillet, 76801 Saint-Etienne-Du-Rouvray, France.

\noindent
 *** Universit\'e Paris 13, Sorbonne Paris Cit\'e, LAGA, CNRS (UMR 7539),
    Villetaneuse, France.

\noindent
 **** Bolyai Institute, University of Szeged,
     Aradi v\'ertan\'uk tere 1, H--6720 Szeged, Hungary.

\noindent e--mails: barczy@math.u-szeged.hu (M. Barczy), \\
\phantom{e--mails:\,} mohamed.ben-alaya@univ-rouen.fr (M. Ben Alaya), \\
\phantom{e--mails:\,} kebaier@math.univ-paris13.fr (A. Kebaier), \\
\phantom{e--mails:\,} papgy@math.u-szeged.hu (G. Pap).

\noindent $\diamond$ Corresponding author.

\renewcommand{\thefootnote}{}
\vspace*{0.2cm}
\footnote{This research is supported by Laboratory of Excellence MME-DII,
          Grant no.\ ANR11-LBX-0023-01 (\texttt{http://labex-mme-dii.u-cergy.fr/}).
M\'aty\'as Barczy was supported between September 2016 and January 2017 by the
 ''Magyar \'Allami E\"otv\"os \"Oszt\"ond\'{\i}j 2016'' Grant no.\ 75141 funded by the Tempus
 Public Foundation, and from September 2017 by the J\'anos Bolyai Research Scholarship of the
 Hungarian Academy of  Sciences.
Ahmed Kebaier benefited from the support of the chair
 “Risques Financiers”, Fondation du Risque.}

\vspace*{-5mm}

\begin{abstract}
We consider a stable Cox--Ingersoll--Ross process driven by a standard Wiener process and a
 spectrally positive strictly stable L\'evy process, and we study asymptotic properties of the
 maximum
 likelihood estimator (MLE) for its growth rate based on continuous time observations.
We distinguish three cases: subcritical, critical and supercritical.
In all cases we prove strong consistency of the MLE in question, in the subcritical case
 asymptotic normality, and in the supercritical case asymptotic mixed normality are shown as
 well.
In the critical case the description of the asymptotic behavior of the MLE in question remains
 open.
\end{abstract}

\section{Introduction}

We consider a jump-type Cox-Ingersoll-Ross (CIR) process driven by a standard Wiener process
 and a spectrally positive strictly \ $\alpha$-stable L\'evy process given by the SDE
 \begin{align}\label{stable_CIR}
  \dd Y_t = (a - b Y_t) \, \dd t + \sigma \sqrt{Y_t} \, \dd W_t
            + \delta \sqrt[\alpha]{Y_{t-}} \, \dd L_t ,
  \qquad  t \in [0, \infty) ,
 \end{align}
 with an almost surely non-negative initial value \ $Y_0$, \ where \ $a \in [0, \infty)$,
 \ $b \in \RR$, \ $\sigma \in [0, \infty)$, \ $\delta \in (0, \infty)$, \ $\alpha \in (1, 2)$,
 \ $(W_t)_{t \in [0, \infty)}$ \ is a 1-dimensional standard Wiener process, and
 \ $(L_t)_{t\in[0,\infty)}$ \ is a spectrally positive \ $\alpha$-stable L\'evy
 process such that the characteristic function of \ $L_1$ \ takes the form
 \begin{align}\label{char_L}
  \EE(\ee^{\ii\theta L_1})
  = \exp\left\{\int_0^\infty
                (\ee^{\ii \theta z} - 1 - \ii \theta z) C_\alpha z^{-1-\alpha}
                \, \dd z\right\} ,
  \qquad \theta \in \RR ,
 \end{align}
 where \ $C_\alpha := (\alpha \Gamma(-\alpha))^{-1}$ \ and \ $\Gamma$ \ denotes the Gamma
 function.
In fact, \ $(L_t)_{t\in[0,\infty)}$ \ is a strictly \ $\alpha$-stable L\'evy process,
 see, e.g., Sato \cite[part (vi) of Theorem 14.7]{Sat}.
We suppose that \ $Y_0$, \ $(W_t)_{t\in[0,\infty)}$ \ and \ $(L_t)_{t\in[0,\infty)}$ \ are
 independent.
Under the given conditions together with \ $\EE(Y_0) < \infty$, \ there is a (pathwise) unique
 strong solution of the SDE \eqref{stable_CIR} with
 \ $\PP(\text{$Y_t \in [0, \infty)$ for all $t \in [0, \infty)$}) = 1$.
\ As a matter of fact, the SDE \eqref{stable_CIR} is a special case of the SDE (1.8) in Fu and
 Li \cite{FuLi} (with the special choice \ $z_1 \equiv 0$), \ for which the existence of a
 pathwise unique non-negative strong solution has been proved (see Fu and Li
 \cite[Corollary 6.3]{FuLi}).
Eventually, the process \ $(Y_t)_{t\in[0,\infty)}$ \ given by the SDE \eqref{stable_CIR} is a
 continuous state and continuous time branching process with immigration (CBI process), see (ii)
 of Proposition \ref{Pro_stable_CIR}.
We call \ $Y$ \ an $\alpha$-stable CIR process (or Alpha-CIR process), which is a
 generalization of the usual CIR process (given by the SDE \eqref{stable_CIR} formally with
 \ $\delta = 0$).

Stable CIR processes become more and more popular in stochastic modelling, and it is an
 interesting class of CBI processes on its own right as well.
Carr and Wu \cite[equation (31)]{CarWu} considered a stochastic process admitting an
 infinitesimal generator which coincides with the corresponding one of an $\alpha$-stable CIR
 process with \ $\sigma = 0$, \ see (iv) of Proposition \ref{Pro_stable_CIR}.

Li and Ma \cite{LiMa} proved exponential ergodicity for the process \ $(Y_t)_{t\in[0,\infty)}$
 \ provided that \ $a \in (0, \infty)$ \ and \ $b \in (0, \infty)$, \ for more details,
 see (ii) of Theorem \ref{Ergodicity}.
Li and Ma \cite{LiMa} also described the asymptotic behavior of the conditional least squares
 estimator (LSE) and weighted conditional LSE of the drift parameters \ $(a, b)$ \ of an
 $\alpha$-stable CIR process given by the SDE \eqref{stable_CIR} with \ $\sigma = 0$, \ based on
 (discretely observed) low frequency observations in the subcritical case (i.e., when
 \ $b \in (0, \infty)$).
\ In the region \ $\alpha \in (1, \frac{1+\sqrt{5}}{2})$, \ Li and Ma \cite{LiMa} showed that
 the normalizing factor for the LSE of \ $(a, b)$ \ is \ $n^{(\alpha-1)/\alpha^2}$, \ which is
 quite different from the \ $\sqrt{n}$-normalization being quite usual for subcritical models.
On the top of it all, Li and Ma \cite{LiMa} also proved that the corresponding normalizing
 factor for the weighted LSE of \ $(a, b)$ \ is \ $n^{(\alpha-1)/\alpha}$ \ (being different
 from the one for the (usual) LSE) in the whole region \ $\alpha \in (1, 2)$.

Jiao et al.~\cite{JiaMaSco} investigated several properties of $\alpha$-stable CIR processes
 such as integral representations, branching property in the pathwise sense, necessary and
 sufficient conditions for strictly positiveness and they made an analysis of the jumps of the
 process.
Further, they used $\alpha$-stable CIR processes for interest rate modelling and pricing by
 pointing out that these processes can describe some recent phenomena on sovereign bond market
 such as large fluctuations at a local extent together with the usual small oscillations, for
 more details, see the Introduction of Jiao et al.~\cite{JiaMaSco}.
Very recently, Jiao et al.~\cite{JiaMaScoSga} have proposed concrete examples of applications and investigated a factor model for
 electricity prices, where $\alpha$-stable CIR processes may appear as factors of the model in
 question.

Peng \cite{Pen} introduced and studied a so-called $\alpha$-stable CIR process with restart,
 by which one means that the process in question behaves as an $\alpha$-stable CIR process
 given by the SDE \eqref{stable_CIR} with \ $\sigma = 0$, \ it is killed at the boundary \ $0$
 \ of \ $[0, \infty)$, \ and according to an exponential clock it jumps to a new point in
 \ $[0, \infty)$ \ according to a given probability distribution on \ $[0, \infty)$.
\ As it was pointed out in Peng \cite{Pen}, restart phenomenon appears in internet congestion
 as well: whenever a web page takes too much time to appear, it is useful to press the reload
 button and then usually the web page appears immediately.

Yang \cite{Yan} studied $\alpha$-stable CIR processes with small \ $\alpha$-stable noises
 given by the SDE
 \begin{align}\label{stable_CIR_small}
  \dd Y^\vare_t
  = (a - b Y^\vare_t) \, \dd t + \delta \vare \sqrt[q]{Y^\vare_{t-}} \, \dd L_t ,
  \qquad  t \in [0, \infty) ,
 \end{align}
 with a non-negative deterministic initial value \ $Y^\vare_0 = y_0 \in [0, \infty)$, \ where
 \ $q \in \bigl(0, \frac{1}{1-1/\alpha}\bigr)$ \ and \ $\vare \in (0,\infty)$.
\ The asymptotic behavior of an approximate maximum likelihood estimator (MLE) of
 \ $(a, b, \delta)$ \ has been described based on discrete time observations at \ $n$
 \ regularly spaced time points \ $\frac{k}{n}$, \ $k \in \{1, \ldots, n\}$, \ on a fixed
 time interval \ $[0, 1]$.
\ Tending \ $\vare$ \ to \ $0$ \ and \ $n \to \infty$ \ at a given rate, for some restricted
 parameter set, Yang \cite[Theorem 2.4]{Yan} proved asymptotic normality of the approximate
 MLE in question.
In some sense it is surprising, since this restricted parameter set contains parameters
 belonging to critical \ ($b = 0$) \ and supercritical \ ($b \in (-\infty, 0)$) \ models
 as well both with normal limit distributions, and for critical models, the limit distribution, in general,
  is not even mixed normal.

Ma and Yang \cite{MaYan} investigated asymptotic behavior of the LSE of \ $a$ \ for the model
 \eqref{stable_CIR_small} (all the other parameters are supposed to be known) based on
 discrete time observations as in Yang \cite{Yan} described above.
They described the asymptotic behavior of the LSE in question and derived large and moderate
 deviation inequalities for it as well, see Ma and Yang \cite[Theorems 2.1, 2.3--2.5]{MaYan}.

In this paper, supposing that \ $a \in [0, \infty)$, \ $\sigma, \delta \in (0, \infty)$ \ and
 \ $\alpha \in (1, 2)$ \ are known, we study the asymptotic properties of the MLE of
 \ $b \in \RR$ \ based on continuous time observations \ $(Y_t)_{t\in[0,T]}$ \ with
 \ $T \in (0, \infty)$, \ starting the process \ $Y$ \ from some known non-random initial value
 \ $y_0 \in [0, \infty)$.

The paper is organized as follows.
Section \ref{Prel} is devoted to some preliminaries.
First, we recall some useful properties of the stable CIR process \ $(Y_t)_{t\in[0,\infty)}$
 \ given by the SDE \eqref{stable_CIR} such as the existence of a non-negative pathwise unique
 strong solution, the forms of the Laplace transform and the infinitesimal generator or
 conditions on the strictly positiveness of the process or the integrated process, see
 Proposition \ref{Pro_stable_CIR}.
We derive a so-called Grigelionis form of the semimartingale \ $(Y_t)_{t\in[0,\infty)}$, \ see
 Proposition \ref{Grigelionis}.
Based on the asymptotic behavior of the expectation of \ $Y_t$ \ as \ $t \to \infty$, \ we
 distinguish subcritical, critical or supercritical cases according to \ $b \in (0, \infty)$,
 \ $b = 0$ \ or \ $b \in (-\infty, 0)$, \ see Proposition \ref{Pro_moments} and Definition
 \ref{Def_criticality}.
In Proposition \ref{Pro_moments} it also turns out that the parameter \ $b$ \ can be interpreted
 as a growth rate of the model.
We recall a result about the existence of a unique stationary distribution for the process
 \ $(Y_t)_{t\in[0,\infty)}$ \ in the subcritical and critical cases, and about its exponential
 ergodicity in the subcritical case, due to Li \cite{Li}, Li and Ma \cite{LiMa} and
 Jin et al.\ \cite{JinKreRud}, see Theorem \ref{Ergodicity}.
We call the attention that there exists a unique stationary distribution for
 \ $(Y_t)_{t\in[0,\infty)}$ \ in the critical case as well.
Remark \ref{Rem_ergodic_alternative_proof} is devoted to give an alternative proof for
 the weak convergence of \ $Y_t$ \ as \ $t \to \infty$ \ in Theorem \ref{Ergodicity}
 in case of \ $\sigma \in (0,\infty)$, \ giving more insight as well.
In Remark \ref{Thm_MLE_cons_sigma}, we give a statistic for \ $\sigma^2$ \ using continuous time
 observations \ $(Y_t)_{t\in[0,T]}$ \ with an arbitrary \ $T \in (0, \infty)$, \ and due to this
 result we do not consider the estimation of the parameter \ $\sigma$, \ it is supposed to be
 known.
In Section \ref{section_Laplace_crit}, we derive a formula for the joint Laplace transform of
 \ $Y_t$ \ and \ $\int_0^t Y_s \, \dd s$, \ where \ $t \in [0, \infty)$, \ using Theorem 4.10 in
 Keller-Ressel \cite{KR2}, see Theorem \ref{Thm_Laplace}.
We note that this form of the joint Laplace transform in question is a consequence of Theorem 5.3
 in Filipovi\'{c} \cite{Fil}, a special case of Proposition 3.3 in Jiao et al.\ \cite{JiaMaSco} as
 well, and it is used for describing the
 asymptotic behavior of the MLE of \ $b$ \ in question in the critical and supercritical cases.
Section \ref{section_EUMLE} is devoted to prove the existence and uniqueness of the MLE of \ $b$
 \ (provided that \ $\sigma \in (0, \infty)$) \ deriving an explicit formula for it as
 well, see Proposition \ref{LEMMA_MLEb_exist}.
In Remark \ref{Rem_L_measurability}, under the additional assumption
 \ $a \in \bigl[\frac{\sigma^2}{2}, \infty\bigr)$, \ we prove that \ $L_t$ \ is a
 measurable function of \ $(Y_u)_{u\in[0,T]}$ \ for all \ $t \in [0, T]$ \ with any
 \ $T \in (0, \infty)$.
\ In Section \ref{section_MLE_subcritical}, provided that \ $a \in (0, \infty)$, \ we
 prove strong consistency and asymptotic normality of the MLE of \ $b$ \ in the subcritical
 case, see Theorem \ref{Thm_MLEb_subcritical}.
The asymptotic normality in question holds with a usual square root normalization \ ($\sqrt{T}$),
 \ but as usual, the asymptotic variance depends on the unknown parameter \ $b$, \ as well.
To get around this problem, we also replace the normalization \ $\sqrt{T}$ \  by a random one
 \ $\frac{1}{\sigma}\left(\int_0^t Y_s\,\dd s\right)^{1/2}$ \ (depending only on the observation,
 but not on the parameter \ $b$) \ with the advantage that the MLE of \ $b$ \ with this random scaling
 is asymptotically standard normal,
 so one can give asymptotic confidence intervals for the unknown parameter \ $b$, \ which is desirable for practical purposes.
\ Section \ref{section_MLE_critical} is devoted to prove the strong consistency of the MLE of
 \ $b$ \ in the critical case, provided that \ $a \in (0, \infty)$, \ (see Theorem
 \ref{Thm_MLE_critical}) using the limit behavior of the unique locally bounded solution of the
 differential equation \eqref{psi_DE} at infinity described in Proposition
 \ref{Pro_critical_psi_asymp}.
We call the attention to the fact that for the \ $\alpha$-stable CIR process \ $(Y_t)_{t\in[0,\infty)}$,
 \ the critical case \ ($b = 0$) \ is somewhat special (compared to the original CIR process
 with \ $b = 0$), \ since there still exists a unique stationary distribution for
 \ $(Y_t)_{t\in[0,\infty)}$, \ however its expectation is infinite unless \ $a = 0$ \ (see
 Theorem \ref{Ergodicity}), and surprisingly, we can prove strong consistency of the MLE in
 question not only weak consistency usually proved for critical models.
In the critical case the description of the asymptotic behavior of the MLE remains open.
In Section \ref{section_MLE_supercritical}, for the supercritical case, provided that
 \ $a \in (0, \infty)$, \ we prove that the MLE of \ $b$ \ is strongly consistent and
 asymptotically mixed normal with the deterministic scaling \ $\ee^{-bT/2}$, \
 and it is asymptotically standard normal with the random scaling \ $\frac{1}{\sigma}\left(\int_0^t Y_s\,\dd s\right)^{1/2}$,
 \ see Theorem \ref{Thm_MLE_supercritical}.
We point out that the limit mixed normal law in question is characterized in a somewhat
 complicated way, namely in its description a positive random variable \ $V$ \ comes into play
 of which the Laplace transform contains a function related to the branching mechanism of the
 CBI process \ $(Y_t)_{t\in[0,\infty)}$, \ see Theorem \ref{Thm_supercritical_convergence}.
We give two proofs for the derivation of the Laplace transform of \ $V$, \ and the second one is
 heavily based on the general theory of CBI processes, for which we will refer to Li \cite{Li}.
We close the paper with three Appendices, where we recall certain sufficient conditions for the
 absolute continuity of probability measures induced by semimartingales together with a
 representation of the Radon--Nikodym derivative (Appendix \ref{App_LR}), some limit theorems
 for continuous local martingales (Appendix \ref{App_clm}) and in case of a $\frac{3}{2}$-stable
 CIR process we present some explicit formulae for the Laplace transform of the unique
 stationary distribution in the subcritical and critical cases, of \ $Y_t$,
 \ $t \in [0, \infty)$, \ in all the cases of \ $b \in \RR$, \ and of \ $V$ \ in the
 supercritical case, respectively (Appendix \ref{App_alpha=3/2}).

Finally, we summarize the novelties of the paper.
According to our knowledge, maximum likelihood estimation based on continuous time observations
 has never been studied before for the \ $\alpha$-stable CIR process \ $(Y_t)_{t\in[0,\infty)}$,
 \ and since these processes become more and more popular in financial mathematics and market
 models for electricity prices, the problem of estimating its parameters is an important
 question as well.
Further, in the critical case, somewhat surprisingly, we can prove strong consistency of the MLE
 of \ $b$, \ which can be considered as a new phenomenon, since for other critical financial
 models, such as for the usual CIR process or for the Heston process, only weak consistency is
 proved in the critical case, see Overbeck \cite[Theorem 2, parts (iii) and (iv)]{Ove} and
 Barczy and Pap \cite[Remark 4.4]{BarPap2}, respectively.

\section{Preliminaries}
\label{Prel}

Let \ $\NN$, \ $\ZZ_+$, \ $\RR$, \ $\RR_+$, \ $\RR_{++}$, \ $\RR_-$, \ $\RR_{--}$ \ and \ $\CC$
 \ denote the sets of positive integers, non-negative integers, real numbers, non-negative real
 numbers, positive real numbers, non-positive real numbers, negative real numbers and complex
 numbers, respectively.
For \ $x, y \in \RR$, \ we will use the notations \ $x \land y := \min(x, y)$ \ and
\ $x \lor y := \max(x, y)$.
\ The integer part of a real number \ $x \in \RR$ \ is denoted by \ $\lfloor x \rfloor$.
\ By \ $\|x\|$ \ and \ $\|A\|$, \ we denote the Euclidean norm of a vector \ $x \in \RR^d$ \ and
 the induced matrix norm of a matrix \ $A \in \RR^{d \times d}$, \ respectively.
By \ $\cB(\RR_+)$, \ we denote the Borel $\sigma$-algebra on \ $\RR_+$.
\ We will denote the convergence in probability, in distribution and almost surely, and almost
 sure equality by \ $\stoch$, \ $\distr$, \ $\as$ \ and \ $\ase$, \ respectively.
By \ $C^2_c(\RR_+, \RR)$ \ and \ $C^{\infty}_c(\RR_+, \RR)$, \ we denote the set of twice
 continuously differentiable real-valued functions on \ $\RR_+$ \ with compact support and the
 set of infinitely differentiable real-valued functions on \ $\RR_+$ \ with compact support,
 respectively.

Let \ $\big(\Omega, \cF, (\cF_t)_{t\in\RR_+}, \PP\big)$ \ be a filtered probability space
 satisfying the usual conditions, i.e.,  \ $(\Omega,\cF,\PP)$ \ is complete, the filtration
 \ $(\cF_t)_{t\in\RR_+}$ \ is right-continuous, \ $\cF_0$ \ contains all the $\PP$-null sets
 in \ $\cF$, \ and \ $\cF = \sigma\big(\bigcup_{t\in\RR_+} \cF_t\big)$.
\ Let \ $(W_t)_{t\in\RR_+}$ \ be a standard Wiener process with respect to the filtration
 \ $(\cF_t)_{t\in\RR_+}$, \ and \ $(L_t)_{t\in\RR_+}$ \ be a spectrally positive strictly
 \ $\alpha$-stable L\'evy process with respect to the filtration \ $(\cF_t)_{t\in\RR_+}$ \ such
 that the characteristic function of \ $L_1$ \ is given by \eqref{char_L}.
We assume that \ $W$ \ and \ $L$ \ are independent.
Recall that the L\'evy-It\^{o}'s representation of \ $L$ \ takes the form
 \begin{align}\label{help_Levy_Ito_decomposition}
  L_t = \int_{(0,t]} \int_{(0,\infty)} z \, \tmu^L(\dd s,\dd z)
      = \gamma t + \int_{(0,t]} \int_{(0, 1]} z \, \tmu^L(\dd s,\dd z)
                 + \int_{(0,t]} \int_{(1, \infty)} z \, \mu^L(\dd s,\dd z)
 \end{align}
 for \ $t \in \RR_+$, \ where
 \ $\mu^L(\dd s,\dd z)
    := \sum_{u\in\RR_+} \bbone_{\{\Delta L_u\ne 0\}} \, \vare_{(u,\Delta L_u)}(\dd s,\dd z)$
 \ is the integer-valued Poisson random measure on \ $\RR_{++}^2$ \ associated with the jumps
 \ $\Delta L_u := L_u - L_{u-}$, \ $u \in \RR_{++}$, \ $\Delta L_0 := 0$, \ of the process
 \ $L$, \ and \ $\vare_{(u,x)}$ \ denotes the Dirac measure at the point \ $(u, x) \in \RR_+^2$,
 \ $\tmu^L(\dd s,\dd z) := \mu^L(\dd s,\dd z) - \dd s \, m(\dd z)$, \ where
 \ $m(\dd z) := C_\alpha z^{-1-\alpha} \bbone_{(0,\infty)}(z) \, \dd z$, \ and
 \ $\gamma := - \int_{(1, \infty)} z \, \dd s \, m(\dd z)
            = - C_\alpha \int_1^\infty z^{-\alpha} \, \dd z = \frac{C_\alpha}{1-\alpha}$.
\ The measure \ $m$ \ is nothing else but the L\'evy measure of \ $L$.
 \ We also note that \ $(L_t)_{t\in\RR_+}$ \ is a martingale and consequently \ $\EE(L_t) = 0$,
 \ $t\in\RR_+$.

The next proposition is about the existence and uniqueness of a strong solution of the SDE
 \eqref{stable_CIR} stating also that \ $Y$ \ is a CBI process with explicitly given branching
 and immigration mechanisms and we also collect some other useful properties \ of \ $Y$ \ based
 on Dawson and Li \cite{DawLi}, Fu and Li \cite{FuLi}, Li \cite{Li} and
 Jiao et al.\ \cite{JiaMaSco}.

\begin{Pro}\label{Pro_stable_CIR}
Let \ $\eta_0$ \ be a random variable independent of \ $(W_t)_{t\in\RR_+}$ \ and
 \ $(L_t)_{t\in\RR_+}$ \ satisfying \ $\PP(\eta_0 \in \RR_+) = 1$ \ and
 \ $\EE(\eta_0) < \infty$.
\ Let \ $a \in \RR_+$, \ $b \in \RR$, \ $\sigma\in\RR_+$, \ $\delta \in \RR_{++}$, \ and \ $\alpha\in(1,2)$.
\ Then the following statements hold.
 \renewcommand{\labelenumi}{{\rm(\roman{enumi})}}
 \begin{enumerate}
  \item
   There exists a pathwise unique strong solution \ $(Y_t)_{t\in\RR_+}$ \ of the SDE
    \eqref{stable_CIR} such that \ $\PP(Y_0 = \eta_0) = 1$ \ and
    \ $\PP(\text{$Y_t \in \RR_+$ \ for all \ $t \in \RR_+$}) = 1$.
  \item
   The process \ $(Y_t)_{t\in\RR_+}$ \ is a CBI process having branching mechanism
    \[
      R(z) = \frac{\sigma^2}{2} z^2 + \frac{\delta^\alpha}{\alpha}z^\alpha + bz ,
      \qquad z \in \RR_+ ,
    \]
    and immigration mechanism
    \[
      F(z) = a z , \qquad z \in \RR_+ .
    \]
  \item
   For all \ $t \in \RR_+$ \ and \ $y_0 \in \RR_+$, \ the Laplace transform of \ $Y_t$ \ takes
    the form
    \begin{align}\label{Laplace_Y_t1}
      \EE(\ee^{-\lambda Y_t} \mid Y_0 = y_0)
      &= \exp\left\{- y_0 v_t(\lambda) - \int_0^t F(v_s(\lambda)) \, \dd s \right\}
    \end{align}
    for all \ $\lambda \in \RR_+$, \ where \ $\RR_+ \ni t \mapsto v_t(\lambda) \in \RR_+$ \ is
     the unique locally bounded solution to
     \begin{equation}\label{LiLaplace}
      \frac{\partial}{\partial t} v_t(\lambda) = - R(v_t(\lambda)) , \qquad
      v_0(\lambda) = \lambda .
     \end{equation}
   If \ $t \in \RR_+$, \ $y_0 \in \RR_+$ \ and \ $\lambda \in \RR_{++} \setminus \{\theta_0\}$
    \ with \ $\theta_0 := \inf\{z \in \RR_{++} : R(z) \in \RR_+\} \in \RR_+$, \ then we have
    \begin{equation}\label{Laplace_Y_t2}
     \EE(\ee^{-\lambda Y_t} \mid Y_0 = y_0)
     = \exp\left\{- y_0 v_t(\lambda)
                  + \int_\lambda^{v_t(\lambda)} \frac{F(z)}{R(z)} \, \dd z \right\} .
    \end{equation}
   Especially, \eqref{Laplace_Y_t2} holds for all \ $\lambda \in \RR_{++}$ \ whenever
    \ $b \in \RR_+$.
  \item
   The infinitesimal generator of \ $Y$ \ takes the form
    \begin{align}\label{infgen}
     \begin{split}
      (\cA f)(y) = (a - b y) f'(y) + \frac{\sigma^2}{2} y f''(y)
                   + \delta^\alpha y
                     \int_0^\infty
                      \Big(f(y + z) - f(y) - z f'(y)\Big) C_\alpha z^{-1-\alpha} \, \dd z ,
     \end{split}
    \end{align}
    where \ $y \in \RR_+$, \ $f \in C^2_c(\RR_+, \RR)$, \ and \ $f'$ \ and \ $f''$ \ denote the
    first and second order partial derivatives of \ $f$.
  \item
   If, in addition, \ $\PP(\eta_0 \in \RR_{++}) = 1$ \ or \ $a \in \RR_{++}$, \ then
    \ $\PP\bigl(\int_0^t Y_s \, \dd s \in \RR_{++}\bigr) = 1$ \ for all \ $t \in \RR_{++}$.
  \item
   If, in addition, \ $\sigma \in \RR_{++}$ \ and \ $a \geq \frac{\sigma^2}{2}$, \ then
    \ $\PP\bigl(\text{$Y_t \in \RR_{++}$ \ for all \ $t \in \RR_{++}$}\bigr) = 1$.
  \item
   If, in addition, \ $\PP(\eta_0 \in \RR_{++})=1$, \ $a = 0$ \ and \ $b \in \RR_+$, \ then
    \ $\PP(\tau_0 < \infty) = 1$, \ where \ $\tau_0 := \inf\{s \in \RR_+ : Y_s = 0\}$, \ and
    \ $\PP(\text{$Y_t = 0$ \ for all \ $t \geq \tau_0$}) = 1$.
 \end{enumerate}
\end{Pro}

\noindent{\bf Proof.}
For the existence of a pathwise unique non-negative strong solution satisfying
 \ $\PP(Y_0 = \eta_0) = 1$ \ and \ $\PP(\text{$Y_t \in \RR_+$ \ for all \ $t \in \RR_+$}) = 1$,
  \ see Fu and Li \cite[Corollary 6.3]{FuLi}, which yields (i).

Further, Theorem 6.2 in Dawson and Li \cite{DawLi} together with
 \[
   \int_0^\infty (z \land z^2 ) C_\alpha z^{-1-\alpha} \, \dd z
   = C_\alpha \int_0^1 z^{1-\alpha} \, \dd z + C_\alpha \int_1^\infty z^{-\alpha} \, \dd z
   = C_\alpha \left(\frac{1}{2-\alpha} + \frac{1}{\alpha-1}\right)
   < \infty
 \]
 and
 \begin{align}\label{help_branching_mechanism}
  \int_0^\infty (\ee^{-zx} - 1 + z x) C_\alpha x^{-1-\alpha} \, \dd x
  = \frac{1}{\alpha} \frac{\alpha(\alpha-1)}{\Gamma(2-\alpha)}
    \int_0^\infty (\ee^{-zx} - 1 + z x) x^{-1-\alpha} \, \dd x
  = \frac{z^\alpha}{\alpha}
 \end{align}
 for \ $z \in \RR_+$ \ (see, e.g., Li \cite[Example 1.9]{Li}) imply that \ $Y$ \ is a CBI
 process having branching and immigration mechanisms given in (ii).

For formula \eqref{Laplace_Y_t1} and, in case of \ $b \in \RR_+$, \ formula \eqref{Laplace_Y_t2}
 see Li \cite[formula (3.29) and page 67]{Li}.
Next we check that
 \begin{equation}\label{substitution}
  - \int_0^t F(v_s(\lambda)) \, \dd s = \int_\lambda^{v_t(\lambda)} \frac{F(z)}{R(z)} \, \dd z
 \end{equation}
 for all \ $t \in \RR_+$ \ and \ $\lambda \in \RR_{++} \setminus \{\theta_0\}$.
\ It is enough to verify that the continuously differentiable function
 \ $(0, t) \ni s \mapsto v_s(\lambda)$ \ is strictly monotone for all
 \ $\lambda \in \RR_{++} \setminus \{\theta_0\}$, \ since then, by the substitution
 \ $z = v_s(\lambda)$, \ we obtain
 \[
   - \int_0^t F(v_s(\lambda)) \, \dd s
   = - \int_\lambda^{v_t(\lambda)}
        \frac{F(z)}{\frac{\partial}{\partial s}v_s(\tv_z(\lambda))} \, \dd z
   = \int_\lambda^{v_t(\lambda)}
      \frac{F(z)}{R(v_s(\tv_z(\lambda)))} \, \dd z
 \]
 and hence \eqref{substitution}, where
 \ $(v_0(\lambda) \land v_t(\lambda), v_0(\lambda) \lor v_t(\lambda)) \ni z \mapsto \tv_z(\lambda)$
 \ denotes the inverse of \ $(0, t) \ni s \mapsto v_s(\lambda)$.
\ By Li \cite[Proposition 3.1]{Li}, the function
 \ $\RR_+ \ni \lambda \mapsto v_s(\lambda) \in \RR_+$ \ is strictly increasing for all
 \ $s \in \RR_+$.
\ We have \ $v_s(\theta_0) = \theta_0$ \ for all \ $s \in \RR_+$, \ since \ $R(\theta_0) = 0$
 \ yields that this constant function is the unique locally bounded solution to the differential
 equation \eqref{LiLaplace} with initial value \ $\theta_0$.
\ If \ $b \in \RR_+$, \ then \ $\theta_0 = 0$, \ thus \ $\lambda \in \RR_{++}$ \ implies
 \ $v_s(\lambda) > v_s(0) = 0$ \ for all \ $s \in \RR_+$.
\ In this case, using the differential equation \eqref{LiLaplace} and the inequality
 \ $R(z) > 0$ \ for all \ $z \in \RR_{++}$, \ we obtain
 \ $\frac{\partial}{\partial s}v_s(\lambda) = -R(v_s(\lambda)) < 0$ \ for all \ $s \in \RR_+$,
 \ hence the function \ $(0, t) \ni s \mapsto v_s(\lambda)$ \ is strictly decreasing, thus we
 conclude \eqref{substitution} for \ $b \in \RR_+$.
\ If \ $b \in \RR_{--}$, \ then \ $\theta_0 \in \RR_{++}$.
\ Consequently, in case of \ $b \in \RR_{--}$ \ and \ $\lambda \in (0, \theta_0)$ \ we have
 \ $v_s(\lambda) < v_s(\theta_0) = \theta_0$ \ for all \ $s \in \RR_+$.
\ In this case, using the differential equation \eqref{LiLaplace} and the inequality
 \ $R(z) < 0$ \ for all \ $z \in (0, \theta_0)$, \ we obtain
 \ $\frac{\partial}{\partial s}v_s(\lambda) = -R(v_s(\lambda)) > 0$ \ for all \ $s \in \RR_+$,
 \ hence the function \ $(0, t) \ni s \mapsto v_s(\lambda)$ \ is strictly increasing, thus we
 conclude \eqref{substitution} for \ $b \in \RR_{--}$ \ and \ $\lambda \in (0, \theta_0)$.
In a similar way, in case of \ $b \in \RR_{--}$ \ and \ $\lambda \in (\theta_0, \infty)$ \ we
 have \ $v_s(\lambda) > v_s(\theta_0) = \theta_0$ \ for all \ $s \in \RR_+$.
\ In this case, using the differential equation \eqref{LiLaplace} and the inequality
 \ $R(z) > 0$ \ for all \ $z \in (\theta_0, \infty)$, \ we obtain
 \ $\frac{\partial}{\partial s}v_s(\lambda) = -R(v_s(\lambda)) < 0$ \ for all \ $s \in \RR_+$,
 \ hence the function \ $(0, t) \ni s \mapsto v_s(\lambda)$ \ is strictly decreasing, thus we
 conclude \eqref{substitution} for \ $b \in \RR_{--}$ \ and \ $\lambda \in (\theta_0, \infty)$
 \ as well.

The form of the infinitesimal generator \eqref{infgen} can be checked similarly as in the
 proof of Theorem 2.1 of Barczy et al.\ \cite{BarDorLiPap2}, implying (iv).

For (v), let us fix \ $t \in \RR_{++}$ \ and put
 \[
   A_t := \bigl\{\text{$\omega \in \Omega$
                   : $[0, t] \ni s \mapsto Y_s(\omega)$ \ is c\`adl\`ag and
                     \ $Y_s(\omega) \in \RR_+$ \ for all \ $s \in [0, t]$}\bigr\} .
 \]
Then, by (i), \ $\PP(A_t) = 1$ \ and for all \ $\omega \in A_t$,
 \ $\int_0^t Y_s(\omega) \, \dd s = 0$ \ if and only if \ $Y_s(\omega) = 0$ \ for all
 \ $s \in [0, t)$.
\ By \eqref{stable_CIR},
 \[
   Y_s = Y_0 + a s - b \int_0^s Y_u \, \dd u + \sigma \int_0^s \sqrt{Y_u} \, \dd W_u
         + \delta \int_0^s \sqrt[\alpha]{Y_{u-}} \, \dd L_u , \qquad \ s \in \RR_+ ,
 \]
 holds \ $\PP$-almost surely.
The stochastic integrals on the right hand side can be approximated as
 \begin{align*}
  &\sup_{s\in[0,t]}
    \Biggl|\sum_{i=1}^\ns \sqrt{Y_{\frac{i-1}{n}}} (W_{\frac{i}{n}} - W_{\frac{i-1}{n}})
          - \int_0^s \sqrt{Y_u} \, \dd W_u\Biggr|
   \stoch 0 \qquad \text{as}\qquad n \to \infty , \\
  &\sup_{s\in[0,t]}
    \Biggl|\sum_{i=1}^\ns
           \sqrt[\alpha]{Y_{\frac{i-1}{n}-}} (L_{\frac{i}{n}} - L_{\frac{i-1}{n}})
          - \int_0^s \sqrt[\alpha]{Y_{u-}} \, \dd L_u\Biggr|
   \stoch 0 \qquad \text{as}\qquad n \to \infty  ,
 \end{align*}
 see Jacod and Shiryaev \cite[Theorem I.4.44]{JSh}.
Hence there exists a sequence \ $(n_k)_{k\in\NN}$ \ of positive integers such that
 \begin{align*}
  &\sup_{s\in[0,t]}
    \Biggl|\sum_{i=1}^\nks \sqrt{Y_{\frac{i-1}{n_k}}} (W_{\frac{i}{n_k}} - W_{\frac{i-1}{n_k}})
          - \int_0^s \sqrt{Y_u} \, \dd W_u\Biggr|
    \to 0 \qquad \text{as}\qquad k \to \infty , \\
  &\sup_{s\in[0,t]}
    \Biggl|\sum_{i=1}^\nks
           \sqrt[\alpha]{Y_{\frac{i-1}{n_k}-}} (L_{\frac{i}{n_k}} - L_{\frac{i-1}{n_k}})
          - \int_0^s \sqrt[\alpha]{Y_{u-}} \, \dd L_u\Biggr|
   \to 0 \qquad \text{as}\qquad k \to \infty
 \end{align*}
 hold \ $\PP$-almost surely.
Let us denote by \ $\tA_t$ \ the event on which the above two  \ $\PP$-almost sure convergences hold.
Consequently, with the notation
 \[
   \ttA_t := \biggl\{\omega \in \Omega : \int_0^t Y_s(\omega) \, \dd s = 0\biggr\} ,
 \]
 we have
 \begin{align*}
  \ttA_t \cap \tA_t \cap A_t
  &\subset \ttA_t
           \cap \biggl\{\omega \in \Omega
                  : \text{$\biggl(\int_0^s \sqrt{Y_u} \, \dd W_u\biggr)(\omega) = 0$,
                          \ $\biggl(\int_0^s
                                     \sqrt[\alpha]{Y_{u-}} \, \dd L_u\biggr)(\omega) = 0$
                          \ for all \ $s \in [0, t)$}\biggr\} \\
  &\subset \ttA_t
           \cap \bigl\{\omega \in \Omega
                  : \text{$Y_s(\omega) = Y_0(\omega) + a s$ \ for all
                          \ $s \in [0, t)$}\bigr\}\\
  &\subset \ttA_t
           \cap \biggl\{\omega \in \Omega
                  : \text{$\int_0^s (Y_0(\omega) + a u) \, \dd u = 0$
                          \ for all \ $s \in [0, t)$}\biggr\} \\
  &\subset \ttA_t
           \cap \Bigl\{\omega \in \Omega
                  : \text{$Y_0(\omega) s + \frac{as^2}{2} = 0$
                          \ for all \ $s \in [0, t)$}\Bigr\} \\
  &\subset \ttA_t
           \cap \Bigl\{\omega \in \Omega
                  : \text{$Y_0(\omega) = - \frac{as}{2}$ \ for all \ $s \in [0, t)$}\Bigr\} ,
 \end{align*}
 where the last event has probability \ $0$, \ implying
 \ $\PP\bigl(\int_0^t Y_s(\omega) \, \dd s = 0\bigr) = 0$.
\ Thus  \ $\PP\bigl(\int_0^t Y_s(\omega) \, \dd s \in \RR_{++}\bigr) = 0$, \ and hence we have
 (v).

For (vi), see Proposition 3.7 in Jiao et al.\ \cite{JiaMaSco}.

Finally, we prove part (vii).
First note that in case of \ $a = 0$, \ $(Y_t)_{t\in\RR_+}$ \ is a continuous time branching
 process (without immigration).
If \ $b \in \RR_+$, \ then by Corollary 3.9 in Li \cite{Li},
 \ $\PP(\tau_0 < \infty \mid Y_0 = y_0) = 1$ \ for all \ $y_0 \in \RR_{++}$, \ since Condition
 3.6 in Li \cite{Li} holds for all \ $\theta > 0$ \ due to
 \ $\int_\theta^\infty \frac{1}{R(z)} \, \dd z
    \leq \int_\theta^\infty \frac{2}{\sigma^2z^2} \, \dd z < \infty$.
\ The last statement follows from the fact that in case of \ $a = 0$ \ and \ $\PP(Y_0 = 0) = 1$,
 \ the pathwise unique non-negative strong solution of the SDE \eqref{stable_CIR} is \ $Y_t = 0$
 \ for all \ $t \in \RR_+$.
\proofend

Note that, by Proposition \ref{Pro_stable_CIR}, the process \ $(Y_t)_{t\in\RR_+}$ \ is a
 semimartingale, see, e.g., Jacod and Shiryaev \cite[I.4.33]{JSh}.
Now we derive a so-called Grigelionis form for the semimartingale \ $(Y_t)_{t\in\RR_+}$, \ see,
 e.g., Jacod and Shiryaev \cite[III.2.23]{JSh} or Jacod and Protter
 \cite[Theorem 2.1.2]{JacPro}.

\begin{Pro}\label{Grigelionis}
Let \ $\eta_0$ \ be a random variable independent of \ $(W_t)_{t\in\RR_+}$ \ and
 \ $(L_t)_{t\in\RR_+}$ \ satisfying \ $\PP(\eta_0 \in \RR_+) = 1$ \ and
 \ $\EE(\eta_0) < \infty$.
\ For \ $a \in \RR_+$, \ $b \in \RR$, \ $\sigma\in\RR_+$, \ $\delta \in \RR_{++}$, \ and \ $\alpha\in(1,2)$,
 \ let \ $(Y_t)_{t\in\RR_+}$ \ be the unique strong solution of the SDE \eqref{stable_CIR} satisfying
 \ $\PP(Y_0 = \eta_0) = 1$.
\ Then the Grigelionis form of \ $(Y_t)_{t\in\RR_+}$ \ takes the form
 \begin{equation}\label{YGrigelionis}
  \begin{aligned}
   Y_t &= Y_0 + \int_0^t (a-bY_u + \gamma \delta \sqrt[\alpha]{Y_u}) \, \dd u
          + \int_0^t
             \left(\int_\RR
                    (h(z\delta\sqrt[\alpha]{Y_u}) - \delta \sqrt[\alpha]{Y_u} h(z))
                    \, m(\dd z)\right)
             \dd u \\
       &\quad
          + \sigma \int_0^t \sqrt{Y_u} \, \dd W_u \\
       &\quad
          + \int_0^t \int_\RR h( z\delta\sqrt[\alpha]{Y_{u-}}) \, \tmu^L(\dd u,\dd z)
          + \int_0^t
             \int_\RR
              (z\delta\sqrt[\alpha]{Y_{u-}}  - h(z\delta\sqrt[\alpha]{Y_{u-}}))
              \, \mu^L(\dd u,\dd z)
  \end{aligned}
 \end{equation}
 for \ $t \in \RR_+$, \ where \ $h : \RR \to [-1, 1]$, \ $h(z):= z \bbone_{[-1,1]}(z)$,
 \ $z \in \RR$.
\end{Pro}

\noindent{\bf Proof.}
Using \eqref{help_Levy_Ito_decomposition} and Proposition II.1.30 in Jacod and Shiryaev
 \cite{JSh}, we obtain
 \begin{align*}
  Y_t &= Y_0 + \int_0^t (a-bY_u) \, \dd u + \int_0^t \sigma \sqrt{Y_u} \, \dd W_u
         + \delta \int_0^t \sqrt[\alpha]{Y_{u-}} \, \dd L_u \\
      &= Y_0 + \int_0^t (a-bY_u) \, \dd u + \int_0^t \sigma \sqrt{Y_u} \, \dd W_u
         + \gamma \delta \int_0^t \sqrt[\alpha]{Y_{u-}} \, \dd u \\
      &\quad
         + \delta \int_0^t \int_\RR \sqrt[\alpha]{Y_{u-}} h(z) \, \tmu^L(\dd u,\dd z)
         + \delta \int_0^t \int_\RR \sqrt[\alpha]{Y_{u-}} (z - h(z)) \, \mu^L(\dd u,\dd z)
 \end{align*}
 for \ $t \in \RR_+$.
\ In order to prove the statement, it is enough to show
 \begin{gather}\label{G1}
  \delta
  \int_0^t
   \int_\RR \sqrt[\alpha]{Y_{u-}} h(z) \, \big(\mu^L(\dd u,\dd z) - \dd u \, m(\dd z)\big)
  = I_1 - I_2 , \\
  \delta \int_0^t \int_\RR \sqrt[\alpha]{Y_{u-}} (z - h(z)) \, \mu^L(\dd u,\dd z)
  = I_3 + I_4 , \label{G2}
 \end{gather}
 with
 \begin{align*}
  I_1 &:= \int_0^t \int_\RR
           h(z\delta\sqrt[\alpha]{Y_{u-}} ) \, \bigl(\mu^L(\dd u,\dd z) - \dd u
           \, m(\dd z)\bigr) , \\
  I_2 &:= \int_0^t \int_\RR
           (h(z\delta\sqrt[\alpha]{Y_{u-}}) - \delta\sqrt[\alpha]{Y_{u-}} h(z))
           \, \bigl(\mu^L(\dd u,\dd z) - \dd u \, m(\dd z)\bigr) , \\
  I_3 &:= \int_0^t \int_\RR
           (z\delta\sqrt[\alpha]{Y_{u-}}  - h(z\delta\sqrt[\alpha]{Y_{u-}}))
           \, \mu^L(\dd u,\dd z) , \\
  I_4 &:= \int_0^t \int_\RR
           (h(z \delta\sqrt[\alpha]{Y_{u-}} )- \delta\sqrt[\alpha]{Y_{u-}}h(z))
           \, \mu^L(\dd u,\dd z) ,
 \end{align*}
 and the equality
 \begin{equation}\label{drift}
  I_4 - I_2 = I_5 \qquad \text{with} \qquad
  I_5 := \int_0^t
          \left(\int_\RR
                 (h(z\delta\sqrt[\alpha]{Y_u}) - \delta\sqrt[\alpha]{Y_u}h(z))
                 \, m(\dd z)\right) \dd u .
 \end{equation}
For the equations \eqref{G1}, \eqref{G2} and \eqref{drift}, it suffices to check the
 existence of \ $I_2$, \ $I_3$ \ and \ $I_5$.

First note that for every \ $s \in (0, \infty)$ \ we have
 \begin{align}\label{help_ITO_atiras}
  h(s z)- sh(z)
  = \begin{cases}
     s z \bbone_{\{1<|z|\leq\frac{1}{s}\}}
      & \text{if \ $s \in (0, 1)$, \ $z \in \RR$,} \\
     0 & \text{if \ $s = 1$, \ $z \in \RR$,} \\
     -s z \bbone_{\{\frac{1}{s}<|z|\leq1\}}
      & \text{if \ $s \in (1, \infty)$, \ $z \in \RR$.}
    \end{cases}
 \end{align}
The existence of \ $I_2$ \ will be a consequence of \ $I_2 = I_{2,1} - I_{2,2} - I_{2,3}$ \ with
 \begin{align*}
  I_{2,1} &:= \int_0^t \int_\RR
               \delta \sqrt[\alpha]{Y_{u-}}
               z \bbone_{\{1<|z|\leq\frac{1}{\delta\sqrt[\alpha]{Y_{u-}}}\}}
               \bbone_{\{\delta\sqrt[\alpha]{Y_{u-}}\in(0,1)\}}
               \, \mu^L(\dd u,\dd z) , \\
  I_{2,2} &:= \int_0^t \int_\RR
               \delta \sqrt[\alpha]{Y_{u-}}
               z \bbone_{\{1<|z|\leq\frac{1}{\delta\sqrt[\alpha]{Y_{u-}}}\}}
               \bbone_{\{\delta\sqrt[\alpha]{Y_{u-}}\in(0,1)\}}
               \, \dd u \, m(\dd z) , \\
  I_{2,3} &:= \int_0^t \int_\RR
               \delta \sqrt[\alpha]{Y_{u-}}
               z \bbone_{\{\frac{1}{\delta\sqrt[\alpha]{Y_{u-}}}<|z|\leq1\}}
               \bbone_{\{\delta\sqrt[\alpha]{Y_{u-}}\in(1,\infty)\}}
               \, \bigl(\mu^L(\dd u,\dd z) - \dd u \, m(\dd z)\bigr) ,
 \end{align*}
 since on the set \ $\{Y_{u-} = 0\}$, \ the integrand
 \ $h(z\delta\sqrt[\alpha]{Y_{u-}}) - \delta\sqrt[\alpha]{Y_{u-}} h(z)$ \ in \ $I_2$ \ takes
 value \ $0$.
\ Here we have
 \[
   |I_{2,1}|
   \leq
   \int_0^t \int_\RR
    |\delta \sqrt[\alpha]{Y_{u-}} z| \bbone_{\{1<|z|\leq\frac{1}{\delta\sqrt[\alpha]{Y_{u-}}}\}}
    \bbone_{\{\delta\sqrt[\alpha]{Y_{u-}}\in(0,1)\}}
    \, \mu^L(\dd u,\dd z)
   \leq \int_0^t \int_\RR \bbone_{\{1<|z|\}} \, \mu^L(\dd u,\dd z)
   < \infty
 \]
  \ $\PP$-almost surely, see, e.g., Sato \cite[Lemma 20.1]{Sat}.
Moreover,
 \begin{align*}
  |I_{2,2}|
  &\leq
   \int_0^t \int_\RR
    |\delta\sqrt[\alpha]{Y_{u-}} z| \bbone_{\{1<|z|\leq\frac{1}{\delta\sqrt[\alpha]{Y_{u-}}}\}}
    \bbone_{\{\delta\sqrt[\alpha]{Y_{u-}}\in(0,1)\}}
    \, \dd u \, m(\dd z) \\
  &\leq
   \int_0^t \int_\RR \bbone_{\{1<|z|\}} \, \dd u \, m(\dd z)
   = t m(\{z \in \RR : |z| > 1 \})
   < \infty .
 \end{align*}
Further, the function \ $\Omega \times \RR_+ \times \RR \ni (\omega, t, z) \mapsto h(z)$
 \ belongs to \ $G_{\mathrm{loc}}(\mu^L)$, \ see Jacod and Shiryaev
 \cite[Definitions II.1.27, Theorem II.2.34]{JSh}.
We have
 \ $|z \bbone_{\{\frac{1}{\delta\sqrt[\alpha]{Y_{u-}}}<|z|\leq1\}}
     \bbone_{\{\delta\sqrt[\alpha]{Y_{u-}}\in(1,\infty)\}}|
    \leq |h(z)|$,
 \ hence, by the definition of \ $G_{\mathrm{loc}}(\mu^L)$, \ the function
 \ $\Omega \times \RR_+ \times \RR \ni (\omega, t, z)
    \mapsto z \bbone_{\{\frac{1}{\delta\sqrt[\alpha]{Y_{u-}}}<|z|\leq1\}}
    \bbone_{\{\delta\sqrt[\alpha]{Y_{u-}}\in(1,\infty)\}}$
 \ also belongs to \ $G_{\mathrm{loc}}(\mu^L)$.
\ By Jacod and Shiryaev \cite[Proposition II.1.30]{JSh}, we conclude that the function
 \ $\Omega \times \RR_+ \times \RR \ni (\omega, t, z)
    \mapsto \delta \sqrt[\alpha]{Y_{u-}}
            z \bbone_{\{\frac{1}{\delta\sqrt[\alpha]{Y_{u-}}}<|z|\leq1\}}
            \bbone_{\{\delta\sqrt[\alpha]{Y_{u-}}\in(1,\infty)\}}$
 \ also belongs to \ $G_{\mathrm{loc}}(\mu^L)$, \ thus the integral \ $I_{2,3}$ \ exists, and
 hence we obtain the existence of \ $I_2$, \ and hence that of \ $I_1$.

Next observe that for the process
 \ $\zeta_t := \delta \int_0^t \sqrt[\alpha]{Y_{u-}} \, \dd L_u$, \ $t \in \RR_+$, \  we have
 \ $\Delta \zeta_t = \delta \sqrt[\alpha]{Y_{t-}} \Delta L_t$, \ $t \in \RR_+$, \ following
 from \eqref{help_Levy_Ito_decomposition} and Jacod and Shiryaev
 \cite[Definitions II.1.27]{JSh}.
Consequently,
 \begin{align*}
  I_3
  &= \int_0^t \int_\RR
      z \delta \sqrt[\alpha]{Y_{u-}} \bbone_{\{|z\delta\sqrt[\alpha]{Y_{u-}}|>1\}}
      \, \mu^L(\dd u,\dd z)
   = \sum_{u\in[0,t]}
      \Delta L_u \delta \sqrt[\alpha]{Y_{u-}}
      \bbone_{\{|\Delta L_u \delta\sqrt[\alpha]{Y_{u-}}|>1\}} \\
  &= \sum_{u\in[0,t]} \Delta \zeta_u \bbone_{\{|\Delta \zeta_u|>1\}}
 \end{align*}
 is a finite sum, since the process \ $(\zeta_t)_{t\in[0,\infty)}$ \ admits c\`adl\`ag
 trajectories, hence there can be at most finitely many points \ $u \in [0, t]$ \ at which the
 absolute value \ $|\Delta \zeta_u|$ \ of the jump size \ $\Delta \zeta_u$ \ exceeds 1, see,
 e.g., Billingsley \cite[page 122]{Bil}.
Thus we obtain the existence of \ $I_3$, \ and hence that of \ $I_4$.

Finally, we have
 \begin{align*}
  |I_5|
  &\leq
   \int_0^t
    \left(\int_\RR
           |\delta\sqrt[\alpha]{Y_u} z| \bbone_{\{1<|z|\leq\frac{1}{\delta\sqrt[\alpha]{Y_u}}\}}
            \bbone_{\{\delta\sqrt[\alpha]{Y_u}\in(0,1)\}}
           \, m(\dd z)\right) \dd u \\
  &\phantom{\leq}
    + \int_0^t
       \left(\int_\RR
              |\delta\sqrt[\alpha]{Y_u} z|
              \bbone_{\{\frac{1}{\delta\sqrt[\alpha]{Y_u}}<|z|\leq1\}}
              \bbone_{\{\delta \sqrt[\alpha]{Y_u}\in(1,\infty)\}}
              \, m(\dd z)\right) \dd u \\
  &\leq
   \int_0^t
    \left(\int_\RR \bbone_{\{1<|z|\}} \, m(\dd z)\right) \dd u
   + \int_0^t
      \left(\int_\RR |\delta\sqrt[\alpha]{Y_u} z|^2 \bbone_{\{|z|\leq1\}} \, m(\dd z)\right)
      \dd u \\
  &= t m(\{z \in \RR : |z| > 1\})
     + \int_0^t \delta^2 Y_u^{\frac{2}{\alpha}} \, \dd u \int_{-1}^1 |z|^2 \, m(\dd z)
   < \infty ,
 \end{align*}
 since
 \ $\int_{-1}^1 |z|^2 \, m(\dd z) = \int_0^1 z^2 C_\alpha z^{-1-\alpha} \, \dd z
    = \frac{C_\alpha}{2-\alpha} \in \RR_{++}$,
 \ hence we conclude the existence of \ $I_5$.
\proofend

Next we present a result about the first moment of \ $(Y_t)_{t\in\RR_+}$.

\begin{Pro}\label{Pro_moments}
Let \ $a \in \RR_+$, \ $b \in \RR$, \ $\sigma\in\RR_+$, \ $\delta \in \RR_{++}$, \ and \ $\alpha\in(1,2)$.
\ Let \ $(Y_t)_{t\in\RR_+}$ \ be the unique strong solution of the SDE \eqref{stable_CIR}
 satisfying \ $\PP(Y_0 \in \RR_+) = 1$ \ and \ $\EE(Y_0) < \infty$.
\ Then
 \begin{align}\label{expectation}
  \EE(Y_t)
  = \begin{cases}
     \ee^{-bt} \bigl(\EE(Y_0) - \frac{a}{b}\bigr) + \frac{a}{b}
      & \text{if \ $b \ne 0$,} \\[1mm]
     \EE(Y_0) + a t
      & \text{if \ $b = 0$,}
    \end{cases}
  \qquad t \in \RR_+ .
 \end{align}
Consequently, if \ $b \in \RR_{++}$, \ then
 \begin{equation}\label{expectation_sub}
  \lim_{t\to\infty} \EE(Y_t) = \frac{a}{b} ,
 \end{equation}
 if \ $b = 0$, \ then
 \[
   \lim_{t\to\infty} t^{-1} \EE(Y_t) = a ,
 \]
 if \ $b \in \RR_{--}$, \ then
 \[
   \lim_{t\to\infty} \ee^{bt} \EE(Y_t)
   = \EE(Y_0) - \frac{a}{b} .
 \]
\end{Pro}

\noindent{\bf Proof.}
By Proposition \ref{Pro_stable_CIR}, \ $(Y_t)_{t\in\RR_+}$ \ is CBI process with an
 infinitesimal generator given in \eqref{infgen}.
By the notations of Barczy et al.\ \cite{BarLiPap}, this CBI process has parameters
 \ $(d, c, \beta, B, \nu, \mu)$, \ where \ $d = 1$, \ $c = \frac{1}{2} \sigma^2$,
 \ $\beta = a$, \ $B = -b - \int_0^\infty (z-1)^+ \, \mu(\dd z)$, \ $\nu = 0$ \ and
 \ $\mu = \delta^\alpha m$.
\ Since \ $\EE(Y_0) < \infty$ \ and the moment condition
 \ $\int_{\RR\setminus\{0\}} |z| \bbone_{\{|z|\geq1\}} \, \nu(\dd z) < \infty$
 \ trivially holds, we may apply formula (3.1.11) in Li \cite{Li2} or Lemma 3.4 and (2.14) in
 Barczy et al.\ \cite{BarLiPap} with \ $\tB = B +  \int_0^\infty (z-1)^+ \, \mu(\dd z) = - b$
 \ and \ $\tbeta = \beta + \int_{\RR\setminus\{0\}} z \, \nu(\dd z) = a$ \ yielding that
 \[
   \EE(Y_t) = \ee^{t\tB} \EE(Y_0) + \left(\int_0^t \ee^{u\tB} \, \dd u\right) \tbeta.
 \]
This implies \eqref{expectation} and the other parts of the assertion.
\proofend

Based on the asymptotic behavior of the expectations \ $\EE(Y_t)$ \ as \ $t \to \infty$, \ we
 introduce a classification of the stable CIR model given by the SDE \eqref{stable_CIR}.

\begin{Def}\label{Def_criticality}
Let \ $a \in \RR_+$, \ $b \in \RR$, \ $\sigma\in\RR_+$, \ $\delta \in \RR_{++}$, \ and \ $\alpha\in(1,2)$.
\ Let \ $(Y_t)_{t\in\RR_+}$ \ be the unique strong solution of the SDE \eqref{stable_CIR}
 satisfying \ $\PP(Y_0 \in \RR_+) = 1$ \ and \ $\EE(Y_0) < \infty$.
\ We call \ $(Y_t)_{t\in\RR_+}$ \ subcritical, critical or supercritical if \ $b \in \RR_{++}$,
 \ $b = 0$ \ or \ $b \in \RR_{--}$, \ respectively.
\end{Def}

The following result states the existence of a unique stationary distribution for the process
 \ $(Y_t)_{t\in\RR_+}$ \ in the subcritical and critical cases, and the exponential ergodicity
 in the subcritical case.

\begin{Thm}\label{Ergodicity}
Let \ $a \in \RR_+$, \ $b \in \RR_+$, \ $\sigma\in\RR_+$, \ $\delta \in \RR_{++}$, \ and \ $\alpha\in(1,2)$.
\ Let \ $(Y_t)_{t\in\RR_+}$ \ be the unique strong solution of the SDE \eqref{stable_CIR}
 satisfying \ $\PP(Y_0 \in \RR_+) = 1$ \ and \ $\EE(Y_0) < \infty$.
 \renewcommand{\labelenumi}{{\rm(\roman{enumi})}}
 \begin{enumerate}
  \item
   Then \ $(Y_t)_{t\in\RR_+}$ \ converges in law to its unique stationary distribution \ $\pi$
    \ having Laplace transform
    \begin{equation}\label{help_limit_MLE_critical2_spec}
      \int_0^\infty \ee^{-\lambda y} \, \pi(\dd y)
      = \exp\biggr\{- \int_0^\lambda \frac{F(x)}{R(x)} \, \dd x\biggl\}
      = \exp\biggr\{- \int_0^\lambda
                       \frac{ax}
                            {\frac{\sigma^2}{2}x^2 + \frac{\delta^\alpha}{\alpha}x^\alpha +bx}
                       \, \dd x\biggl\}
    \end{equation}
    for \ $\lambda \in \RR_+$.
   \ Especially, in case of \ $b = 0$ \ and \ $\sigma = 0$, \
   \ $\pi$ \ is a strictly \ $(2-\alpha)$-stable distribution with no negative jumps.
   Moreover, the expectation of \ $\pi$ \ is given by
    \begin{align}\label{help_stac_expectation}
     \int_0^\infty y \, \pi(\dd y)
     = \begin{cases}
        0 & \text{if \ $a = 0$ \ and \ $b = 0$,} \\
        \frac{a}{b} \in \RR_+ & \text{if \ $b \in \RR_{++}$,} \\
        +\infty & \text{if \ $a \in \RR_{++}$ \ and \ $b = 0$.}
       \end{cases}
    \end{align}
  \item
   If, in addition, \ $a \in \RR_{++}$ \ and \ $b \in \RR_{++}$, \ then the process
    \ $(Y_t)_{t\in\RR_+}$ \ is \ exponentially ergodic, i.e., there exist constants
    \ $C \in \RR_{++}$ \ and \ $D \in \RR_{++}$ \ such that
    \[
      \|\PP_{Y_t \mid Y_0 = y} - \pi\|_{\mathrm{TV}}
      \leq C (y + 1) \ee^{- D t} , \qquad t \in \RR_+ , \qquad y \in \RR_+ ,
    \]
    where \ $\|\mu\|_{\mathrm{TV}}$ \ denotes the total-variation norm of a signed measure
    \ $\mu$ \ on \ $\RR_+$ \ defined by
    \ $\|\mu\|_{\mathrm{TV}} := \sup_{A\in\cB(\RR_+)} |\mu(A)|$,
    \ and \ $\PP_{Y_t \mid Y_0 = y}$ \ is the conditional distribution of \ $Y_t$ \ with respect
    to the condition \ $Y_0 = y$.
   \ As a consequence, for all Borel measurable functions \ $f : \RR_+ \to \RR$ \ with
    \ $\int_0^\infty |f(y)| \, \pi(\dd y) < \infty$, \ we have
    \begin{equation}\label{ergodic}
     \frac{1}{T} \int_0^T f(Y_s) \, \dd s \as \int_0^\infty f(y) \, \pi(\dd y) \qquad
     \text{as \ $T \to \infty$.}
    \end{equation}
 \end{enumerate}
\end{Thm}

\noindent{\bf Proof.}
The weak convergence of \ $Y_t$ \ towards \ $\pi$ \ as \ $t \to \infty$, \ and the
 fact that \ $\pi$ \ is a stationary distribution for \ $(Y_t)_{t\in\RR_+}$ \ follow
 immediately from Li \cite[Theorem 3.20 and the paragraph after Corollary 3.21]{Li},
 since
 \ $R(z) = \frac{\sigma^2}{2} z^2 + \frac{\delta^\alpha}{\alpha} z^\alpha + b z
         \in \RR_{++}$,
 \ $z \in \RR_{++}$, \ and condition (3.30) in Li \cite{Li} is satisfied.
Indeed, for all \ $\lambda \in \RR_{++}$,
 \begin{align*}
  \int_0^\lambda \frac{F(z)}{R(z)} \, \dd z
  = \int_0^\lambda
     \frac{az}{\frac{\sigma^2}{2}z^2+\frac{\delta^\alpha}{\alpha}z^\alpha+bz} \, \dd z
  \leq \frac{a\alpha}{\delta^\alpha} \int_0^\lambda z^{1-\alpha} \, \dd z
  = \frac{a\alpha\lambda^{2-\alpha}}{\delta^\alpha(2-\alpha)}
  < \infty .
 \end{align*}
We note that Li and Ma \cite[Proposition 2.2]{LiMa} contains the above considerations in case of
 \ $b \in \RR_{++}$.
\ The uniqueness of a stationary distribution in (i) follows from, e.g., page 80 in Keller-Ressel
 \cite{KR}.
Namely, let us assume that there exists another stationary distribution \ $\pi'$ \ for
 \ $(Y_t)_{t\in\RR_+}$, \ and let \ $(Y_t')_{t\in\RR_+}$ \ be the unique strong solution of the
 SDE \eqref{stable_CIR} with \ $a \in \RR_+$, \ $b \in \RR_+$, \ $\sigma \in \RR_+$, \ and
 \ $\delta \in \RR_{++}$ \ satisfying \ $\cL(Y_0') = \pi'$, \ where \ $\cL(Y_0')$ \ denotes the
 law of \ $Y_0'$.
\ Then, by part (iii) of Proposition \ref{Pro_stable_CIR}, for all \ $\lambda \in \RR_+$,
 \begin{align*}
  \lim_{t\to\infty} \EE(\ee^{-\lambda Y_t'})
  &= \lim_{t\to\infty} \EE(\EE(\ee^{-\lambda Y_t'} \mid Y_0'))
   = \lim_{t\to\infty}
      \EE\left(\exp\left\{- Y_0' v_t(\lambda)
                          + \int_\lambda^{v_t(\lambda)}
                             \frac{F(z)}{R(z)} \, \dd z\right\}\right) \\
  &= \EE\left(\exp\left\{- \int_0^\lambda \frac{F(z)}{R(z)} \, \dd z\right\}\right)
   = \exp\left\{- \int_0^\lambda \frac{F(z)}{R(z)} \, \dd z \right\} ,
 \end{align*}
 where the last but one step follows by the dominated convergence theorem and the fact that
 \ $\lim_{t\to\infty} v_t(\lambda) = 0$, \ $\lambda \in \RR_+$ \ (see, e.g., the proof of
 Theorem 3.20 in Li \cite{Li}).
Since \ $\cL(Y_t') = \pi'$, \ $t \in \RR_+$, \ we have
 \ $\int_0^\infty \ee^{-\lambda z} \, \pi'(\dd z)
    = \exp\big\{- \int_0^\lambda \frac{F(z)}{R(z)} \, \dd z\big\}$,
 \ $\lambda \in \RR_+$, \ yielding that
 \ $\int_0^\infty \ee^{-\lambda z} \, \pi'(\dd z)
    = \int_0^\infty \ee^{-\lambda z} \, \pi(\dd z)$, \ $\lambda \in \RR_+$.
\ By the uniqueness of Laplace transform, we get \ $\pi = \pi'$, \ as desired.

Further, in case of \ $b = 0$ \ and \ $\sigma = 0$, \ by
 \eqref{help_limit_MLE_critical2_spec}, we have
 \begin{align*}
  \int_0^\infty \ee^{-\lambda y} \, \pi(\dd y)
  = \exp\left\{- \int_0^\lambda
                \frac{ax}{\frac{\delta^\alpha}{\alpha}x^\alpha} \, \dd x\right\}
  = \exp\left\{- \frac{a\alpha}{(2-\alpha)\delta^\alpha} \lambda^{2-\alpha} \right\} ,
  \qquad \lambda \in \RR_+ ,
 \end{align*}
 so \ $\pi$ \ is a strictly \ $(2-\alpha)$-stable distribution with no negative jumps.
Finally, again by \eqref{help_limit_MLE_critical2_spec},
 \begin{align*}
  \int_0^\infty y \, \pi(\dd y)
  = \lim_{\lambda\downarrow0}
     \frac{a\lambda}
          {\frac{\sigma^2}{2}\lambda^2+\frac{\delta^\alpha}{\alpha}\lambda^\alpha+b\lambda}
  = \begin{cases}
        0 & \text{if \ $a = 0$ \ and \ $b = 0$,} \\
        \frac{a}{b} \in \RR_+ & \text{if \ $b \in \RR_{++}$,} \\
        +\infty & \text{if \ $a \in \RR_{++}$ \ and \ $b = 0$.}
    \end{cases}
 \end{align*}
For part (ii), we can use Theorem 2.5 in Li and Ma \cite{LiMa}.
We only have to check Condition 2.1 of Li and Ma \cite{LiMa}, namely, we have to show the
 existence of some constant \ $\theta \in \RR_{++}$ \ such that \ $R(z) \in \RR_{++}$ \ for
 \ $z > \theta$ \ and \ $\int_\theta^{\infty} \frac{1}{R(z)} \, \dd z < \infty$.
\ Here \ $R(z) \in \RR_{++}$ \ for all \ $z \in \RR_{++}$ \ (due to \ $b \in \RR_{++}$), \ and,
 e.g., with \ $\theta = 1$,
 \begin{align*}
  \int_1^{\infty} \frac{1}{R(z)} \, \dd z
  \leq \frac{\alpha}{\delta^\alpha} \int_1^{\infty} \frac{1}{z^\alpha} \, \dd z
  = \frac{\alpha}{(\alpha-1)\delta^\alpha}
  < \infty .
 \end{align*}
In case of \ $\sigma = 0$, \ the exponential ergodicity of \ $(Y_t)_{t\in\RR_+}$ \ also follows
 by Theorem 6.1 in Jin et al.\ \cite{JinKreRud}.
Convergence \eqref{ergodic} follows, e.g., from the discussion after Proposition 2.5 in
 Bhattacharya \cite{Bha}.
\proofend

\begin{Rem}\label{Rem_ergodic_alternative_proof}
In what follows, in case of \ $\sigma \in \RR_{++}$, \ we present another (and more detailed)
 proof for the convergence \ $Y_t \distr \pi$ \ as \ $t \to \infty$ \ in Theorem
 \ref{Ergodicity} giving more insight as well.
It is enough to consider the case of \ $\PP(Y_0 = y_0) = 1$ \ with some \ $y_0 \in \RR_+$.
Using \eqref{Laplace_Y_t1}, we have
 \[
   \EE(\ee^{-\lambda Y_t})
   = \exp\bigg\{- y_0 v_t(\lambda) - a \int_0^t v_s(\lambda) \dd s\bigg\}
 \]
 for \ $t \in \RR_+$ \ and \ $\lambda \in \RR_+$, \ where the function
 \ $\RR_+ \ni \lambda \mapsto v_t(\lambda) \in \RR_+$ \ is given by \eqref{LiLaplace}.
First, we show that \ $\lim_{t\to\infty} v_t(\lambda) = 0$.
\ The proof is based on the following version of the comparison theorem (see, e.g., Lemma C.3.
 in Filipovi\'{c} et al.\ \cite{FilMaySch} or Amann \cite[Lemma 16.4]{Ama}): if
 \ $S : \RR_+ \times \RR \to \RR$ \ is a continuous function which is locally Lipschitz
 continuous in its second variable and \ $p, q : \RR_+ \to \RR$ \ are differentiable functions
 satisfying
 \begin{align*}
  & p'(s) \leq S(s, p(s)), \qquad s \in \RR_+ , \\
  & q'(s) = S(s, q(s)), \qquad s \in \RR_+ , \\
  & p(0) \leq q(0) ,
 \end{align*}
 then \ $p(s) \leq q(s)$ \ for all \ $s \in \RR_+$.
\ By choosing \ $S : \RR_+ \times \RR \to \RR$, \ $S(s, x) := - \frac{\sigma^2}{2} x^2$,
 \ $(s, x) \in \RR_+ \times \RR$, \ the comparison theorem yields that
 \begin{align}\label{help_crit1}
  0 \leq v_s(\lambda) \leq f(s), \qquad s \in \RR_+ ,
 \end{align}
 where \ $f : \RR_+ \to \RR_+$ \ is the unique locally bounded solution to the differential
 equation
 \[
   f'(s) = S(s,f(s)) = - \frac{\sigma^2}{2} f(s)^2 ,
   \qquad s \in \RR_+ \qquad \text{with} \qquad f(0) =  \lambda .
 \]
The solution of this separable differential equation takes the form
 \begin{align}\label{psi_form}
  f(s) = \frac{\lambda}{1+\frac{\sigma^2\lambda}{2}s}, \qquad s \in \RR_+ .
 \end{align}
Hence, using \ $\sigma \in \RR_{++}$, \ we readily have \ $\lim_{t\to\infty} f(t) = 0$,
 \ which, by \eqref{help_crit1}, yields that \ $\lim_{t\to\infty} v_t(\lambda) = 0$ \ for all \ $\lambda\in\RR_+$,
 \ as desired.
Further, by \eqref{Laplace_Y_t2}, we have
 \[
   \lim_{t\to\infty} \int_0^t v_s(\lambda) \, \dd s
   = - \lim_{t\to\infty}
        \int_\lambda^{v_t(\lambda)}
         \frac{z}{\frac{\sigma^2}{2}z^2+\frac{\delta^\alpha}{\alpha}z^\alpha+bz} \, \dd z ,
         \qquad \lambda\in\RR_+.
 \]
Then, by the continuity of the integral upper limit function, we have
 \[
   \lim_{t\to\infty} \int_0^t v_s(\lambda) \, \dd s
   = \int_0^\lambda
      \frac{z}{\frac{\sigma^2}{2}z^2+\frac{\delta^\alpha}{\alpha}z^\alpha+bz} \, \dd z ,
   \qquad \lambda \in \RR_+ ,
 \]
 where the integral on the right hand side is well-defined, since
 \begin{align*}
  \int_0^\lambda
   \frac{z}{\frac{\sigma^2}{2}z^2+\frac{\delta^\alpha}{\alpha}z^\alpha+bz} \, \dd z
  \leq \int_0^\lambda \frac{z}{\frac{\delta^\alpha}{\alpha}z^\alpha} \, \dd z
  = \frac{\alpha}{\delta^\alpha} \int_0^\lambda \frac{1}{z^{\alpha-1}} \, \dd z
  = \frac{\alpha\lambda^{2-\alpha}}{\delta^\alpha(2-\alpha)}
  < \infty .
 \end{align*}
Consequently, by continuity theorem, we have \eqref{help_limit_MLE_critical2_spec} in case of
 \ $\sigma \in \RR_{++}$.
\proofend
\end{Rem}

Next we give a statistic for \ $\sigma^2$ \ using continuous time observations
 \ $(Y_t)_{t\in[0,T]}$ \ with an arbitrary \ $T \in \RR_{++}$.
\ Due to this result we do not consider the estimation of the parameter \ $\sigma$, \ it is
 supposed to be known.
The parameter \ $\sigma$ \ is a parameter for the diffusion part related to \ $W$, \ and, in general,
 the estimation of this kind of parameter is possible using an arbitrarily short (continuous time) observation
 of the underlying process (at least theoretically), and that's why authors suppose this to be known.
In the forthcoming Remark \ref{Thm_MLE_cons_sigma}, we demonstrate that this also holds for our model.

\begin{Rem}\label{Thm_MLE_cons_sigma}
Let \ $a \in \RR_+$, \ $b \in \RR$, \ $\sigma\in\RR_+$, \ $\delta \in \RR_{++}$, \ and \ $\alpha\in(1,2)$.
\ Let \ $(Y_t)_{t\in\RR_+}$ \ be the unique strong solution of the SDE \eqref{stable_CIR}
 satisfying \ $\PP(Y_0 \in \RR_+) = 1$ \ and \ $\EE(Y_0) < \infty$.
\ The Grigelionis representation given in \eqref{YGrigelionis} implies that the continuous
 martingale part \ $Y^\cont$ \ of \ $Y$ \ is
 \ $Y^\cont_t = \sigma \int_0^t \sqrt{Y_u} \, \dd W_u$, \ $t \in \RR_+$,
 \ see Jacod and Shiryaev \cite[III.2.28 Remarks, part 1)]{JSh}.
Consequently, the (predictable) quadratic variation process of \ $Y^\cont$ \ is
 \ $\langle Y^\cont \rangle_t = \sigma^2 \int_0^t Y_u \, \dd u$, \ $t \in \RR_+$.
\ Suppose that we have \ $\PP(Y_0 \in \RR_{++}) = 1$ \ or \ $a \in \RR_{++}$.
\ Then for all \ $T \in \RR_{++}$, \ we have
 \[
   \sigma^2 = \frac{\langle Y^\cont \rangle_T}{\int_0^T Y_u \, \dd u} =: \hsigma^2_T ,
 \]
 since, due to (v) of Proposition \ref{Pro_stable_CIR},
 \ $\PP\bigl(\int_0^T Y_u \, \dd u \in \RR_{++}\bigr) = 1$.
\ We note that \ $\hsigma^2_T$ \ is a statistic, i.e., there exists a measurable function
 \ $\Xi : D([0,T], \RR) \to \RR$ \ such that \ $\hsigma^2_T = \Xi((Y_u)_{u\in[0,T]})$, \ where
 \ $D([0,T], \RR)$ \ denotes the space of real-valued c\`adl\`ag functions defined on \ $[0,T]$,
 \ since
 \begin{equation}\label{sigma}
  \frac{1}{\frac{1}{n} \sum_{i=1}^{\lfloor nT\rfloor} Y_{\frac{i-1}{n}}}
  \Biggl(\sum_{i=1}^{\lfloor nT\rfloor} \bigl(Y_{\frac{i}{n}} - Y_{\frac{i-1}{n}}\bigr)^2
         - \sum_{u\in[0,T]} (\Delta Y_u)^2\Biggr)
  \stoch \hsigma^2_T \qquad \text{as \ $n \to \infty$,}
 \end{equation}
 where the convergence in \eqref{sigma} holds \ $\PP$-almost surely along a suitable subsequence, the
 members of the sequence in \eqref{sigma} are measurable functions of \ $(Y_u)_{u\in[0,T]}$,
 \ and one can use Theorems 4.2.2 and 4.2.8 in Dudley \cite{Dud}.
Next we prove \eqref{sigma}.
By Theorem I.4.47 a) in Jacod and Shiryaev \cite{JSh},
 \[
   \sum_{i=1}^{\lfloor nT\rfloor} \bigl(Y_{\frac{i}{n}} - Y_{\frac{i-1}{n}}\bigr)^2
   \stoch [Y]_T \qquad \text{as \ $n \to \infty$,} \qquad
   T \in  \RR_+,
 \]
 where \ $([Y]_t)_{t\in \RR_+}$ \ denotes the quadratic variation process of the semimartingale
 \ $Y$.
\ By Theorem I.4.52 in Jacod and Shiryaev \cite{JSh},
 \[
   [Y]_T = \langle Y^\cont \rangle_T + \sum_{u\in[0,T]} (\Delta Y_u)^2 , \qquad T \in \RR_+ .
 \]
Consequently, for all \ $T \in \RR_+$, \ we have
 \[
   \sum_{i=1}^{\lfloor nT\rfloor} \bigl(Y_{\frac{i}{n}} - Y_{\frac{i-1}{n}}\bigr)^2
   - \sum_{u\in[0,T]} (\Delta Y_u)^2
   \stoch \langle Y^\cont \rangle_T \qquad \text{as \ $n \to \infty$.}
 \]
Moreover, for all \ $T \in \RR_+$, \ we have
 \[
   \frac{1}{n} \sum_{i=1}^{\lfloor nT\rfloor} Y_{\frac{i-1}{n}}
   \stoch \int_0^T Y_u \, \dd u  \qquad \text{as \ $n \to \infty$,}
 \]
 see Proposition I.4.44 in Jacod and Shiryaev \cite{JSh}.
Hence \eqref{sigma} follows by the fact that convergence in probability is closed under
 multiplication.
Finally, we note that \ $T$ \ is fixed above, and it is enough to know
 any observation \ $(Y_t)_{t\in[0,T]}$ \ to carry out the above calculations, where \ $T>0$ \ can be arbitrarily small.
\proofend
\end{Rem}

\section{Joint Laplace transform of \ $Y_t$ \ and \ $\int_0^t Y_s \, \dd s$}
\label{section_Laplace_crit}

Using Theorem 4.10 in Keller-Ressel \cite{KR2} we derive a formula for the joint Laplace
 transform of \ $Y_t$ \ and \ $\int_0^t Y_s \, \dd s$, \ where \ $t \in \RR_+$.
\ We note that this form of the joint Laplace transform in question is a consequence of Theorem
 5.3 in Filipovi\'{c} \cite{Fil} and a special case of Proposition 3.3 in
 Jiao et al.\ \cite{JiaMaSco} as well.

\begin{Thm}\label{Thm_Laplace}
Let \ $a \in \RR_+$, \ $b \in \RR$, \ $\sigma\in\RR_+$, \ $\delta \in \RR_{++}$, \ and \ $\alpha\in(1,2)$.
\ Let \ $(Y_t)_{t\in\RR_+}$ \ be the unique strong solution of the SDE \eqref{stable_CIR}
 satisfying \ $\PP(Y_0 = y_0) = 1$ \ with some \ $y_0 \in \RR_+$.
\ Then for all \ $u, v \in \RR_-$,
 \[
   \EE\left[\exp\left\{uY_t + v \int_0^t Y_s \, \dd s\right\}\right]
   = \exp\bigg\{y_0 \psi_{u,v}(t) + a \int_0^t \psi_{u,v}(s) \, \dd s \bigg\} ,
   \qquad t \in \RR_+ ,
 \]
 where the function \ $\psi_{u,v} : \RR_+ \to \RR_-$ \ is the unique locally bounded solution to
 the differential equation
 \begin{align}\label{psi_DE}
  \psi'_{u,v}(t)
  = \frac{\sigma^2}{2} \psi_{u,v}(t)^2 + \frac{\delta^\alpha}{\alpha} (-\psi_{u,v}(t))^\alpha
    - b \psi_{u,v}(t) + v ,
  \qquad t \in \RR_+ , \qquad
  \psi_{u,v}(0) = u .
 \end{align}
Further, if \ $(u,v) \ne (0, 0)$, \ then \ $\psi_{u,v}(t) \in \RR_{--}$, \ $t \in \RR_{++}$,
 \ and if \ $(u, v) = (0, 0)$, \ then \ $\psi_{u,v}(t) = 0$, \ $t \in \RR_+$.
\end{Thm}

\noindent{\bf Proof.}
By Theorem 4.10 in Keller-Ressel \cite{KR2},
 \ $\bigl(Y_t, \int_0^t Y_s \, \dd s\bigr)_{t\in\RR_+}$ \ is a
 \ $2$-dimensional CBI process with branching mechanism
 \ $\tR(z_1, z_2) = (\tR_1(z_1, z_2), \tR_2(z_1, z_2))$, \ $z_1, z_2 \in \RR_+$, \ with
 \[
   \tR_1(z_1, z_2) = R(z_1) - z_2, \qquad \tR_2(z_1, z_2) = 0 , \qquad z_1, z_2 \in \RR_+ ,
 \]
 and with immigration mechanism \ $\tF(z_1, z_2) = F(z_1)$, \ $z_1, z_2 \in \RR_+$, \ where
 \ $R$ \ and \ $F$ \ are given in Proposition \ref{Pro_stable_CIR}.
Consequently, by Theorem 2.7 of Duffie et al.\ \cite{DufFilSch}
 (see also Barczy et al.\ \cite[Theorem 2.4]{BarLiPap}), we have
 \begin{align*}
  \EE\left[\exp\left\{uY_t + v \int_0^t Y_s \, \dd s\right\}\right]
  &= \exp\bigg\{y_0 \psi_{u,v}(t)
                - \int_0^t \tF(-\psi_{u,v}(s), -\varphi_{u,v}(s)) \, \dd s \bigg\} \\
  &= \exp\bigg\{y_0 \psi_{u,v}(t) + a \int_0^t \psi_{u,v}(s) \, \dd s \bigg\} ,
 \qquad t \in \RR_+ ,
 \end{align*}
 where the function \ $(\psi_{u,v}, \varphi_{u,v}) : \RR_+ \to \RR_-^2$ \ is the unique locally
 bounded solution to the system of differential equations
 \[
   \begin{cases}
    \psi'_{u,v}(t)
    = \tR_1(-\psi_{u,v}(t),-\varphi_{u,v}(t))
    = \frac{\sigma^2}{2} \psi_{u,v}(t)^2 + \frac{\delta^\alpha}{\alpha} (-\psi_{u,v}(t))^\alpha
      - b \psi_{u,v}(t) + \varphi_{u,v}(t) ,
    \qquad t \in \RR_+ , \\
   \varphi'_{u,v}(t)
   = \tR_2(-\psi_{u,v}(t), -\varphi_{u,v}(t))
   = 0 , \qquad t \in \RR_+ ,
  \end{cases}
 \]
 with initial values \ $\psi_{u,v}(0) = u$, \ $\varphi_{u,v}(0) = v$.
\ Clearly, \ $\varphi_{u,v}(t) = v$, \ $t \in \RR_+$, \ hence we obtain
 \[
   \psi'_{u,v}(t)
   = \frac{\sigma^2}{2} \psi_{u,v}(t)^2 + \frac{\delta^\alpha}{\alpha} (-\psi_{u,v}(t))^\alpha
     - b \psi_{u,v}(t) + v ,
   \qquad t \in \RR_+ , \qquad \psi_{u,v}(0) = u ,
 \]
 as desired.

If \ $(u, v) = (0, 0)$, \ then, since the identically zero function is a (locally bounded)
 solution to \eqref{psi_DE}, by the unicity of such a solution, we have \ $\psi_{0,0}(t) = 0$,
 \ $t \in \RR_+$.
\ If \ $(u, v) \ne (0, 0)$, \ then, on the contrary, let us suppose that there exists an
 \ $t_0 \in \RR_{++}$ \ such that \ $\psi_{u,v}(t_0) = 0$.
\ Let \ $t_* := \inf\{t\in\RR_{++} : \psi_{u,v}(t) = 0\}$.
\ Then \ $t_* < \infty$, \ $\psi_{u,v}(t) < 0$ \ for all \ $t \in [0, t_*)$, \ and
 \ $\psi_{u,v}(t_*) = 0$.
\ If \ $t_* = 0$, \ then \ $0 = \psi_{u,v}(t_*) = \psi_{u,v}(0) = u$, \ and hence
 \ $v \in \RR_{--}$.
\ Further, there exists a sequence \ $(t_n)_{n\in\NN}$ \ such that \ $t_n \in \RR_{++}$,
 \ $n \in \NN$, \ $t_n \downarrow 0 = t_*$ \ as \ $n \to \infty$, \ and \ $\psi_{0,v}(t_n) = 0$,
 \ $n \in \NN$.
\ Consequently, using that a locally bounded solution to \eqref{psi_DE} is unique,
 \ $\psi_{0,v}(kt_n) = \psi_{0,v}(0) = 0$, \ $k, n \in\NN$.
\ Since \ $t_n \downarrow 0$ \ as \ $n \to \infty$, \ for all \ $t \in \RR_{++}$, \ there exists
 a sequence \ $(k_n)_{n\in\NN}$ \ such that \ $k_n \in \NN$, \ $n \in \NN$, \ and
 \ $k_n t_n \to t$ \ as \ $n \to \infty$ \ (one can choose
 \ $k_n: = \lfloor \frac{t}{t_n} \rfloor$, \ $n \in \NN$).
\ Since \ $\psi_{0,v}$ \ is continuous, we have \ $\psi_{0,v}(t) = 0$, \ $t \in \RR_+$,
 \ yielding us to a contradiction (due to \ $v \in \RR_{--}$).
So, if \ $t_* = 0$, \ then \ $\psi_{u,v}(t) < 0$, \ $t > 0$, \ as desired.
In the sequel, let us assume that \ $t_* > 0$.
\ On the one hand, \ $\psi_{u,v}'(t_*) \geq 0$, \ since
 \[
   \psi_{u,v}'(t_*) = \lim_{h\uparrow0} \frac{\psi_{u,v}(t_* + h) - \psi_{u,v}(t_*)}{h}
                    = \lim_{h\uparrow0} \frac{1}{h} \psi_{u,v}(t_* + h) ,
 \]
 where \ $h < 0$ \ and \ $\psi_{u,v}(t_* + h) < 0$ \ yield
 \ $\frac{1}{h} \psi_{u,v}(t_* + h) > 0$.
\ On the other hand, by \eqref{psi_DE}, \ $\psi_{u,v}'(t_*) = v \leq 0$, \ yielding that
 \ $v = 0$.
\ Consequently, if \ $t_* > 0$, \ then \eqref{psi_DE} takes the form
 \begin{align*}
  \psi_{u,v}'(t)
  = \frac{\sigma^2}{2} \psi_{u,v}(t)^2 + \frac{\delta^\alpha}{\alpha} (-\psi_{u,v}(t))^\alpha
    - b\psi_{u,v}(t) ,
  \qquad t \in \RR_+ ,
 \end{align*}
 where \ $\psi_{u,v}(0) = u$ \ and \ $\psi_{u,v}(t_*) = 0$.
\ By the uniqueness of a locally bounded solution of \eqref{psi_DE}, \ $\psi_{u,v}(t) = 0$ \ for
 all \ $t \geq t_*$, \ and hence \ $\psi_{u,v}(t) = 0$ \ for all \ $t \geq 0$.
\ Indeed, let \ $\widetilde\psi_{u,v}(\tau) := \psi_{u,v}(-\tau)$, \ $\tau \leq 0$, \ and
 \ $\tau^*:=-t_*$.
\ Then
 \begin{align*}
  \widetilde\psi_{u,v}'(\tau)
  = - \frac{\sigma^2}{2} \widetilde\psi_{u,v}(\tau)^2
    - \frac{\delta^\alpha}{\alpha} (-\widetilde\psi_{u,v}(\tau))^\alpha
      + b \widetilde\psi_{u,v}(\tau) ,
  \qquad \tau \leq 0 ,
 \end{align*}
 where \ $\widetilde\psi_{u,v}(0) = u$, \ $\widetilde\psi_{u,v}(\tau_*) = 0$, \ and, by the
 uniqueness of a locally bounded solution of the differential equation above,
 \ $\widetilde\psi_{u,v}(\tau) = 0$ \ for \ $\tau\in[\tau_*, 0]$, \ i.e.,
 \ $\psi_{u,v}(-\tau) = 0$ \ for \ $\tau \in [-t_*, 0]$, \ i.e., \ $\psi_{u,v}(t) = 0$ \ for
 \ $t \in [0, t_*]$.
\ This yields us to a contradiction taking into account the definition of \ $t_*$ \ and the fact
 that \ $t_* > 0$.
\proofend

\begin{Rem}\label{Rem_integrals}
\noindent (i)
 By (iii) of Proposition \ref{Pro_stable_CIR}, \ $\psi_{u,0}(t) = - v_t(-u)$ \ for all
  \ $t \in \RR_+$ \ and \ $u \in \RR_-$.

\noindent (ii)
 The differential equation \eqref{psi_DE} is a special case of Polyanin and Zaitsev
  \cite[Section 1.5.1-2/(27)]{PolZai}.
\proofend
\end{Rem}

\section{Existence and uniqueness of MLE}
\label{section_EUMLE}

In this section, we will consider the stable CIR model \eqref{stable_CIR} with known
 \ $a \in \RR_+$, \ $\sigma, \delta \in \RR_{++}$, \ $\alpha \in (1, 2)$, \ and a known
 deterministic initial value \ $Y_0 = y_0 \in \RR_+$, \ and we will consider \ $b \in \RR$ \ as
 an unknown parameter.

Let \ $\PP_b$ \ denote the probability measure induced by \ $(Y_t)_{t\in\RR_+}$ \ on the
 measurable space \ $(D(\RR_+, \RR), \cD(\RR_+, \RR))$ \ endowed with the natural filtration
 \ $(\cD_t(\RR_+, \RR))_{t\in\RR_+}$, \ see Appendix \ref{App_LR}.
Further, for all \ $T \in \RR_{++}$, \ let \ $\PP_{b,T} := \PP_{b\!\!}|_{\cD_T(\RR_+, \RR)}$
 \ be the restriction of \ $\PP_b$ \ to \ $\cD_T(\RR_+, \RR)$.

The next result is about the form of the Radon--Nikodym derivative
 \ $\frac{\dd \PP_{b,T}}{\dd \PP_{\tb,T}}$ \ for \ $b, \tb \in \RR$.
\ We will consider \ $\PP_{\tb,T}$ \  as a fixed reference measure, and we will derive the MLE
 for the parameter \ $b$ \ based on the observations \ $(Y_t)_{t\in[0,T]}$.

\begin{Pro}\label{RNb}
Let \ $a \in \RR_+$, \ $b, \tb \in \RR$, \ $\sigma, \delta \in \RR_{++}$, \ and \ $\alpha\in(1,2)$.
\ Then for all \ $T \in \RR_{++}$, \ the probability measures \ $\PP_{b,T}$ \ and
 \ $\PP_{\tb,T}$ \ are absolutely continuous with respect to each other, and
 \begin{align}\label{RNformulab}
  \log\biggl(\frac{\dd \PP_{b,T}}{\dd \PP_{\tb,T}}(\tY)\biggr)
  = - \frac{b-\tb}{\sigma^2}
      \left(\tY_T - y_0 - a T - \delta \int_0^T \sqrt[\alpha]{\tY_{u-}} \, \dd L_u\right)
    - \frac{b^2-\tb^2}{2\sigma^2} \int_0^T \tY_u \, \dd u
 \end{align}
 holds \ $\PP$-almost surely, where \ $\tY$ \ is the \ $\alpha$-stable CIR process corresponding to the parameter \ $\tb$.
\end{Pro}

\noindent{\bf Proof.}
In what follows, we will apply Theorem III.5.34 in Jacod and Shiryaev \cite{JSh} (see also
 Appendix \ref{App_LR}).
We will work on the canonical space \ $(D(\RR_+, \RR), \cD(\RR_+, \RR))$.
\ Let \ $(\eta_t)_{t\in\RR_+}$ \ denote the canonical process \ $\eta_t(\omega) := \omega(t)$,
 \ $\omega \in D(\RR_+, \RR)$, \ $t \in \RR_+$.
\ Recall that the stable CIR process \eqref{stable_CIR} can be written in the form
 \eqref{YGrigelionis}.
By Proposition \ref{Pro_stable_CIR}, the SDE \eqref{stable_CIR} has a pathwise unique strong
 solution (with the given deterministic initial value \ $y_0 \in \RR_+$), \ and hence, by
 Theorem III.2.26 in Jacod and Shiryaev \cite{JSh}, under the probability measure \ $\PP_b$,
 \ the canonical process \ $(\eta_t)_{t\in\RR_+}$ \ is a semimartingale with semimartingale
 characteristics \ $(B^{(b)}, C, \nu)$ \ associated with the truncation function \ $h$
 \ (introduced in Proposition \ref{Grigelionis}), where
 \begin{align*}
  &B^{(b)}_t
   = \int_0^t
      \left(a - b \eta_u + \gamma \delta \sqrt[\alpha]{\eta_u}
            + \int_\RR
               (h(z\delta\sqrt[\alpha]{\eta_u}) - h(z)\delta\sqrt[\alpha]{\eta_u})
               \, m(\dd z)\right) \dd u ,
   \qquad t \in \RR_+ , \\
  &C_t = \int_0^t (\sigma \sqrt{\eta_u})^2 \, \dd u
       = \sigma^2 \int_0^t \eta_u \, \dd u , \qquad t \in \RR_+ , \\
  &\nu(\dd t, \dd y) = K(\eta_t, \dd y) \, \dd t
 \end{align*}
 with the Borel transition kernel \ $K$ \ from \ $\RR_+ \times \RR$ \ into \ $\RR$ \ given by
 \[
   K(y,R) := \int_{\RR} \!\bbone_{R\setminus\{0\}}(z\delta\sqrt[\alpha]{y}) \, m(\dd z)
   \qquad \text{for \ $y \in \RR_+$ \ and \ $R \in \cB(\RR)$}
 \]
 with \ $m(\dd z) = C_\alpha z^{-1-\alpha} \bbone_{(0,\infty)}(z) \, \dd z$.
\ The aim of the following discussion is to check the set of sufficient conditions presented in
 Appendix \ref{App_LR} (of which the notations will be used) in order to have right to apply
 Theorem III.5.34 in Jacod and Shiryaev \cite{JSh}.
First note that \ $(C_t)_{t\in\RR_+}$ \ and \ $\nu(\dd t, \dd y)$ \ do not depend on the unknown
 parameter \ $b$, \ and hence \ $ V^{(\tb, b)} $ \ is identically one and then \eqref{GIR1} and
 \eqref{GIR2} readily hold.
We also have
 \[
   \PP_b\big(\nu(\{t\}\times \RR) = 0\big)
   = \PP_b\left(\int_{\{t\}} K(\eta_s,\RR) \, \dd s = 0\right) =  1 ,
   \qquad t \in \RR_+ , \quad b \in \RR .
 \]
Further, \ $(C_t)_{t\in\RR_+}$ \ can be represented as \ $C_t = \int_0^t c_u \, \dd F_u$,
 \ $t \in \RR_+$, \ where the stochastic processes \ $(c_t)_{t\in\RR_+}$ \ and
 \ $(F_t)_{t\in\RR_+}$ \ are given by \ $c_t := \sigma^2 \eta_t$, \ $t \in \RR_+$, \ and
 \ $F_t = t$, \ $t \in \RR_+$.
\ Consequently, for all \ $b, \tb \in \RR$,
 \begin{align*}
  B^{(b)}_t - B^{(\tb)}_t
  = - (b - \tb) \int_0^t \eta_u \, \dd u
  = \int_0^t c_u  \beta^{(\tb,b)}_u \, \dd F_u
 \end{align*}
 $\PP_b$-almost surely for every \ $t \in \RR_+$, \ where the stochastic process
 \ $(\beta^{(\tb,b)}_t)_{t\in\RR_+}$ \ is given by
 \[
   \beta^{(\tb,b)}_t  = - \frac{b - \tb}{\sigma^2} , \qquad t \in \RR_+ ,
 \]
 which yields \eqref{GIR3}.
Next we check \eqref{GIR4}, i.e.,
 \begin{gather}\label{COND}
  \PP_b\left(\int_0^t \bigl(\beta^{(\tb,b)}_u\bigr)^2 c_u \, \dd F_u < \infty\right)
  = 1 ,  \qquad t \in \RR_+ .
 \end{gather}
We have
 \[
   \int_0^t \bigl(\beta^{(\tb,b)}_u\bigr)^2 c_u \, \dd F_u
   = \frac{(b-\tb)^2}{\sigma^2} \int_0^t \eta_u \, \dd u , \qquad t \in \RR_+ .
 \]
Since for each \ $\omega \in D(\RR_+, \RR)$, \ the trajectory
 \ $[0, t] \ni u \mapsto \eta_u(\omega)$ \ is c\`adl\`ag, hence bounded (see, e.g., Billingsley
 \cite[(12.5)]{Bil}), we have \ $\int_0^t \eta_u(\omega) \, \dd u < \infty$, \ hence we obtain
 \eqref{COND}.

Next, we check that, under the probability measure \ $\PP_{b}$, \ local uniqueness holds for the
 martingale problem on the canonical space corresponding to the triplet \ $(B^{(b)}, C, \nu)$
 \ with the given initial value \ $y_0$ \ with \ $\PP_b$ \ as its unique solution
 (for the definition of local uniqueness in question, see Definition III.2.27 in Jacod and Shiryaev \cite{JSh}).
By Proposition \ref{Pro_stable_CIR}, the SDE \eqref{stable_CIR} has a pathwise unique strong
 solution (with the given deterministic initial value \ $y_0 \in  \RR_+$), \ and hence Theorem
 III.2.26 in Jacod and Shiryaev \cite{JSh} yields that the set of all solutions to the
 martingale problem on the canonical space corresponding to \ $(B^{(b)}, C, \nu)$ \ has only one
 element \ $(\PP_b)$ \ yielding the desired local uniqueness.
We also mention that Theorem III.4.29 in Jacod and Shiryaev \cite{JSh} implies that under the
 probability measure \ $\PP_b$, \ all local martingales have the integral representation
 property relative to \ $\eta$.

By Theorem III.5.34 in Jacod and Shiryaev \cite{JSh} (see also Appendix \ref{App_LR}),
 \ $\PP_{b,T}$ \ and \ $\PP_{\tb,T}$ \ are equivalent (one can change the roles of \ $b$ \ and
 \ $\tb$), \ and we have
 \[
   \frac{\dd \PP_{b,T}}{\dd \PP_{\tb,T}}(\eta)
   = \exp\bigg\{\int_0^T \beta^{(\tb,b)}_u \, \dd (\eta^\cont)^{(\tb)}_u
                -\frac{1}{2}
                 \int_0^T \bigl(\beta^{(\tb,b)}_u\bigr)^2 c_u \, \dd u\bigg\} ,
   \qquad T \in \RR_{++}
 \]
 holds \ $\PP_{\tb}$-almost surely, where \ $((\eta^\cont)^{(\tb)}_t)_{t\in\RR_+}$ \ denotes the continuous (local) martingale part
 of \ $(\eta_t)_{t\in\RR_+}$ \ under \ $\PP_{\tb}$.
\ Using Remarks III.2.28 in Jacod and Shiryaev \cite{JSh} and \eqref{YGrigelionis},
 the continuous (local) martingale part \ $(\tY^\cont_t)_{t\in\RR_+}$ \ of
 \ $(\tY_t)_{t\in\RR_+}$ \ takes the form
 \ $\tY^\cont_t = \sigma \int_0^t \sqrt{\tY_u} \, \dd W_u$, \ $t \in \RR_+$,
 \ and, by \eqref{stable_CIR}, we have
 \[
   \dd \tY^\cont_t
   = \dd \tY_t - (a - \tb \tY_t) \, \dd t - \delta \sqrt[\alpha]{\tY_{t-}}\dd L_t , \qquad
   t \in \RR_+ .
 \]
Hence
 \begin{align*}
  \log\biggl(\frac{\dd \PP_{b,T}}{\dd \PP_{\tb,T}}(\tY)\bigg)
  &= \int_0^T
      \!\biggl(-\frac{b-\tb}{\sigma^2}\biggr)
      \Bigl(\dd \tY_u - \delta \sqrt[\alpha]{\tY_{u-}} \dd L_u\Bigr)
     - \int_0^T \!\biggl(-\frac{b-\tb}{\sigma^2}\biggr) (a - \tb \tY_u) \, \dd u \\
  &\quad
     - \frac{1}{2} \int_0^T \!\biggl(-\frac{b-\tb}{\sigma^2}\biggr)^2 \sigma^2 \tY_u \, \dd u \\
  &= - \frac{b-\tb}{\sigma^2} \int_0^T (\dd \tY_u - \delta \sqrt[\alpha]{\tY_{u-}}\dd L_u)
     + \frac{b-\tb}{\sigma^2} \int_0^T a \, \dd u
     - \frac{b^2-\tb^2}{ 2\sigma^2} \int_0^T \tY_u \, \dd u
 \end{align*}
 holds \ $\PP$-almost surely, which yields the statement.
\proofend

Next, using Proposition \ref{RNb}, by considering \ $\PP_{\tb,T}$ \ as a fixed reference
 measure, we derive an MLE for the parameter \ $b$ \ based on the observations
 \ $(Y_t)_{t\in[0,T]}$.
\ Let us denote the right hand side of \eqref{RNformulab} by \ $\Lambda_T(b,\tb)$ \ replacing
 \ $\tY$ \ by \ $Y$.
\ By an MLE \ $\hb_T$ \ of the parameter \ $b$ \ based on the observations
 \ $(Y_t)_{t\in[0,T]}$, \ we mean
 \[
   \hb_T := \argmax_{b\in\RR} \Lambda_T(b, \tb) ,
 \]
 which will turn out to be not dependent on \ $\tb$.
\ Our method for deriving an MLE is one of the known ones in the literature, and it turns out
 that these lead to the same estimator \ $\hb_T$, \ see Remark \ref{Luschgy}.

Next, we formulate a result about the unique existence of MLE \ $\hb_T$ \  of \ $b$ \ for all
 \ $T \in \RR_{++}$.

\begin{Pro}\label{LEMMA_MLEb_exist}
Let \ $a \in \RR_+$, \ $b \in \RR$, \ $\sigma, \delta \in \RR_{++}$, \ $\alpha\in(1,2)$, \  and \ $y_0 \in \RR_+$.
\ If \ $a \in \RR_{++}$ \ or \ $y_0 \in \RR_{++}$, \ then for each \ $T \in \RR_{++}$,
 \ there exists a unique MLE \ $\hb_T$ \ of \ $b$ \ $\PP$-almost surely having the form
 \begin{equation}\label{MLEb}
  \hb_T
  = - \frac{Y_T - y_0 - a T - \delta \int_0^T \sqrt[\alpha]{Y_{u-}} \, \dd L_u}
           {\int_0^T Y_s \, \dd s} ,
 \end{equation}
 provided that \ $\int_0^T Y_s \, \dd s \in \RR_{++}$ \ (which holds \ $\PP$-almost surely due to (v) of
 Proposition \ref{Pro_stable_CIR}).
\end{Pro}

\noindent{\bf Proof.}
Due to (v) of Proposition \ref{Pro_stable_CIR},
 \ $\PP\bigl(\int_0^T Y_s \, \dd s \in\RR_{++} \bigr) = 1$ \ for all \ $T \in \RR_{++}$, \ and
 hence the right hand side of \eqref{MLEb} is well-defined \ $\PP$-almost surely.
The aim of the following discussion is to show that the right hand side of \eqref{MLEb} is a
 measurable function of \ $(Y_u)_{u\in[0,T]}$ \ (i.e., a statistic).
By Proposition \ref{RNb}, for all \ $b, \tb \in \RR$ \ and \ $T \in \RR_{++}$, \ the probability
 measures \ $\PP_{\tb,T}$ \ and \ $\PP_{b,T}$ \ are absolutely continuous with respect to each
 other, and
 \[
   \log\biggl(\frac{\dd \PP_{\tb,T}}{\dd \PP_{b,T}}(Y)\biggr)
   = - \frac{\tb-b}{\sigma^2}
       \left(Y_T - y_0 - a T - \delta \int_0^T \sqrt[\alpha]{Y_{u-}} \, \dd L_u\right)
     - \frac{\tb^2 - b^2}{2\sigma^2} \int_0^T Y_u \, \dd u
 \]
 holds \ $\PP$-almost surely.

The left-hand side of the above equality is measurable with respect to \ $(Y_t)_{t\in[0,T]}$
 \ (see, e.g.,  Jacod and Shiryaev \cite[Theorem III.3.4]{JSh}), and hence its right hand-side
 is also measurable, which yields the measurability  of
 \ $\int_0^T \sqrt[\alpha]{Y_{u-}} \, \dd L_u$ \ with respect to \ $(Y_t)_{t\in[0,T]}$ \ and
 consequently that of \ $\hb_T$.

By Proposition \ref{RNb}, for all \ $b,\tb\in\RR$, \ we have
 \begin{align*}
  &\frac{\partial}{\partial b}
   \log(\Lambda_T(b,\tb))
   = - \frac{1}{\sigma^2}
       \left(Y_T - y_0 - a T - \delta \int_0^T\sqrt[\alpha]{Y_{u-}}\dd L_u \right)
     - \frac{b}{\sigma^2} \int_0^T Y_u \, \dd u , \\
  &\frac{\partial^2}{\partial b^2}
   \log(\Lambda_T(b,\tb))
   = - \frac{1}{\sigma^2} \int_0^T Y_u \, \dd u .
 \end{align*}
Thus the MLE \ $\hb_T$ \ of \ $b$ \ based on a continuous time observation \ $(Y_s)_{s\in[0,T]}$ \ exists \ $\PP$-almost
 surely, and it takes the form \eqref{MLEb} provided that \ $\int_0^T Y_s \, \dd s \in\RR_{++}$.
\proofend

\begin{Rem}\label{Rem_L_measurability}
In what follows, under the assumptions of Proposition \ref{LEMMA_MLEb_exist} and the additional
 assumption \ $a \geq \frac{\sigma^2}{2}$, \ we prove that \ $L_t$ \ is a measurable function of
 \ $(Y_u)_{u\in[0,T]}$ \ for all \ $t \in [0, T]$, \ where \ $T \in \RR_{++}$, \ which (in the
 special case \ $a \geq \frac{\sigma^2}{2}$) \ gives another proof for the fact that the right
 hand side of \eqref{MLEb} is a statistic.
Recalling the notation \ $\zeta_t = \delta\int_0^t \sqrt[\alpha]{Y_{u-}}\,\dd L_u$,
 \ $t \in \RR_+$, \ we have \ $\Delta \zeta_t = \delta \sqrt[\alpha]{Y_{t-}} \Delta L_t$,
 \ $t \in \RR_+$ \ (following from \eqref{help_Levy_Ito_decomposition} and Jacod and Shiryaev
 \cite[Definitions II.1.27]{JSh}) and using the SDE \eqref{stable_CIR}, we get
 \ $\Delta Y_t = \Delta \zeta_t = \delta \sqrt[\alpha]{Y_{t-}}\Delta L_t$, \ $t \in \RR_+$.
\ Hence \ $\Delta L_t = \frac{\Delta Y_t}{\delta \sqrt[\alpha]{Y_{t-}}}$, \ $t \in \RR_{++}$,
 \ since, by (iv) of Proposition \ref{Pro_stable_CIR},
 \ $\PP\bigl(Y_t\in\RR_{++} \ \text{for all \ $t\in\RR_{++}$} \bigr) = 1$.
\ For all \ $t \in [0, T]$ \ and \ $\vare \in (0, 1)$,
 \begin{equation*}
  \begin{aligned}
   \int_{(0,t]} \int_{\{\vare<|z|\leq 1\}} z \, \tmu^L(\dd u,\dd z)
   &= \sum_{u\in[0,t]} \bbone_{\{ \vare < |\Delta L_u|\leq 1 \}} \Delta L_u
      - \int_{(0,t]} \int_{\{\vare<|z|\leq 1\}} z \, \dd u \, m(\dd z) \\
   &= \sum_{u\in[0,t]}
       \bbone_{\big\{\vare < \frac{|\Delta Y_u|}{\delta \sqrt[\alpha]{Y_{u-}}} \leq 1 \big\}}
       \frac{\Delta Y_u}{\delta \sqrt[\alpha]{Y_{u-}}}
      - t \int_{\{\vare<|z|\leq 1\}} z \, m(\dd z) ,
  \end{aligned}
 \end{equation*}
 which is a measurable function of \ $(Y_t)_{t\in[0,T]}$.
\ Similarly, for all \ $t \in [0, T]$,
 \begin{equation*}
  \begin{aligned}
   \int_{(0,t]} \int_{\{|z|> 1\}} z \, \mu^L(\dd u,\dd z)
   = \sum_{u\in[0,t]} \bbone_{\{|\Delta L_u|>1\}} \Delta L_u
   = \sum_{u\in[0,t]}
      \bbone_{\big\{\frac{|\Delta Y_u|}{\delta\sqrt[\alpha]{Y_{u-}}}>1\big\}}
      \frac{\Delta Y_u}{\delta \sqrt[\alpha]{Y_{u-}}},
 \end{aligned}
 \end{equation*}
 which is a measurable function of \ $(Y_u)_{u\in[0,T]}$ \ as well.
Hence, using \eqref{help_Levy_Ito_decomposition}, for all \ $t \in [0, T]$,
 \begin{equation*}
  \sum_{u\in[0,t]}
   \bbone_{\big\{\frac{|\Delta Y_u|}{\delta\sqrt[\alpha]{Y_{u-}}}>\vare\big\}}
   \frac{\Delta Y_u}{\delta\sqrt[\alpha]{Y_{u-}}}
  - t \int_{\{\vare<|z|\leq 1\}} z \, m(\dd z)
  + \gamma t
  \stoch L_t \quad\text{as \ $\vare \downarrow 0$,}
 \end{equation*}
 yielding that \ $L_t$ \ is a measurable function of \ $(Y_t)_{t\in[0,T]}$.
\ Finally, note also that if \ $Y_{t-} = 0$, \ $t \in \RR_+$, \ then, using
 \ $\Delta Y_t = \delta \sqrt[\alpha]{Y_{t-}} \Delta L_t$, \ $t \in \RR_+$, \ we have
 \ $\Delta Y_t = 0$ \ yielding \ $Y_t = Y_{t-} =0$ \ (irrespective of the size of the jump of
 \ $L$ \ at \ $t$).
\proofend
\end{Rem}

\begin{Rem}\label{Luschgy}
In the literature there is another way of deriving an MLE.
S{\o}rensen \cite{SorM} defined an MLE of \ $b$ \ as a solution of the equation
 \ $\dot{\Lambda}_T(b) = 0$, \ where \ $\dot{\Lambda}_T(b)$ \ is the so-called score vector
 given in formula (3.3) in S{\o}rensen \cite{SorM}.
Luschgy \cite{Lus2}, \cite{Lus} called this equation as an estimating equation.
With the notations of the proof of Proposition \ref{RNb}, taking into account of the form of
 \ $\beta^{(\tb,b)}$ \ and the fact that \ $V^{(\tb,b)}$ \ is identically one, we have
 \begin{align*}
  \dot{\Lambda}_T(b)
  &:= \int_0^T \biggl(-\frac{1}{\sigma^2}\biggr) \, \dd Y^\cont_u
   = - \frac{1}{\sigma^2}
       \int_0^T (\dd Y_u - (a - b Y_u) \, \dd u - \delta \sqrt[\alpha]{Y_{u-}} \, \dd L_u) \\
  &= - \frac{1}{\sigma^2}
       \biggl(Y_T - y_0 - a T + b \int_0^T Y_u \, \dd u
              - \delta \int_0^T \sqrt[\alpha]{Y_{u-}} \, \dd L_u\biggr)
 \end{align*}
 for \ $b \in \RR$ \ and \ $T \in \RR_{++}$.
\ The estimating equation \ $\dot{\Lambda}_T(b) = 0$, \ $b \in \RR$, \ has a unique solution
 \ $- \frac{Y_T-y_0-aT-\delta \int_0^T \sqrt[\alpha]{Y_{u-}} \, \dd L_u}{\int_0^T Y_u\,\dd u}$
 \ provided that \ $\int_0^T Y_u\,\dd u$ \ is strictly positive, which holds \ $\PP$-almost surely
 provided that \ $a \in \RR_{++}$ \ or \ $y_0 \in \RR_{++}$ \ (due to (v) of Proposition
 \ref{Pro_stable_CIR}).
Recall that this unique solution coincides with \ $\hb_T$, \ see \eqref{MLEb}.
\proofend
\end{Rem}

\section{Asymptotic behavior of the MLE in the subcritical case}
\label{section_MLE_subcritical}

If \ $a \in \RR_{++}$ \ or \ $y_0 \in \RR_{++}$, \ then, using \eqref{MLEb} and the SDE
 \eqref{stable_CIR}, we get
 \begin{align}\label{MLEb-}
  \hb_T - b
  = - \frac{Y_T - y_0 - a T - \delta \int_0^T \sqrt[\alpha]{Y_{u-}} \, \dd L_u
            + b \int_0^T Y_s \, \dd s}
           {\int_0^T Y_s \, \dd s}
  = - \frac{\sigma \int_0^T \sqrt{Y_s} \, \dd W_s}{\int_0^T Y_s \, \dd s}
 \end{align}
 provided that \ $\int_0^T Y_s \, \dd s \in\RR_{++}$, \ which holds \ $\PP$-almost surely due to (v) of
 Proposition \ref{Pro_stable_CIR}.
Here note that \ $\sigma \int_0^T \sqrt{Y_s} \, \dd W_s = Y^{\mathrm{cont}}_T$, \ $T \in \RR_+$,
 \ due to part 1) of Remarks III.2.28 in Jacod and Shyriaev \cite{JSh} and \eqref{YGrigelionis}.

\begin{Thm}\label{Thm_MLEb_subcritical}
Let \ $a,b,\sigma,\delta \in \RR_{++}$ \ and \ $\alpha\in(1,2)$.
\ Let \ $(Y_t)_{t\in\RR_+}$ \ be the unique strong solution of the SDE \eqref{stable_CIR}
 satisfying \ $\PP(Y_0 = y_0) = 1$ \ with some \ $y_0 \in \RR_+$.
\ Then the MLE \ $\hb_T$ \ of \ $b$ \ is strongly consistent and asymptotically normal, i.e.,
 \ $\hb_T \as b$ \ as \ $T \to \infty$, \ and
 \begin{align}\label{help_conv_subcritical}
  \sqrt{T} (\hb_T - b) \distr \cN\left(0, \frac{\sigma^2b}{a}\right)
  \qquad \text{as \ $T \to \infty$.}
 \end{align}
With a random scaling,
 \[
   \frac{1}{\sigma} \biggl(\int_0^T Y_s \, \dd s\biggr)^{1/2} (\hb_T - b) \distr \cN(0, 1)
   \qquad \text{as \ $T \to \infty$.}
 \]
\end{Thm}

\noindent{\bf Proof.}
By Proposition \ref{LEMMA_MLEb_exist}, there exists a unique MLE \ $\hb_T$ \ of \ $b$ \ for
 all \ $T \in \RR_{++}$, \ which has the form given in \eqref{MLEb}.
By (i) of Theorem \ref{Ergodicity}, \ $(Y_t)_{t\in\RR_+}$ \ has a unique stationary
 distribution \ $\pi$ \ with \ $\int_0^\infty y \, \pi(\dd y) = \frac{a}{b} \in \RR_{++}$.
\ By (ii) of Theorem \ref{Ergodicity}, \ we have
 \ $\frac{1}{T} \int_0^T Y_s \, \dd s \as \int_0^\infty y \, \pi(\dd y)$ \ as
 \ $T \to \infty$, \ implying also \ $\int_0^T Y_s \, \dd s \as \infty$ \ as \ $T \to \infty$.
\ Since the quadratic variation process of the square integrable martingale
 \ $\bigl(\int_0^t \sqrt{Y_s} \, \dd W_s\bigr)_{t\in\RR_+}$ \ takes the form
 \ $\bigl(\int_0^t Y_s \, \dd s\bigr)_{t\in\RR_+}$, \ using \eqref{MLEb-} and Theorem
 \ref{DDS_stoch_int}, we have \ $\hb_T \as b$ \ as \ $T \to \infty$.
\ Moreover, using Theorem \ref{THM_Zanten} with
 \ $\eta := \left(\int_0^\infty y \, \pi(\dd y)\right)^{1/2}$ \ and Slutsky's lemma, we have
 \[
   \sqrt{T} (\hb_T - b)
   = -\sigma \frac{\frac{1}{\sqrt{T}} \int_0^T \sqrt{Y_s} \, \dd W_s}
                  {\frac{1}{T} \int_0^T Y_s \, \dd s}
   \distr -\sigma
           \frac{\left(\int_0^\infty y\,\pi(\dd y)\right)^{1/2}\,\cN(0,1)}
                {\int_0^\infty y\,\pi(\dd y)}
   = \cN\biggl(0, \frac{\sigma^2}{\int_0^\infty y\,\pi(\dd y)}\biggr)
 \]
 as \ $T \to \infty$, \ hence, by \eqref{help_stac_expectation}, we obtain
 \eqref{help_conv_subcritical}.
Further, Slutsky's lemma yields
 \begin{align*}
  \frac{1}{\sigma} \left(\int_0^T Y_s \, \dd s \right)^{1/2} (\hb_T - b)
  &= \frac{1}{\sigma} \left(\frac{1}{T} \int_0^T Y_s \, \dd s \right)^{1/2}
     \sqrt{T} (\hb_T - b) \\
  &\distr \frac{1}{\sigma} \left(\int_0^\infty y\,\pi(\dd y)\right)^{1/2} \,
          \cN\biggl(0, \frac{\sigma^2}{\int_0^\infty y\,\pi(\dd y)}\biggr)
   = \cN(0, 1)
 \end{align*}
 as \ $T \to \infty$.
\proofend

\section{Consistency of the MLE in the critical case}
\label{section_MLE_critical}

First, we describe the asymptotic behavior of the solution \ $\psi_{u,v}$ \ of the differential
 equation \eqref{psi_DE} as \ $t \to \infty$ \ in case of \ $b = 0$ \ using a so-called
 separator technique.

Note that, for \ $b,\sigma \in \RR_+$, \ $\delta \in \RR_{++}$ \ and \ $\alpha\in(1,2)$, \ the
 function
 \ $\RR_{-} \ni x \mapsto \tR(x) := R(-x)
    = \frac{\sigma^2}{2} x^2 + \frac{\delta^\alpha}{\alpha} (-x)^\alpha -b x$
 \ is strictly monotone decreasing, continuous, convex
 and \ $\lim_{x\to-\infty} \tR (x) = + \infty$.
\ Indeed, \ $\tR'(x) = \sigma^2 x - \delta^\alpha(-x)^{\alpha-1} - b < 0$, \ $x \in \RR_{--}$,
 \ and \ $\tR''(x) = \sigma^2 + \delta^\alpha (\alpha-1) (-x)^{\alpha-2} > 0$,
 \ $x \in \RR_{--}$.
\ Hence, for all \ $v \in \RR_{--}$, \ the equation
 \ $R(-x) = \frac{\sigma^2}{2} x^2 + \frac{\delta^\alpha}{\alpha} (-x)^\alpha - b x = - v$,
 \ $x \in \RR_{-}$, \ has a unique negative solution, denoted by \ $c_v \in \RR_{--}$.
\ Further, if \ $x \in (c_v, 0]$, \ then \ $R(-x) < - v$, \ i.e., \ $R(-x) + v < 0$; \ and if
 \ $x \in (-\infty, c_v)$, \ then \ $R(-x) > - v$, \ i.e., \ $R(-x) + v > 0$.

\begin{Pro}\label{Pro_critical_psi_asymp}
Let \ $a \in \RR_+$, \ $b = 0$, \ $\sigma \in \RR_+$, \ $\delta \in \RR_{++}$, \ and \ $\alpha\in(1,2)$.
\ Then for all \ $u \in \RR_{-}$ \ and \ $v \in \RR_{--}$, \ the unique locally bounded solution
 \ $\psi_{u,v}$ \ of the differential equation \eqref{psi_DE} satisfies
  \ $\lim_{t\to\infty} \psi_{u,v}(t) = c_v$.
\ Further, \ $\psi_{u,v}(t) \in (c_v, 0]$, \ $t \in \RR_+$ \ if \ $c_v < u \leq 0$;
 \ $\psi_{u,v}(t) \in (-\infty, c_v)$, \ $t \in \RR_+$ \ if \ $u < c_v$; \ and
 \ $\psi_{u,v}(t) = c_v$, \ $t \in \RR_+$ \ if \ $u = c_v$.
\end{Pro}

\noindent{\bf Proof.}
Let \ $v \in \RR_{--}$ \ be fixed.
By Theorem \ref{Thm_Laplace}, \ $\psi_{u,v}(t) \in \RR_{--}$ \ for all \ $t \in \RR_{++}$.
\ Let
 \ $Q(x) := R(-x) + v = \frac{\sigma^2}{2} x^2 + \frac{\delta^\alpha}{\alpha} (-x)^\alpha + v$,
 \ $x \in \RR_{-}$.
\ Note that \ $Q$ \ is continuously differentiable on \ $\RR_{--}$.
\ Further, \ $Q(x) = 0$, \ $x \in \RR_{-}$, \ holds if and only if \ $x = c_v$, \ and
 \ $Q(x) < 0$ \ for \ $x \in (c_v, 0]$ \ and \ $Q(x) > 0$ \ for \ $x \in (-\infty, c_v)$.

If \ $\psi_{u,v}(0) = u = c_v$, \ then the unique locally bounded solution of the differential
 equation \eqref{psi_DE} takes the form \ $\psi_{u,v}(t) = c_v$, \ $t \in \RR_+$, \ since in this
 case \ $\psi_{u,v}'(t) = 0$, \ $t \in \RR_+$, \ and
 \ $\frac{\sigma^2}{2} (\psi_{u,v}(t))^2 + \frac{\delta^\alpha}{\alpha} (-\psi_{u,v}(t))^\alpha
    + v
    = \frac{\sigma^2}{2}c_v^2 + \frac{\delta^\alpha}{\alpha}(-c_v)^\alpha + v
    = Q(c_v) = 0$,
 \ $t \in \RR_+$, \ and hence \eqref{psi_DE} holds.

If \ $\psi_{u,v}(0) = u > c_v$, \ then \ $\psi_{u,v}(t) > c_v$ \ for all \ $t \in \RR_+$.
\ Indeed, on the contrary, let us suppose that there exists \ $t_0 \in \RR_{++}$ \ such that
 \ $\psi_{u,v}(t_0) = c_v$ \ (which can be supposed due to the continuity of \ $\psi_{u,v}$).
\ Then \ $\psi_{u,v}(t) = c_v$ \ would hold for all \ $t \in \RR_+$, \ since it is known that
 if two maximal solutions of an autonomous ordinary differential equation (with a continuously
 differentiable function on the right hand side) coincide at some points, then their ranges
 coincide, see, e.g., Arnol'd \cite[Corollary on page 118]{Arn},
 and the identically \ $c_v$ \ function is a solution of \eqref{psi_DE} (without the initial
 value).
Since \ $\psi_{u,v}(0) > c_v$, \ this leads us to a contradiction.
Consequently, by \eqref{psi_DE}, \ $\psi_{u,v}'(t) = Q(\psi_{u,v}(t)) < 0$, \ $t \in \RR_+$,
 \ yielding that \ $\psi_{u,v}$ \ is monotone decreasing.
Since \ $\psi_{u,v}$ \ is bounded below by \ $c_v$, \ there exists an
 \ $\widetilde c_v \in \RR_{-}$ \ such that \ $\widetilde c_v \geq c_v$ \ and
 \ $\lim_{t\to\infty} \psi_{u,v}(t) = \widetilde c_v$.
\ It remains to check that \ $\widetilde c_v = c_v$.
\ We show that \ $Q(\widetilde c_v) = 0$, \ yielding \ $\widetilde c_v = c_v$, \ since \ $c_v$
 \ is the only root of \ $Q(x) = 0$, \ $x \in \RR_{-}$.
\ On the contrary, let us suppose that \ $Q(\widetilde c_v) > 0$ \ (the case
 \ $Q(\widetilde c_v) < 0$ \ can be handled similarly).
Since \ $Q$ \ is continuous at \ $\widetilde c_v$, \ there exists \ $\kappa > 0$ \ such that
 \ $Q(x) > \frac{Q(\widetilde c_v)}{2}$ \ for all \ $x \in \RR_{-}$ \ satisfying
 \ $|x - \widetilde c_v| < \kappa$.
\ Since \ $\lim_{t\to\infty} \psi_{u,v}(t) = \widetilde c_v$, \ there exists \ $T > 0$ \ such
 that \ $|\psi_{u,v}(t) - \widetilde c_v| < \kappa$ \ for \ $t \geq T$.
\ Hence
 \[
   \psi_{u,v}'(t) = Q(\psi_{u,v}(t)) > \frac{Q(\widetilde c_v)}{2} , \qquad t \geq T .
 \]
Integrating over \ $[T, t]$, \ we have
 \[
  \psi_{u,v}(t) - \psi_{u,v}(T) \geq \frac{Q(\widetilde c_v)}{2}(t-T) , \qquad t \geq T .
 \]
Since, by assumption, \ $Q(\widetilde c_v) > 0$, \ taking the limit \ $t \to \infty$, \ we get
 \ $\lim_{t\to\infty}\psi_{u,v}(t) = \infty$, \ yielding us to a contradiction.

The case \ $\psi_{u,v}(0) = u < c_v$ \ can be handled similarly, and the other parts of the
 proposition follow as well.
\proofend

\begin{Thm}\label{Thm_MLE_critical}
Let \ $a \in \RR_{++}$, \ $b = 0$, \ $\sigma, \, \delta \in \RR_{++}$, \ and \ $\alpha\in(1,2)$.
\ Let \ $(Y_t)_{t\in\RR_+}$ \ be the unique strong solution of the SDE \eqref{stable_CIR}
 satisfying \ $\PP(Y_0 = y_0) = 1$ \ with some \ $y_0 \in \RR_+$.
\ Then the MLE of \ $b$ \ is strongly consistent, i.e., \ $\hb_T \as b$ \ as \ $T \to \infty$.
\end{Thm}

\noindent{\bf Proof.}
Since \ $\big(\int_0^t Y_s \, \dd s\big)_{t\in\RR_+}$ \ is monotone increasing \ $\PP$-almost surely,
 there exists an \ $[0, \infty]$-valued random variable \ $\xi$ \ such that
 \ $\int_0^t Y_s \, \dd s \as \xi$ \ as \ $t \to \infty$, \ and consequently,
 by the dominated convergence theorem, \ $\lim_{t\to\infty} \EE\big[\exp\big\{v \int_0^t Y_s \, \dd s\big\}\big]
   = \EE[\ee^{v\xi}]$ \ for any \ $v\in\RR_{-}$.
\ By Theorem \ref{Thm_Laplace} with \ $b = 0$, \ we have
 \[
   \EE\left[\exp\left\{v \int_0^t Y_s \, \dd s\right\}\right]
      = \exp\bigg\{y_0 \psi_{0,v}(t) + a \int_0^t \psi_{0,v}(s) \dd s\bigg\},
      \qquad t\in\RR_+,\quad v\in\RR_{-}.
 \]
First we check that
 \begin{equation}\label{help_DE2}
  \int_0^t \psi_{0,v}(s) \, \dd s
  = \int_0^{\psi_{0,v}(t)}
     \frac{x}{\frac{\sigma^2}{2}x^2 + \frac{\delta^\alpha}{\alpha}(-x)^\alpha + v} \, \dd x
 \end{equation}
 for all \ $t \in \RR_{++}$ \ and \ $v \in \RR_{--}$, \ where the function
 \ $\psi_{0,v} : \RR_+ \to \RR_-$ \ is given by \eqref{psi_DE}.
\ Recall that, for all \ $v \in \RR_{--}$, \ $c_v \in \RR_{--}$ \ denotes the unique negative
 solution of the equation
 \ $\frac{\sigma^2}{2} x^2 + \frac{\delta^\alpha}{\alpha} (-x)^\alpha + v = 0$, \ we have
 \ $\frac{\sigma^2}{2}x^2 + \frac{\delta^\alpha}{\alpha}(-x)^\alpha + v < 0$ \ for all
 \ $x \in (c_v, 0]$, \ and, by Proposition \ref{Pro_critical_psi_asymp},
 \ $\psi_{0,v}(t) \in (c_v, 0]$ \ for all \ $t \in \RR_+$.
\ Consequently, by \eqref{psi_DE}, the function
 \ $[0, t] \ni s \mapsto \psi_{0,v}(s) \in (c_v, 0]$ \ is strictly decreasing and continuously
 differentiable, hence, for all \ $t \in \RR_{++}$, \ by the substitution \ $x = \psi_{0,v}(s)$
 \ we obtain
 \[
   \int_0^t \psi_{0,v}(s) \, \dd s
   = \int_0^{\psi_{0,v}(t)} \frac{x}{\psi_{0,v}'(\psi_{0,v}^{-1}(x))} \, \dd x
 \]
 and hence \eqref{help_DE2}, where \ $\psi_{0,v}^{-1}$ \ denotes the inverse of \ $\psi_{0,v}$.
\ By \eqref{help_DE2}, we have
 \begin{align*}
  \EE\left[\exp\left\{v \int_0^t Y_s \, \dd s\right\}\right]
   = \exp\bigg\{y_0 \psi_{0,v}(t)
                + a \int_0^{\psi_{0,v}(t)}
                     \frac{x}
                          {\frac{\sigma^2}{2}x^2+\frac{\delta^\alpha}{\alpha}(-x)^\alpha+v}
                     \, \dd x\bigg\}
 \end{align*}
 for \ $t \in \RR_{++}$ \ and \ $v \in \RR_{--}$.
\ Then, by Proposition \ref{Pro_critical_psi_asymp}, \ $\lim_{t\to\infty} \psi_{0,v}(t) = c_v$,
 \ and hence
 \[
   \lim_{t\to\infty}
    \int_0^{\psi_{0,v}(t)}
     \frac{x}{\frac{\sigma^2}{2}x^2+\frac{\delta^\alpha}{\alpha}(-x)^\alpha+v} \, \dd x
   = \int_0^{c_v}
      \frac{x}{\frac{\sigma^2}{2}x^2+\frac{\delta^\alpha}{\alpha}(-x)^\alpha+v} \, \dd x
   = - \infty .
 \]
Here the last step can be checked as follows.
We have
 \[
   \lim_{x\to c_v} \frac{\frac{\sigma^2}{2}x^2+\frac{\delta^\alpha}{\alpha}(-x)^\alpha+v}{x-c_v}
   = \lim_{x\to c_v} \frac{(R(-x)+v)-(R(-c_v)+v)}{x-c_v}
   = - R'(-c_v)
   = \sigma^2 c_v - \delta^\alpha (-c_v)^{\alpha-1}
   < 0 ,
 \]
 thus there exists \ $x_0 \in (c_v, 0)$ \ such that
 \ $\frac{\frac{\sigma^2}{2}x^2+\frac{\delta^\alpha}{\alpha}(-x)^\alpha+v}{x-c_v}
    \leq \frac{\sigma^2 c_v - \delta^\alpha (-c_v)^{\alpha-1}}{2}$
 \ for all \ $x \in (c_v, x_0)$.
\ Hence
 \begin{align*}
  \int_0^{c_v} \frac{x}{\frac{\sigma^2}{2}x^2+\frac{\delta^\alpha}{\alpha}(-x)^\alpha+v} \, \dd x
  &\leq \int_{x_0}^{c_v}
         \frac{x}{\frac{\sigma^2}{2}x^2+\frac{\delta^\alpha}{\alpha}(-x)^\alpha+v} \, \dd x \\
  &\leq \frac{2x_0}{\sigma^2 c_v - \delta^\alpha (-c_v)^{\alpha-1}}
        \int_{x_0}^{c_v} \frac{1}{x-c_v} \, \dd x
   = - \infty ,
 \end{align*}
 as desired.
Hence, since \ $a \in \RR_{++}$, \ we have
 \begin{align*}
  \lim_{t\to\infty} \EE\left[\exp\left\{v \int_0^t Y_s \, \dd s\right\}\right]
  = \exp\left\{y_0 c_v + a (-\infty)\right\}
  = 0 , \qquad v \in \RR_{--} ,
 \end{align*}
 which yields that \ $\EE(\ee^{v\xi}) = 0$, \ $v \in \RR_{--}$.
\ Since, for all \ $v \in \RR_{--}$,
 \begin{align*}
  \EE(\ee^{v\xi})
  &= \EE(\ee^{v\xi} \mid \xi = \infty) \PP(\xi = \infty)
     + \EE(\ee^{v\xi} \mid \xi < \infty) \PP(\xi < \infty) \\
  &= 0 \cdot \PP(\xi = \infty) + \EE(\ee^{v\xi} \mid \xi < \infty) \PP(\xi < \infty) ,
 \end{align*}
 we have \ $0 = \EE(\ee^{v\xi} \mid \xi < \infty) \PP(\xi < \infty)$, \ yielding that
 \ $\PP(\xi < \infty) = 0$, \ i.e., \ $\PP(\xi = \infty) = 1$.
\ That is, we have proved that \ $\int_0^t Y_s \, \dd s \as \infty$ \ as \ $t \to \infty$.
\ Since the quadratic variation process of the square integrable martingale
 \ $\bigl(\int_0^t \sqrt{Y_s} \, \dd W_s\bigr)_{t\in\RR_+}$ \ takes the form
 \ $\bigl(\int_0^t Y_s \, \dd s\bigr)_{t\in\RR_+}$, \ using \eqref{MLEb-} and Theorem
 \ref{DDS_stoch_int}, we have \ $\hb_T \as b$ \ as \ $T \to \infty$, \ as desired.
\proofend

In the critical case the description of the asymptotic behavior of the MLE in question
 remains open.

\section{Asymptotic behavior of the MLE in the supercritical case}
\label{section_MLE_supercritical}

\begin{Thm}\label{Thm_supercritical_convergence}
Let \ $a \in \RR_+$, \ $b \in \RR_{--}$, \ $\sigma \in \RR_+$, \ $\delta \in \RR_{++}$, \ and \ $\alpha\in(1,2)$.
\ Let \ $(Y_t)_{t\in\RR_+}$ \ be the unique strong solution of the SDE \eqref{stable_CIR}
 satisfying \ $\PP(Y_0 = y_0) = 1$ \ with some \ $y_0 \in \RR_+$.
\renewcommand{\labelenumi}{{\rm(\roman{enumi})}}
\begin{enumerate}
 \item
  Then there exists a random variable \ $V$ \ with \ $\PP(V \in \RR_+) = 1$ \ such that
   \[
     \ee^{bt} Y_t \as V \qquad \text{and} \qquad
     \ee^{bt} \int_0^t Y_u \, \dd u \as -\frac{V}{b}
     \qquad \text{as \ $t \to \infty$.}
   \]
 \item
  Moreover, the Laplace transform of \ $V$ \ takes the form
   \begin{align}\label{Laplace_supercritical_limit}
    \EE(\ee^{u V})
    = \exp\biggl\{y_0 \psi_u^*  + \int_0^{-\psi_u^*} \frac{F(z)}{R(z)} \, \dd z\biggr\} ,
    \qquad u \in \RR_- ,
   \end{align}
   where the functions \ $F$ \ and \ $R$ \ are given in Proposition \ref{Pro_stable_CIR}, and
   \ $\psi_u^*:= \lim_{t\to\infty} \psi_{u\ee^{bt},0}(t)$, \ where the function
   \ $\psi_{u\ee^{bt},0} : \RR_+ \to \RR_-$ \ is given by \eqref{psi_DE}.
 \item
  Further, \ $\psi_0^* = 0$ \ and  \ $\psi_u^* = - K^{-1}(-u)$ \ for all \ $u \in \RR_{--}$,
   \ where \ $K^{-1}$ \ denotes the inverse of the strictly increasing function
   \ $K : (0, \theta_0) \to \RR_{++}$ \ given by
   \[
     K(\lambda)
     := \lambda
        \exp\biggl\{\int_0^\lambda \biggl(\frac{b}{R(z)} - \frac{1}{z}\biggr) \dd z\biggr\} ,
     \qquad \lambda \in (0, \theta_0) ,
   \]
   where \ $\theta_0 = \inf\{z \in \RR_{++} : R(z) \in \RR_+\} \in \RR_{++}$.
 \item
  If, in addition, \ $a \in \RR_{++}$, \ then \ $\PP(V \in \RR_{++}) = 1$.
\end{enumerate}
\end{Thm}

In the next remark we present more properties of \ $\psi_u^*$, \ $u \in \RR_{--}$.

\begin{Rem}\label{Rem_super_rep}
\noindent
(i) For all \ $\lambda, \theta \in (0, \theta_0)$, \ we have
 \[
   \frac{K(\theta)}{K(\lambda)}
   = \frac{\theta}{\lambda}
     \exp\biggl\{\int_\lambda^\theta \biggl(\frac{b}{R(z)} - \frac{1}{z}\biggr) \dd z\biggr\}
   = \frac{\theta}{\lambda}
     \exp\biggl\{\int_\lambda^\theta \frac{b}{R(z)} \, \dd z
                 - (\log(\theta) - \log(\lambda))\biggr\}
   = \exp\biggl\{\int_\lambda^\theta \frac{b}{R(z)} \, \dd z\biggr\} ,
 \]
 hence
 \[
   \int_\lambda^\theta \frac{b}{R(z)} \, \dd z
   = \log\biggl(\frac{K(\theta)}{K(\lambda)}\biggr) .
 \]
Consequently, for all \ $\lambda \in (0, \theta_0)$ \ and \ $u \in \RR_{--}$, \ we conclude
 \[
   \int_\lambda^{-\psi_u^*} \frac{b}{R(z)} \, \dd z
   = \log\biggl(\frac{K(-\psi_u^*)}{K(\lambda)}\biggr)
   = \log\biggl(\frac{K(K^{-1}(-u))}{K(\lambda)}\biggr)
   = \log\biggl(-\frac{u}{K(\lambda)}\biggr) .
 \]

\noindent
(ii) Using the formula for the derivative of an inverse function, we have
 \[
   \frac{\dd}{\dd u}\psi_u^*
    = \frac{1}{K'(K^{-1}(-u))}
    = \frac{1}{K'(-\psi_u^*)}
    = \frac{R(-\psi_u^*)}{bK(-\psi_u^*)} , \qquad
    u \in \RR_{--} . \\[-6mm]
 \]
\proofend
\end{Rem}

\noindent{\bf Proof of Theorem \ref{Thm_supercritical_convergence}.}
First we check that the function \ $K$ \ is well-defined and strictly increasing on
 \ $(0, \theta_0)$.
\ Observe that
 \ $R(z) = \frac{\sigma^2}{2} z^2 + \frac{\delta^\alpha}{\alpha}z^\alpha + b z \in \RR_{--}$
 \ for all \ $z \in (0, \theta_0)$.
\ Moreover, we have
 \begin{equation}\label{b/R(z)-1/z}
  \frac{b}{R(z)} - \frac{1}{z}
    = - \frac{\frac{\sigma^2}{2}z+\frac{\delta^\alpha}{\alpha}z^{\alpha-1}}{R(z)} \in \RR_{++} ,
  \qquad z \in (0, \theta_0) .
 \end{equation}
Further,
 \ $\lim_{z\downarrow0} \bigl(\frac{b}{R(z)} - \frac{1}{z}\bigr) z^{2-\alpha}
    = - \frac{\delta^\alpha}{b\alpha} \in \RR_{++}$,
 \ thus there exists \ $z_1 \in (0, \theta_0)$ \ such that
 \ $\bigl(\frac{b}{R(z)} - \frac{1}{z}\bigr) z^{2-\alpha} \leq - \frac{2\delta^\alpha}{b\alpha}$
 \ for all \ $z \in (0, z_1)$.
\ Hence
 \ $\int_0^{z_1} \bigl(\frac{b}{R(z)} - \frac{1}{z}\bigr) \dd z
    \leq - \frac{2\delta^\alpha}{b\alpha} \int_0^{z_1} z^{\alpha-2} \dd z < \infty$.
\ The function \ $(0, \theta_0) \ni z \mapsto \frac{b}{R(z)} - \frac{1}{z} \in \RR_{++}$ \ is
 continuous, thus the integral
 \ $\int_{z_1}^\lambda \bigl(\frac{b}{R(z)} - \frac{1}{z}\bigr) \dd z$ \ exists and finite for all
 \ $\lambda \in (0, \theta_0)$, \ and hence the function \ $K$ \ is well-defined.
Note that the existence and finiteness of the integral
 \ $\int_0^\lambda \bigl(\frac{b}{R(z)} - \frac{1}{z}\bigr) \dd z$ \ follows also from
 Proposition 3.14 and its proof in Li \cite{Li}, since
 \ $\int_1^\infty z \log(z) \, m(\dd z) < \infty$, \ where \ the measure \ $m$ \ takes the form
 \ $m(\dd z) = \delta^\alpha C_\alpha z^{-1-\alpha} \bbone_{\RR_{++}}(z) \, \dd z$.
\ Indeed, by partial integration,
 \begin{equation}\label{PI}
  \begin{aligned}
   &\int_1^\infty z \log(z) \, m(\dd z)
    = \delta^\alpha C_\alpha \int_1^\infty \frac{\log(z)}{z^\alpha} \, \dd z \\
   &= \delta^\alpha C_\alpha \!
      \left(\frac{1}{1-\alpha} \lim_{z\to\infty} \frac{\log(z)}{z^{\alpha-1}}
            - \int_1^\infty \frac{z^{-\alpha}}{1-\alpha} \, \dd z\right)
    = \frac{\delta^\alpha C_\alpha}{(1-\alpha)^2}
    < \infty .
 \end{aligned}
 \end{equation}
The function \ $K$ \ is strictly increasing on \ $(0, \theta_0)$, \ since
 \ $\frac{b}{R(z)} - \frac{1}{z} \in \RR_{++}$ \ for all \ $z \in (0, \theta_0)$, \ and we have
 \ $\lim_{\lambda\downarrow0} K(\lambda) = 0$ \ and
 \ $\lim_{\lambda\uparrow\theta_0} K(\lambda) = + \infty$, \ yielding that the range of the
 function \ $K$ \ is \ $\RR_{++}$, \ hence the inverse \ $K^{-1}$ \ is defined on \ $\RR_{++}$.
\ Indeed,
 \ $\lim_{z\uparrow\theta_0} \frac{R(z)}{z-\theta_0}
    = \lim_{z\uparrow\theta_0} \frac{R(z)-R(\theta_0)}{z-\theta_0} = R'(\theta_0)$.
\ We have \ $R'(\theta_0) \in \RR_{++}$, \ since \ $R(0) = 0$ \ and \ $R(\theta_0) = 0$ \ yields
 the existence of \ $\theta_1 \in (0, \theta_0)$ \ with \ $R'(\theta_1) = 0$, \ and the
 function \ $R'$ \ is strictly increasing on \ $\RR_+$.
\ Thus there exists \ $z_2 \in (0, \theta_0)$ \ such that
 \ $\frac{R(z)}{z-\theta_0} \leq 2 R'(\theta_0)$ \ for all \ $z \in (z_2, \theta_0)$.
\ Hence, by \eqref{b/R(z)-1/z}, for all \ $\lambda \in (z_2, \theta_0)$, \ we have
 \[
   \int_0^\lambda \biggl(\frac{b}{R(z)} - \frac{1}{z}\biggr) \dd z
   \geq - \frac{\frac{\sigma^2}{2}z_2+\frac{\delta^\alpha}{\alpha}z_2^{\alpha-1}}
               {2 R'(\theta_0)}
          \int_{z_2}^\lambda \frac{1}{z-\theta_0} \, \dd z
   \to + \infty
   \qquad \text{as \ $\lambda \uparrow \theta_0$.}
 \]

(i) We prove the existence of an appropriate non-negative random variable \ $V$.
\ We check that
 \begin{align*}
  \EE(Y_t \mid \cF^Y_s)
  = \EE(Y_t \mid Y_s)
  = \ee^{-b(t-s)} Y_s + a \ee^{bs} \int_s^t \ee^{-bx} \, \dd x
 \end{align*}
 for all \ $s, t \in \RR_+$ \ with \ $0 \leq s \leq t$, \ where
 \ $\cF_t^Y := \sigma(Y_s, s \in [0, t])$, \ $t \in \RR_+$.
\ The first equality follows from the Markov property of the process \ $(Y_t)_{t\in\RR_+}$.
\ The second equality is a consequence of the time-homogeneity of the Markov process \ $Y$ \ and
 the fact that
 \begin{align*}
  \EE(Y_t \mid Y_s = y_0)
  = \EE(Y_{t-s} \mid Y_0 = y_0)
  = \ee^{-b(t-s)} y_0 + a \int_0^{t-s} \ee^{-bx} \, \dd x ,
  \qquad t \in \RR_+ , \qquad y_0 \in \RR_+ ,
 \end{align*}
 following from  Proposition \ref{Pro_moments}.
Then
 \[
   \EE(\ee^{bt} Y_t \mid \cF^Y_s) = \ee^{bs} Y_s + a \ee^{b(s+t)} \int_s^t \ee^{-bx} \, \dd x
   \geq \ee^{bs} Y_s
 \]
 for all \ $s, t \in \RR_+$ \ with \ $0 \leq s \leq t$, \ consequently, the process
 \ $(\ee^{bt} Y_t)_{t\in\RR_+}$ \ is a non-negative submartingale with respect to the filtration
 \ $(\cF_t^Y)_{t\in\RR_+}$.
\ Moreover,
 \begin{align*}
  \EE(\ee^{bt} Y_t)
  = y_0 + a \ee^{bt} \int_0^t \ee^{-bx} \, \dd x
  = y_0 + a \int_0^t \ee^{bx} \, \dd x
  \leq y_0 +  a \int_0^\infty \ee^{bx} \, \dd x
  = y_0 - \frac{a}{b} < \infty , \qquad t \in \RR_+ ,
 \end{align*}
 hence, by the submartingale convergence theorem, there exists a non-negative random variable
 \ $V$ \ such that
 \begin{align}\label{lim_Y}
  \ee^{bt} Y_t \as V \qquad \text{as \ $t \to \infty$.}
 \end{align}
Further, if \ $\omega \in \Omega$ \ such that \ $\ee^{bt} Y_t(\omega) \to V(\omega)$ \ as
 \ $t \to \infty$, \ then, by the integral Toeplitz lemma (see K\"uchler and S{\o}rensen
 \cite[Lemma B.3.2]{KucSor}), we have
 \[
   \frac{1}{\int_0^t \ee^{-bu} \, \dd u} \int_0^t Y_u(\omega) \, \dd\omega
   = \frac{1}{\int_0^t \ee^{-bu} \, \dd u} \int_0^t \ee^{-bu} (\ee^{bu} Y_u(\omega)) \, \dd u
   \to V(\omega)  \qquad \text{as \ $t \to \infty$.}
 \]
Here \ $\int_0^t \ee^{-bu} \, \dd u = \frac{\ee^{-bt} - 1}{-b}$, \ $t \in \RR_+$, \ thus we
 conclude
 \begin{align}\label{lim_intY}
  \ee^{bt} \int_0^t Y_u \, \dd u
  = \frac{1-\ee^{bt}}{-b} \frac{\int_0^t Y_u \, \dd u}{\int_0^t \ee^{-bu} \, \dd u}
  \as -\frac{V}{b} \qquad
  \text{as \ $t \to \infty$.}
 \end{align}

{\sl First proof of \textup{(ii)}, \textup{(iii)} and \textup{(iv)}.}
We readily have \eqref{Laplace_supercritical_limit} for \ $u = 0$, \ since then the
 unique locally bounded solution to the differential equation \eqref{psi_DE} is
 \ $\psi_{0,0}(t) = 0$, \ $t \in \RR_+$, \ implying \ $\psi_0^* = 0$.
\ Convergence \ $\ee^{bt} Y_t \as V$ \ as \ $t \to \infty$ \ implies \ $\ee^{bt} Y_t \distr V$
 \ as \ $t \to \infty$, \ and hence, by the continuity theorem,
 \ $\lim_{t\to\infty} \EE(\exp\{u \ee^{bt} Y_t\}) = \EE(\ee^{uV})$ \ for all \ $u \in \RR_-$.
\ By Theorem \ref{Thm_Laplace}, we have
 \begin{equation}\label{Laplace_supercritical_limit2}
  \EE(\exp\{u \ee^{bt} Y_t\})
  = \exp\bigg\{y_0 \psi_{u\ee^{bt},0}(t) + a \int_0^t \psi_{u\ee^{bt},0}(s) \, \dd s\bigg\} , \qquad t \in \RR_+ ,
 \end{equation}
 thus the limit
 \[
   \lim_{t\to\infty}
    \exp\bigg\{y_0 \psi_{u\ee^{bt},0}(t) + a \int_0^t \psi_{u\ee^{bt},0}(s) \, \dd s\bigg\}
   = \EE(\ee^{uV}) \in (0, 1]
 \]
 exists.
Note that the functions \ $\psi_{u\ee^{bt},0}$, \ $u \in \RR_-$, \ do not depend on the
 values of \ $a$ \ and \ $y_0$.
\ Consequently, with \ $a = 0$ \ and with some \ $y_0 \in \RR_{++}$, \ we obtain that the limit
 \ $\lim_{t\to\infty} \exp\big\{\psi_{u\ee^{bt},0}(t)\big\}$ \ exists, and hence the limit
 \ $\lim_{t\to\infty} \psi_{u\ee^{bt},0}(t) = \psi_u^*$ \ exists as well.
Using \eqref{substitution} and \ $\psi_{u,0}(s) = - v_s(-u)$ \ for all \ $s \in \RR_+$ \ and
 \ $u \in \RR_-$ \ (see part (i) of Remark \ref{Rem_integrals}), we obtain
 \[
   a \int_0^t \psi_{u\ee^{bt},0}(s) \, \dd s
   = \int_{-u \ee^{bt}}^{-\psi_{u\ee^{bt},0}(t)} \frac{F(z)}{R(z)} \, \dd z
 \]
 for all \ $t \in \RR_+$ \ and \ $u \in \RR_{--}$ \ satisfying
 \ $-u \ee^{bt} \in (0, \theta_0)$, \ which, together with \eqref{Laplace_supercritical_limit2},
 leads to \eqref{Laplace_supercritical_limit}.

If \ $t \in \RR_+$ \ and \ $u \in \RR_{--}$ \ satisfying \ $-u \ee^{bt} \in (0, \theta_0)$,
 \ then, by the proof of part (iii) of Proposition \ref{Pro_stable_CIR},
 \ $-\psi_{u\ee^{bt},0}(t) \in (0, \theta_0)$, \ hence, by Proposition 3.3 in Li \cite{Li},
 \[
   \int_{-u\ee^{bt}}^{-\psi_{u\ee^{bt},0}(t)} \frac{1}{R(z)} \, \dd z = - t .
 \]
It yields that
 \begin{equation}\label{H}
  \int_{-u\ee^{bt}}^{-\psi_{u\ee^{bt},0}(t)} \biggl(\frac{b}{R(z)} - \frac{1}{z}\biggr) \dd z
  = - b t - \int_{-u\ee^{bt}}^{-\psi_{u\ee^{bt},0}(t)} \frac{1}{z} \, \dd z .
 \end{equation}
The right hand side of \eqref{H} is
 \begin{equation}\label{rhsH}
  \begin{aligned}
   - b t - \int_{-u\ee^{bt}}^{-\psi_{u\ee^{bt},0}(t)} \frac{1}{z} \, \dd z
   &= - b t - (\log(-\psi_{u\ee^{bt},0}(t)) - \log(-u\ee^{bt})) \\
   &= \log(-u) - \log(-\psi_{u\ee^{bt},0}(t)) .
  \end{aligned}
 \end{equation}
We have already proved that \ $-\psi_{u\ee^{bt},0}(s) \in (0, \theta_0)$ \ for all
 \ $s \in \RR_+$ \ yielding \ $- \psi_u^* \in [0, \theta_0]$.
\ Letting \ $t \to \infty$ \ in \eqref{H}, we conclude \ $- \psi_u^* \in (0, \theta_0)$.
\ Indeed, \ $\psi_u^* = 0$ \ is not possible, since then the left hand side of \eqref{H} would
 tend to 0 and the right hand side of \eqref{H} would tend to \ $+ \infty$ \ (see \eqref{rhsH}).
Moreover, \ $\psi_u^* = - \theta_0$ \ is not possible, since then the left hand side of
 \eqref{H} would tend to
 \ $\int_0^{\theta_0} \bigl(\frac{b}{R(z)} - \frac{1}{z}\bigr) \dd z = + \infty$ \ (see the
 beginning of the proof), and the right hand side of \eqref{H} would tend to
 \ $\log(-u) - \log(\theta_0)$ \ (see \eqref{rhsH}).
Thus \ $- \psi_u^* \in (0, \theta_0)$, and letting \ $t \to \infty$ \ in \eqref{H}, we obtain
 \[
   \int_0^{-\psi_u^*} \biggl(\frac{b}{R(z)} - \frac{1}{z}\biggr) \dd z
   = \log(-u) - \log(-\psi_u^*)
 \]
 for all \ $u \in \RR_{--}$.
\ This can be written in the form
 \[
   u = \psi_u^*
       \exp\biggl\{\int_0^{-\psi_u^*} \biggl(\frac{b}{R(z)} - \frac{1}{z}\biggr) \dd z\biggr\} .
 \]
Consequently, \ $u = - K(-\psi_u^*)$ \ for all \ $u \in \RR_{--}$.
 \ The function \ $K$ \ is strictly increasing on \ $(0, \theta_0)$, \ see the beginning of the
 proof, hence we conclude \ $\psi_u^* = - K^{-1}(-u)$ \ for all \ $u \in \RR_{--}$.

Next, we check that if, in addition, \ $a \in \RR_{++}$, \ then \ $\PP(V \in \RR_{++}) = 1$.
\ The monotone convergence theorem yields
 \ $\EE(\ee^{uV}) \downarrow \EE(\bbone_{\{V=0\}}) = \PP(V = 0)$ \ as \ $u \to -\infty$.
\ We have \ $\lim_{u\to-\infty} \psi_u^* = - \lim_{u\to-\infty} K^{-1}(-u) = - \theta_0$,
 \ since \ $\lim_{\lambda\uparrow\theta_0} K(\lambda) = + \infty$, \ see the beginning of the
 proof.
Consequently, by \eqref{Laplace_supercritical_limit}, we obtain
 \[
   \lim_{u\to-\infty} \EE(\ee^{uV})
   = \exp\bigg\{- y_0 \theta_0  + \int_0^{\theta_0} \frac{F(z)}{R(z)} \, \dd z\bigg\}
   = 0
 \]
 and hence \ $\PP(V \in \RR_{++}) = 1$, \ since
 \ $\int_0^{\theta_0} \frac{F(z)}{R(z)} \, \dd z = - \infty$.
\ Indeed, as at the beginning of the proof, there exists \ $z_2 \in (0, \theta_0)$ \ such that
 \ $\frac{R(z)}{z-\theta_0} \leq 2 R'(\theta_0)$ \ for all \ $z \in (z_2, \theta_0)$ \ with
 \ $R'(\theta_0) \in \RR_{++}$, \ hence
 \ $\int_0^{\theta_0} \frac{F(z)}{R(z)} \, \dd z
    \leq \int_{z_2}^{\theta_0} \frac{F(z)}{R(z)} \, \dd z
    \leq \frac{az_2}{2R'(\theta_0)} \int_{z_2}^{\theta_0} \frac{1}{z-\theta_0} \, \dd z
    = -\infty$.

{\sl  Second proof of \textup{(ii)}, \textup{(iii)} and \textup{(iv)}.}
The idea of this proof is due to Cl\'ement Foucart.
First, we need to introduce some notations based on Li \cite{Li}.
For all \ $t \in \RR_+$, \ let \ $\overline{v}_t := \lim_{\lambda\to\infty} v_t(\lambda)$,
 \ which exists in \ $(0, \infty]$, \ see Li \cite[Theorem 3.5]{Li}, where \ $v_t(\lambda)$,
 \ $t, \lambda \in \RR_+$, \ is given in (iii) of Proposition \ref{Pro_stable_CIR}.
Let \ $\overline v := \lim_{t\to\infty} \overline v_t$, \ which exists in \ $\RR_+$, \ and it is
 known that \ $\overline v$ \ is the largest root of the equation \ $R(z) = 0$, \ $z \in \RR_+$,
 \ see Li \cite[Theorem 3.8]{Li}.
Indeed, Condition 3.6 in Li \cite{Li} holds in our case, since
 \[
   R(z) \geq b z + \frac{\delta^\alpha}{\alpha} z^\alpha
        \geq \frac{\delta^\alpha}{2\alpha}z^\alpha
        > 0
   \quad \text{for} \quad
   z > \left(\frac{-2b\alpha}{\delta^\alpha}\right)^{\frac{1}{\alpha-1}}
     =: \ttheta > 0 ,
 \]
 and
 \[
   \int_{\ttheta}^\infty \frac{1}{R(z)} \, \dd z
   < \frac{2\alpha}{\delta^\alpha} \int_{\ttheta}^\infty z^{-\alpha} \, \dd z
   = \frac{2\alpha\ttheta^{1-\alpha}}{\delta^\alpha(\alpha-1)}
   < \infty .
 \]
Further, by Li \cite[page 63]{Li}, \ $\overline v = \theta_0$ \ and
 \ $\theta_0 \leq \big(\frac{-2b\alpha}{\delta^\alpha}\big)^{\frac{1}{\alpha-1}}$.
\ Since \ $R'(0) = b <0$, \ $R''(z) = \sigma^2 + \delta^\alpha (\alpha-1)z^{\alpha-2} \geq 0$,
 \ $z \in \RR_{++}$, \ and \ $\lim_{z\to\infty} R(z) = \infty$, \ we get there is a positive
 root of \ $R$ \ yielding that \ $\overline v = \theta_0 > 0$.
\ For all \ $t \in \RR_+$, \ let \ $[0, \overline{v}_t) \ni q \mapsto \eta_t(q)$ \ be the
 inverse function of \ $\RR_+ \ni \lambda \mapsto v_t(\lambda)$, \ which exists and is strictly
 monotone increasing, due to the fact that \ $\RR_+ \ni \lambda \mapsto v_t(\lambda)$ \ is
 strictly monotone increasing, see Li \cite[Proposition 3.1]{Li}.
By Li \cite[Proposition 3.14]{Li}, we have
 \begin{align}\label{eta_K}
  \lim_{t\to\infty} \frac{\eta_t(\lambda)}{\ee^{bt}} = K(\lambda) \in \RR_{++} ,
  \qquad \lambda \in (0, \theta_0) .
 \end{align}
Indeed, \ $b \in \RR_{--}$, \ and by \eqref{PI},
 \ $\int_1^\infty z \log(z) \, m(\dd z) < \infty$, \ where the measure \ $m$ \ takes the form
 \ $m(\dd z) = \delta^\alpha C_\alpha z^{-1-\alpha} \bbone_{\RR_{++}}(z) \, \dd z$.
\ The form of the limit \ $K(\lambda)$ \ is
 \[
   K(\lambda)
   = \lambda
     \exp\biggl\{\int_0^\lambda \biggl(\frac{b}{R(z)} - \frac{1}{z}\biggr) \dd z\biggr\} ,
   \qquad \lambda \in (0, \theta_0) ,
 \]
 which follows by the proof of Proposition 3.14 in Li \cite{Li}.
Using \eqref{lim_Y} and \eqref{eta_K}, we have for all \ $\lambda \in (0, \theta_0)$,
 \begin{align}\label{help_eta_t}
  \eta_t(\lambda) Y_t \as K(\lambda) V =: U_\lambda \qquad \text{as \ $t \to \infty$.}
 \end{align}
Using the same ideas as in the proofs of Theorems 3.13 and 3.15 in Li \cite{Li}, we show that
 \begin{align}\label{help_Laplace_W_lambda}
  \EE(\ee^{u U_\lambda})
  = \exp\bigg\{-y_0 f(-u, \lambda) + \int_0^{f(-u,\lambda)} \frac{F(z)}{R(z)} \, \dd z \bigg\} ,
  \qquad u \in \RR_- , \quad \lambda \in (0, \theta_0) ,
 \end{align}
 where \ $f(-u, \lambda) := \lim_{t\to\infty} v_t(-u \eta_t(\lambda)) \in [0, \theta_0)$ \ and,
 in case of \ $u \in \RR_{--}$, \ $f(-u,\lambda)\in(0,\theta_0)$ \ and it satisfies
 \begin{align}\label{help_W_lambda_2}
  \int_\lambda^{f(-u,\lambda)} \frac{1}{R(z)} \, \dd z = \frac{\log(-u)}{b} .
 \end{align}
The case \ $u = 0$ \ is trivial, since in this case
 \ $f(-u, \lambda) = f(0, \lambda) = \lim_{t\to\infty} v_t(0) = 0$, \ because of \ $v_t(0) = 0$,
 \ $t \in \RR_+$ \ for \ $u = 0$.
\ So we can assume that \ $u \in \RR_{--}$.

Note that the integral \ $\int_0^{f(-u,\lambda)} \frac{F(z)}{R(z)} \, \dd z$ \ is well-defined.
First we check that
 \begin{equation}\label{integrals}
  \int_\lambda^0 \frac{1}{R(z)} \, \dd z = + \infty , \qquad
  \int_\lambda^{\theta_0} \frac{1}{R(z)} \, \dd z = - \infty , \qquad
  \lambda \in (0, \theta_0) .
 \end{equation}
Indeed, we have \ $\lim_{z\downarrow0} \frac{R(z)}{z} = b \in \RR_{--}$, \ thus there exists
 \ $z_3 \in (0, \lambda)$ \ such that \ $\frac{R(z)}{z} \leq \frac{b}{2}$ \ for all
 \ $z \in (0, z_3)$.
\ Hence
 \ $\int_\lambda^0 \frac{1}{R(z)} \, \dd z \geq \int_{z_3}^0 \frac{2}{bz} \, \dd z = + \infty$.
\ Moreover, as at the beginning of the proof, there exists \ $z_2 \in (\lambda, \theta_0)$ \ such
 that \ $\frac{R(z)}{z-\theta_0} \leq 2 R'(\theta_0)$ \ for all \ $z \in (z_2, \theta_0)$ \ with
 \ $R'(\theta_0) \in \RR_{++}$.
\ Hence
 \ $\int_\lambda^{\theta_0} \frac{1}{R(z)} \, \dd z
    \leq \frac{1}{2R'(\theta_0)} \int_{z_2}^{\theta_0} \frac{1}{z-\theta_0} \, \dd z
    = - \infty$.
\ Consequently, \ $f(-u, \lambda)$ \ can not be a root of the equation \ $R(z) = 0$,
 \ $z \in \RR_+$, \ i.e., \ $f(-u, \lambda) \notin \{0, \theta_0\}$, \ since otherwise, by
 \eqref{integrals}, the left hand side of
 \eqref{help_W_lambda_2} would be \ $\pm \infty$, \ but the right hand side of
 \eqref{help_W_lambda_2} is a real number.
Hence, using the same argument as in the proof of \eqref{substitution}, we have
 \ $f(-u, \lambda) \in (0, \theta_0)$.
\ The integrand \ $\frac{F}{R}$ \ is continuous on \ $[0, f(-u, \lambda)]$, \ since in
 case of \ $a = 0$ \ the integrand is zero, and in case of \ $a \in \RR_{++}$, \ by the
 definition of \ $\theta_0$, \ we have \ $R(z) < 0$ \ for all \ $z \in (0, \theta_0)$ \ and
 hence for all \ $(0, f(-u, \lambda)]$, \ and
 \ $\lim_{z\downarrow0} \frac{F(z)}{R(z)} = \frac{a}{b} \in \RR_{--}$.

Next, by \eqref{help_eta_t}, we have \ $\eta_t(\lambda) Y_t \distr U_\lambda$ \ as
 \ $t \to \infty$ \ for all \ $\lambda \in (0, \theta_0)$, \ and, by continuity theorem and
 \eqref{Laplace_Y_t2}, for all \ $u \in \RR_-$ \ and \ $\lambda \in (0, \theta_0)$, \ we get
 \begin{align}\label{help_W_lambda}
  \EE(\ee^{u U_\lambda})
  = \lim_{t\to\infty} \EE\left(\exp\left\{u \eta_t(\lambda) Y_t \right\}\right)
  = \lim_{t\to\infty}
     \exp\bigg\{-y_0 v_t(-u \eta_t(\lambda))
                + \int_{-u\eta_t(\lambda)}^{v_t(-u\eta_t(\lambda))}
                   \frac{F(z)}{R(z)} \, \dd z \bigg\} ,
 \end{align}
 since \ $\eta_t(\lambda)\downarrow 0$ \ as \ $t \to \infty$ \ (see Li
 \cite[proof of Proposition 3.14]{Li}), and hence for all \ $u \in \RR_{--}$, \ we have
 \ $- u \eta_t(\lambda) \in (0, \overline v) = (0, \theta_0)$ \ for sufficiently large \ $t$.
\ Recall that, by formula (3.23) in Li \cite{Li}, if \ $\eta_t(\lambda)$ \ and
 \ $- u \eta_t(\lambda)$ \ belong to \ $(0, \theta_0) = (0, \overline v)$, \ where
 \ $\lambda \in (0, \overline v)$, \ then
 \[
   \int_\lambda^{v_t(-u\eta_t(\lambda))} \frac{1}{R(z)} \, \dd z
   = \int_{\eta_t(\lambda)}^{-u\eta_t(\lambda)} \frac{1}{R(z)} \, \dd z .
 \]
Since \ $\eta_t(\lambda)\downarrow 0$ \ as \ $t \to \infty$, \ for all \ $u \in \RR_{--}$, \ we
 have \ $\eta_t(\lambda), - u \eta_t(\lambda) \in (0, \overline v)$ \ for sufficiently large \ $t$,
 \ and hence for all \ $u \in \RR_{--}$,
 \begin{align*}
  \lim_{t\to\infty} \int_{\lambda}^{v_t(-u\eta_t(\lambda))} \frac{1}{R(z)} \, \dd z
  = \lim_{t\to\infty} \int_{\eta_t(\lambda)}^{-u\eta_t(\lambda)} \frac{1}{R(z)} \, \dd z
  = \lim_{t\to\infty} \int_{\eta_t(\lambda)}^{-u\eta_t(\lambda)} \frac{1}{bz} \, \dd z
  = \lim_{t\to\infty} \frac{\log(-u)}{b}
  = \frac{\log(-u)}{b} ,
 \end{align*}
 where we used that \ $\lim_{z\downarrow 0} \frac{R(z)}{bz} = 1$.
\ The function
 \ $(0, \overline v) \ni x \mapsto \int_{\lambda}^x \frac{1}{R(z)} \, \dd z =: G(x)$ \ is
 continuous and strictly decreasing, hence its inverse is also continuous and strictly decreasing,
 implying that for all \ $u \in \RR_{--}$ \ and \ $\lambda \in (0, \overline v)$, \ the
 limit
 \[
   f(-u, \lambda)
   = \lim_{t\to\infty} v_t(- u \eta_t(\lambda))
   = \lim_{t\to\infty} G^{-1}(G(v_t(- u \eta_t(\lambda))))
   = G^{-1}\biggl(\frac{\log(-u)}{b}\biggr)
 \]
 exists and it satisfies \eqref{help_W_lambda_2}, where \ $G^{-1}$ \ denotes the inverse of
 \ $G$, \ since, by \eqref{integrals}, the range of \ $G$ \ is \ $\RR$.
\ Consequently, using the continuity of the integral upper limit function, the fact that
 \ $\eta_t(\lambda) \downarrow 0$ \ as \ $t \to \infty$, \ and \eqref{help_W_lambda}, we have
 \eqref{help_Laplace_W_lambda}, as desired.
Using that \ $K(\lambda) \in \RR_{++}$ \ and \ $V = U_\lambda/K(\lambda)$, \ we have
 \ $\EE(\ee^{u V}) = \EE(\ee^{u U_\lambda/ K(\lambda)})$, \ $u \in \RR_{-}$, \ and then
 \eqref{help_Laplace_W_lambda} yields \eqref{Laplace_supercritical_limit}.

We point out that, in the second proof of \eqref{Laplace_supercritical_limit}, we were
 not able to use directly \eqref{lim_Y}, and that's why the usage of \ $\eta_t(\lambda)$ \ in
 the argument above is really essential for us.

Next, we check that if, in addition, \ $a \in \RR_{++}$, \ then \ $\PP(V \in \RR_{++}) = 1$.
\ By the law of total probability,
 \ $\EE(\ee^{u U_\lambda})
    = 1 \cdot \PP(\ee^{-U_\lambda} = 1)
      + \EE(\ee^{u U_\lambda} \mid \ee^{-U_\lambda} \ne 1) \PP(\ee^{-U_\lambda} \ne 1)$ ,
 \ $u \in \RR_{-}$, \ and hence, by the dominated convergence theorem and
 \eqref{help_Laplace_W_lambda},
 \begin{align*}
  \PP(\ee^{-U_\lambda} = 1)
  &= \lim_{u\to-\infty} \EE(\ee^{u U_\lambda})
   = \lim_{u\to-\infty}
      \exp\bigg\{- y_0 f(-u, \lambda)
                 + \int_0^{f(-u,\lambda)} \frac{F(z)}{R(z)} \, \dd z\bigg\} \\
  &= \exp\bigg\{- y_0 \theta_0 + \int_0^{\theta_0} \frac{F(z)}{R(z)} \, \dd z\bigg\} ,
 \end{align*}
 where we used that, by \eqref{integrals},
 \ $\lim_{u\to-\infty} f(-u, \lambda) = \lim_{u\to-\infty} G^{-1}\bigl(\frac{\log(-u)}{b}\bigr)
    = \lim_{y\to-\infty} G^{-1}(y) = \theta_0 = \overline v$ \ (see also Li
 \cite[proof of Theorem 3.15]{Li}).
In case of \ $a \in \RR_{++}$, \ we have
 \ $\int_0^{\theta_0} \frac{F(z)}{R(z)} \, \dd z = - \infty$ \ (see the end of the first proof of
 (iv)), \ and hence \ $\PP(\ee^{-U_\lambda} = 1) = 0$.
\ This yields that, in case of \ $a \in \RR_{++}$, \ $\PP(U_{\lambda} = 0) = 0$ \ and, since
 \ $K(\lambda) \in \RR_{++}$, \ we have \ $\PP(V = 0) = 0$, \ i.e., \ $\PP(V \in \RR_{++}) = 1$.
\proofend

In the next remark we specialize Theorem \ref{Thm_supercritical_convergence} to the case
 \ $\sigma = 0$.

\begin{Rem}\label{Thm_supercritical_convergence_sigma=0}
Under the conditions of Theorem \ref{Thm_supercritical_convergence}, in case of \ $\sigma = 0$,
 \ using part (ii) of Theorem \ref{Thm_supercritical_convergence}, we derive the Laplace transform
 of \ $V$, \ which results in an explicit expression.
Recall that the function \ $\psi_{u\ee^{bt},0} : \RR_+ \to \RR_-$ \ is the unique locally bounded
 solution to the differential equation
 \begin{align}\label{DE_Bernoulli}
  \psi_{u\ee^{bt},0}'(s)
  = \frac{\delta^\alpha}{\alpha} (-\psi_{u\ee^{bt},0}(s))^\alpha - b \psi_{u\ee^{bt},0}(s) ,
  \qquad s \in \RR_+ , \qquad \psi_{u\ee^{bt},0}(0) = u \ee^{bt} .
 \end{align}
We have to determine the limit \ $\lim_{t\to\infty} \psi_{u\ee^{bt},0}(t)$.
\ If \ $u = 0$, \ then the unique locally bounded solution of the differential equation
 \eqref{DE_Bernoulli} is \ $\psi_{u\ee^{bt},0}(s) = 0$, \ $s \in \RR_+$, \ and hence in this case
 \ $\lim_{t\to\infty} \psi_{u\ee^{bt},0}(t) = 0$.
\ In what follows, let us suppose that \ $u \in \RR_{--}$.
\ The unique solution of the differential equation \eqref{DE_Bernoulli} (which can be transformed
 into a Bernoulli differential equation) is
 \[
   \psi_{u\ee^{bt},0}(s)
   = - \left(\left((-u\ee^{bt})^{1-\alpha} + \frac{\delta^\alpha}{b\alpha}\right)
             \ee^{b(\alpha-1)s}
             - \frac{\delta^\alpha}{b\alpha}\right)^{\frac{1}{1-\alpha}} ,
   \qquad s \in \RR_+ ,
 \]
 and consequently, by part (ii) of Theorem \ref{Thm_supercritical_convergence},
 \[
   \psi_u^* = \lim_{t\to\infty} \psi_{u\ee^{bt},0}(t)
   = - \left((-u)^{1-\alpha} - \frac{\delta^\alpha}{b\alpha}\right)^{\frac{1}{1-\alpha}} ,
   \qquad u \in \RR_{--} ,
 \]
 and hence
 \[
   \EE(\ee^{uV})
   = \exp\biggl\{- y_0
                   \biggl((-u)^{1-\alpha}
                          - \frac{\delta^\alpha}{b\alpha}\biggr)^{\frac{1}{1-\alpha}}
                 + a \int_0^{\bigl((-u)^{1-\alpha}
                                   - \frac{\delta^\alpha}{b\alpha}\bigr)^{\frac{1}{1-\alpha}}}
                      \frac{1}{\frac{\delta^\alpha}{\alpha}z^{\alpha-1}+b} \, \dd z\biggr\} ,
   \qquad u \in \RR_{--} .
 \]
We can derive the above formula for \ $\psi_u^*$, \ $u \in \RR_{--}$, \ using part (iii) of Theorem
 \ref{Thm_supercritical_convergence} as well.
We have \ $\theta_0 = \left(\frac{-b\alpha}{\delta^\alpha}\right)^{\frac{1}{\alpha-1}}$ \ and,
 by \eqref{b/R(z)-1/z},
 \begin{align*}
  \int_0^\lambda \biggl(\frac{b}{R(z)} - \frac{1}{z}\biggr) \dd z
  &= - \int_0^\lambda
        \frac{\frac{\delta^\alpha}{\alpha}z^{\alpha-2}}
             {\frac{\delta^\alpha}{\alpha}z^{\alpha-1}+b}
        \, \dd z
   = - \frac{1}{\alpha-1}
       \biggl(\log\biggl(\frac{\delta^\alpha}{\alpha} \lambda^{\alpha-1} + b\biggr)
              - \log(b)\biggr) \\
  &= \log\biggl[\biggl(\frac{\frac{\delta^\alpha}{\alpha}\lambda^{\alpha-1}+b}
                            {b}\biggr)^{-\frac{1}{\alpha-1}}\biggr]
 \end{align*}
 for all \ $\lambda \in (0, \theta_0)$.
\ Consequently,
 \begin{align*}
  K(\lambda)
  = \lambda
    \biggl(\frac{\delta^\alpha}{b\alpha} \lambda^{\alpha-1} + 1\biggr)^{-\frac{1}{\alpha-1}}
  = \biggl(\frac{\delta^\alpha}{b\alpha} + \lambda^{1-\alpha}\biggr)^{\frac{1}{1-\alpha}} ,
  \qquad \lambda \in (0, \theta_0) ,
 \end{align*}
 thus
 \[
   \psi_u^*
   = - K^{-1}(-u)
   = - \biggl((-u)^{1-\alpha} - \frac{\delta^\alpha}{b\alpha}\biggr)^{\frac{1}{1-\alpha}} ,
   \qquad u \in \RR_{--} .
 \]
\proofend
\end{Rem}

\begin{Thm}\label{Thm_MLE_supercritical}
Let \ $a \in \RR_{++}$, \ $b \in \RR_{--}$, \ $\sigma, \delta \in \RR_{++}$, \ and \ $\alpha\in(1,2)$.
\ Let \ $(Y_t)_{t\in\RR_+}$ \ be the unique strong solution of the SDE \eqref{stable_CIR}
 satisfying \ $\PP(Y_0 = y_0) = 1$ \ with some \ $y_0 \in \RR_+$.
\ Then the MLE of \ $b$ \ is strongly consistent and asymptotically mixed normal,
 i.e., \ $\hb_T \as b$ \ as \ $T \to \infty$, \ and
 \begin{align*}
  \ee^{-bT/2} (\hb_T  - b)
  \distr \sigma Z \left(-\frac{V}{b}\right)^{-1/2}
  \qquad \text{as \ $T \to \infty$,}
 \end{align*}
 where \ $V$ \ is a positive random variable having Laplace transform given in
 \eqref{Laplace_supercritical_limit}, and \ $Z$ \ is a standard normally distributed random
 variable, independent of \ $V$.

With a random scaling, we have
 \[
   \frac{1}{\sigma} \biggl(\int_0^T Y_s \, \dd s\biggr)^{1/2} (\hb_T  - b)
   \distr \cN(0, 1) \qquad \text{as \ $T \to \infty$.}
 \]
\end{Thm}

\noindent{\bf Proof.}
By Proposition \ref{LEMMA_MLEb_exist}, there exists a unique MLE \ $\hb_T$ \ of \ $b$ \ for
 all \ $T \in \RR_{++}$ \ which takes the form given in \eqref{MLEb}.
By Theorem \ref{Thm_supercritical_convergence},
 \ $\ee^{bt} \int_0^t Y_s \, \dd s \as -\frac{V}{b}$ \ as \ $t \to \infty$, \ where
 \ $\PP(V \in \RR_{++}) = 1$, \ and hence
 \[
   \int_0^t Y_s \, \dd s
   = \ee^{-bt} \ee^{bt} \int_0^t Y_s \, \dd s
   \as \infty\cdot \left(-\frac{V}{b}\right) = \infty
   \qquad \text{as \ $t\to\infty$.}
 \]
Since the quadratic variation process of the square integrable martingale
 \ $\bigl(\int_0^t \sqrt{Y_s} \, \dd W_s\bigr)_{t\in\RR_+}$ \ takes the form
 \ $\bigl(\int_0^t Y_s \, \dd s\bigr)_{t\in\RR_+}$, \ using \eqref{MLEb-} and Theorem
 \ref{DDS_stoch_int}, we have \ $\hb_T \as b$ \ as \ $T \to \infty$.
Further, by \eqref{MLEb-},
 \begin{align*}
  \ee^{-bT/2}(\hb_T  - b)
   = - \sigma \frac{\ee^{bT/2}\int_0^T\sqrt{Y_s}\,\dd W_s}{\ee^{bT}\int_0^T Y_s\,\dd s} ,
  \qquad T \in \RR_{++} .
 \end{align*}
Again, by Theorem \ref{Thm_supercritical_convergence},
 \ $\ee^{bT} \int_0^T Y_s \, \dd s \as -\frac{V}{b}$ \ as \ $T \to \infty$, \ and using
 Theorem \ref{THM_Zanten} with \ $\eta := \left(-\frac{V}{b}\right)^{1/2}$ \ and
 \ $v = -\frac{V}{b}$ \ we have
 \begin{align}\label{help_supercritical}
  \left(\ee^{bT/2} \int_0^T \sqrt{Y_s} \, \dd W_s, -\frac{V}{b}\right)
  \distr \left(\left(-\frac{V}{b}\right)^{1/2} Z, -\frac{V}{b}\right)
  \qquad \text{as \ $T \to \infty$.}
 \end{align}
Consequently, by the continuous mapping theorem,
 \ $\ee^{-bT/2} (\hb_T - b) \distr -\sigma Z \left(\frac{-V}{b}\right)^{-1/2}$ \ as
 \ $T \to \infty$, \ yielding the first assertion.

Applying again \eqref{help_supercritical} and the continuous mapping theorem, we obtain
 \begin{align*}
  \frac{1}{\sigma} \biggl(\int_0^T Y_s \, \dd s\biggr)^{1/2} (\hb_T  - b)
  &= - \biggl(\ee^{bT} \int_0^T Y_s \, \dd s\biggr)^{-1/2}
       \ee^{bT/2} \int_0^T \sqrt{Y_s} \, \dd W_s \\
  &\distr - \left(-\frac{V}{b}\right)^{-1/2} \left(-\frac{V}{b}\right)^{1/2} Z
   = -Z \distre \cN(0, 1)
  \qquad \text{as \ $T\to\infty$,}
 \end{align*}
 as desired.
\proofend

\begin{Rem}
Under the conditions of Theorem \ref{Thm_supercritical_convergence}, in case of
 \ $a \in \RR_{++}$, \ $- \log\bigl(\frac{Y_{t+1}}{Y_t}\bigr)$, \ $t \in \RR_{++}$, \ and
 \ $-\frac{Y_t}{\int_0^t Y_s\,\dd s}$, \ $t \in \RR_{++}$, \ are strongly consistent estimators
 of \ $b$ \ as well.
Indeed, by Theorem \ref{Thm_supercritical_convergence}, using that \ $\PP(V \in \RR_{++}) = 1$
 \ and \ $\PP(\int_0^t Y_s \, \dd s \in\RR_{++}) = 1$, \ $t \in \RR_{++}$ \ (see (v) of
 Proposition \ref{Pro_stable_CIR}), in case of \ $a \in \RR_{++}$,
 \[
   -\log\left(\frac{Y_{t+1}}{Y_t}\right)
   = -\log\left(\ee^{-b} \frac{\ee^{b(t+1)}Y_{t+1}}{\ee^{bt}Y_t}\right)
   \as -\log\left(\ee^{-b} \frac{V}{V} \right)
   = b \qquad \text{as \ $t \to \infty$,}
 \]
 and
 \[
   -\frac{Y_t}{\int_0^t Y_s\,\dd s}
    = - \frac{\ee^{bt}Y_t}{\ee^{bt}\int_0^t Y_s\,\dd s}
    \as -\frac{V}{-\frac{V}{b}}
    = b
    \qquad \text{as \ $t\to\infty$.}
 \]
\proofend
\end{Rem}

\vspace*{3mm}

\appendix

\vspace*{5mm}

\noindent{\bf\Large Appendices}

\section{Likelihood-ratio process}\label{App_LR}

Based on Jacod and Shiryaev \cite{JSh}, see also Jacod and M\'emin \cite{JM}, S{\o}rensen
 \cite{SorM} and Luschgy \cite{Lus}, we recall certain sufficient conditions for the absolute
 continuity of probability measures induced by semimartingales together with a representation of
 the corresponding Radon--Nikodym derivative (likelihood-ratio process).

Let \ $D(\RR_+, \RR^d)$ \ denote the space of \ $\RR^d$-valued c\`adl\`ag functions defined on
 \ $\RR_+$.
\ Let \ $(\eta_t)_{t\in\RR_+}$ \ denote the canonical process \ $\eta_t(\omega) := \omega(t)$,
 \ $\omega \in D(\RR_+, \RR^d)$, \ $t \in \RR_+$.
\ Put \ $\cF_t^\eta := \sigma(\eta_s, s \in [0, t])$, \ $t \in \RR_+$, \ and
 \[
   \cD_t(\RR_+, \RR^d)
   := \bigcap_{\vare \in \RR_{++}} \cF_{t+\vare}^\eta , \quad
   t \in \RR_+ , \qquad
   \cD(\RR_+, \RR^d)
   := \sigma\Biggl(\,\bigcup_{t \in \RR_+} \cF_t^\eta\Biggr) .
 \]
Let \ $\Psi \subset \RR^k$ \ be an arbitrary non-empty set, and let \ $\PP_\bpsi$,
 \ $\bpsi \in \Psi$, \ are probability measures on the canonical space
 \ $(D(\RR_+, \RR^d), \cD(\RR_+, \RR^d))$.
\ Suppose that for each \ $\bpsi \in \Psi$, \ under \ $\PP_\bpsi$, \ the canonical process
 \ $(\eta_t)_{t\in\RR_+}$ \ is a semimartingale with semimartingale characteristics
 \ $(B^{(\bpsi)}, C, \nu^{(\bpsi)})$ \ associated with a fixed Borel measurable truncation
 function \ $h:\RR^d \to \RR^d$, \ see Jacod and Shiryaev
 \cite[Definition II.2.6 and Remark II.2.8]{JSh}.
Namely, \ $C_t := \langle (\eta^\cont)^{(\bpsi)} \rangle_t$, \ $t \in \RR_+$, \ where
 \ $(\langle (\eta^\cont)^{(\bpsi)} \rangle_t)_{t\in\RR_+}$ \ denotes the (predictable)
 quadratic variation process (with values in \ $\RR^{d\times d}$) \ of the continuous martingale
 part \ $(\eta^\cont)^{(\bpsi)}$ \ of \ $\eta$ \ under \ $\PP_\bpsi$, \ $\nu^{(\bpsi)}$ \ is the
 compensator of the integer-valued random measure \ $\mu^\eta$ \ on \ $\RR_+ \times \RR^d$
 \ associated with the jumps of \ $\eta$ \ under \ $\PP_\bpsi$ \ given by
 \[
   \mu^\eta(\omega, \dd t, \dd\bx)
   := \sum_{s\in\RR_+}
       \bbone_{\{\Delta\eta_s(\omega)\ne\bzero\}}
       \vare_{(s,\Delta\eta_s(\omega))}(\dd t, \dd\bx) , \qquad
   \omega \in D(\RR_+, \RR^d) ,
 \]
 where \ $\vare_{(t,\bx)}$ \ denotes the Dirac measure at the point
 \ $(t, \bx) \in \RR_+ \times \RR^d$, \ and \ $\Delta\eta_t := \eta_t - \eta_{t-}$,
 \ $t \in \RR_{++}$, \ $\Delta\eta_0 := \bzero$, \ and \ $B^{(\bpsi)}$ \ is the predictable
 process (with values in \ $\RR^d$ \ having finite variation over each finite interval
 \ $[0, t]$, \ $t \in \RR_+$) \ appearing in the canonical decomposition
 \[
   \teta_t = \eta_0 + M^{(\bpsi)}_t + B^{(\bpsi)}_t , \qquad t \in \RR_+ ,
 \]
 of the special semimartingale \ $(\teta_t)_{t\in\RR_+}$ \ under \ $\PP_\bpsi$ \ given by
 \[
   \teta_t := \eta_t - \sum_{s\in[0,t]} (\eta_s - h(\Delta\eta_s)) , \qquad t \in \RR_+ ,
 \]
 where \ $(M^{(\bpsi)}_t)_{t\in\RR_+}$ \ is a local martingale with \ $M^{(\bpsi)}_0 = \bzero$.
\ We call the attention that, by our assumption, the process
 \ $C= \langle (\eta^\cont)^{(\bpsi)} \rangle$ \ does not depend on \ $\bpsi$, \ although
 \ $(\eta^\cont)^{(\bpsi)}$ \ might depend on \ $\bpsi$.
\ In addition, assume that \ $\PP_\bpsi(\nu^{(\bpsi)}(\{t\} \times \RR^d) = 0) = 1$ \ for every
 \ $\bpsi \in \Psi$, \ $t \in \RR_+$, \ and \ $\PP_\bpsi(\eta_0 = \bx_0) = 1$ \ with some
 \ $\bx_0 \in \RR^d$ \ for every \ $\bpsi \in \Psi$.
\ Note that we have the semimartingale representation
 \begin{align*}
  \eta_t &= \bx_0 + B^{(\bpsi)}_t + (\eta^\cont)^{(\bpsi)}_t
            + \int_0^t \int_{\RR^d} h(\bx) \, (\mu^\eta - \nu^{(\bpsi)})(\dd s, \dd\bx) \\
         &\quad
            + \int_0^t \int_{\RR^d} (\bx - h(\bx)) \, \mu^\eta(\dd s, \dd\bx) ,
  \qquad t \in \RR_+ ,
 \end{align*}
 of \ $\eta$ \ under \ $\PP_\bpsi$, \ see Jacod and Shiryaev \cite[Theorem II.2.34]{JSh}.
Moreover, for each \ $\bpsi \in \Psi$, \ let us choose a nondecreasing, continuous, adapted
 process \ $(F_t^{(\bpsi)})_{t\in\RR_+}$ \ with \ $F_0^{(\bpsi)} = 0$ \ and a predictable
 process \ $(c_t^{(\bpsi)})_{t\in\RR_+}$ \ with values in the set of all symmetric positive
 semidefinite \ $d \times d$ \ matrices such that
 \[
   C_t = \int_0^t c_s^{(\bpsi)} \, \dd F_s^{(\bpsi)}
 \]
 $\PP_\bpsi$-almost sure for every \ $t \in \RR_+$.
\ Due to the assumption \ $\PP_\bpsi(\nu^{(\bpsi)}(\{t\} \times \RR^d) = 0) = 1$ \ for every
 \ $\bpsi \in \Psi$, \ $t \in \RR_+$, \ such choices of \ $(F_t^{(\bpsi)})_{t\in\RR_+}$ \ and
 \ $(c_t^{(\bpsi)})_{t\in\RR_+}$ \ are possible, see Jacod and Shiryaev
 \cite[Proposition II.2.9 and Corollary II.1.19]{JSh}.
Let \ $\cP$ \ denote the predictable $\sigma$-algebra on \ $D(\RR_+, \RR^d) \times \RR_+$.
\ Assume also that for every \ $\bpsi, \btpsi \in \Psi$, \ there exist a
 \ $\cP \otimes \cB(\RR^d)$-measurable function
 \ $V^{(\btpsi,\bpsi)} : D(\RR_+, \RR^d) \times \RR_+ \times \RR^d \to \RR_{++}$ \ and a
 predictable \ $\RR^d$-valued process \ $\beta^{(\btpsi,\bpsi)}$ \ satisfying
 \begin{gather}
  \nu^{(\bpsi)}(\dd t, \dd\bx)
  = V^{(\btpsi,\bpsi)}(t, \bx) \nu^{(\btpsi)}(\dd t, \dd\bx) , \label{GIR1} \\
  \int_0^t \int_{\RR^d}
   \Bigl(\sqrt{V^{(\btpsi,\bpsi)}(s, \bx)} - 1\Bigr)^2 \,
   \nu^{(\btpsi)}(\dd s, \dd\bx)
  < \infty , \label{GIR2} \\
  B^{(\bpsi)}_t
  = B^{(\btpsi)}_t
    + \int_0^t c_s^{(\bpsi)} \beta^{(\btpsi,\bpsi)}_s \, \dd F_s^{(\bpsi)}
    + \int_0^t \int_{\RR^d}
       (V^{(\btpsi,\bpsi)}(s, \bx) - 1) h(\bx) \,
       \nu^{(\btpsi)}(\dd s, \dd\bx) , \label{GIR3} \\
  \int_0^t
   (\beta^{(\btpsi,\bpsi)}_s)^\top c_s^{(\bpsi)} \beta^{(\btpsi,\bpsi)}_s \, \dd F_s^{(\bpsi)}
  < \infty , \label{GIR4}
 \end{gather}
 $\PP_\bpsi$-almost sure for every \ $t \in \RR_+$.
\ Further, assume that for each \ $\bpsi \in \Psi$, \ local uniqueness holds for the martingale
 problem on the canonical space corresponding to the triplet \ $(B^{(\bpsi)}, C, \nu^{(\bpsi)})$
 \ with the given initial value \ $\bx_0$ \ with \ $\PP_\bpsi$ \ as its unique solution.
Then for each \ $T \in \RR_+$, \ $\PP_{\bpsi,T}$ \ is absolutely continuous with respect to
 \ $\PP_{\btpsi,T}$, \ where \ $\PP_{\bpsi,T} := \PP_{\bpsi\!\!}|_{\cD_T(\RR_+, \RR^d)}$
 \ denotes the restriction of \ $\PP_\bpsi$ \ to \ $\cD_T(\RR_+, \RR^d)$ \ (similarly for
 \ $\PP_{\btpsi,T}$), \ and, under \ $\PP_{\btpsi,T}$, \ the corresponding
 likelihood-ratio process takes the form
 \begin{align}\label{RN_general}
  \begin{split}
  \log \frac{\dd \PP_{\bpsi,T}}{\dd \PP_{\btpsi,T}}(\eta)
  &= \int_0^T (\beta^{(\btpsi,\bpsi)}_s)^\top \, \dd(\eta^\cont)^{(\btpsi)}_s
     - \frac{1}{2}
       \int_0^T
        (\beta^{(\btpsi,\bpsi)}_s)^\top c_s^{(\bpsi)} \beta^{(\btpsi,\bpsi)}_s
        \, \dd F_s^{(\bpsi)} \\
  &\quad
    + \int_0^T \int_{\RR^d}
       (V^{(\btpsi,\bpsi)}(s, \bx) - 1)
       \, (\mu^\eta - \nu^{(\btpsi)})(\dd s, \dd\bx) \\
  &\quad
    + \int_0^T \int_{\RR^d}
       (\log(V^{(\btpsi,\bpsi)}(s, \bx)) - V^{(\btpsi,\bpsi)}(s, \bx) + 1)
       \, \mu^\eta(\dd s, \dd\bx)
  \end{split}
 \end{align}
 for all \ $T \in \RR_{++}$, \ see Jacod and Shiryaev \cite[Theorem III.5.34]{JSh}.
A detailed proof of \eqref{RN_general} using Jacod and Shiryaev \cite{JSh} can be found in
 Barczy et al.\ \cite[Appendix A]{BarBenKebPap}.

\section{Limit theorems for continuous local martingales}
\label{App_clm}

In what follows we recall some limit theorems for continuous local martingales.
We use these limit theorems for studying the asymptotic behavior of the MLE of \ $b$.
\ First we recall a strong law of large numbers for continuous local martingales.

\begin{Thm}{\bf (Liptser and Shiryaev \cite[Lemma 17.4]{LipShiII})}
\label{DDS_stoch_int}
Let \ $\bigl( \Omega, \cF, (\cF_t)_{t\in\RR_+}, \PP \bigr)$ \ be a filtered probability space
 satisfying the usual conditions.
Let \ $(M_t)_{t\in\RR_+}$ \ be a square-integrable continuous local martingale with respect to
 the filtration \ $(\cF_t)_{t\in\RR_+}$ \ such that \ $\PP(M_0 = 0) = 1$.
\ Let \ $(\xi_t)_{t\in\RR_+}$ \ be a progressively measurable process such that
 \[
   \PP\left( \int_0^t \xi_u^2 \, \dd \langle M \rangle_u < \infty \right) = 1 ,
   \qquad t \in \RR_+ ,
 \]
 and
 \begin{align}\label{SEGED_STRONG_CONSISTENCY2}
  \int_0^t \xi_u^2 \, \dd \langle M \rangle_u \as \infty \qquad \text{as \ $t \to \infty$,}
 \end{align}
 where \ $(\langle M \rangle_t)_{t\in\RR_+}$ \ denotes the quadratic variation process of \ $M$.
\ Then
 \begin{align}\label{SEGED_STOCH_INT_SLLN}
  \frac{\int_0^t \xi_u \, \dd M_u}{\int_0^t \xi_u^2 \, \dd \langle M \rangle_u} \as 0 \qquad
  \text{as \ $t \to \infty$.}
 \end{align}
If \ $(M_t)_{t\in\RR_+}$ \ is a standard Wiener process, the progressive measurability of
 \ $(\xi_t)_{t\in\RR_+}$ \ can be relaxed to measurability and adaptedness to the filtration
 \ $(\cF_t)_{t\in\RR_+}$.
\end{Thm}

The next theorem is about the asymptotic behavior of continuous multivariate
 local martingales, see van Zanten \cite[Theorem 4.1]{Zan}.

\begin{Thm}{\bf (van Zanten \cite[Theorem 4.1]{Zan})}\label{THM_Zanten}
Let \ $\bigl(\Omega, \cF, (\cF_t)_{t\in\RR_+}, \PP\bigr)$ \ be a filtered probability space
 satisfying the usual conditions.
Let \ $(\bM_t)_{t\in\RR_+}$ \ be a $d$-dimensional square-integrable continuous local martingale
 with respect to the filtration \ $(\cF_t)_{t\in\RR_+}$ \ such that \ $\PP(\bM_0 = \bzero) = 1$.
\ Suppose that there exists a function \ $\bQ : [t_0, \infty) \to \RR^{d \times d}$ \ with
 some \ $t_0 \in \RR_+$ \ such that \ $\bQ(t)$ \ is an invertible (non-random) matrix for all
 \ $t \in [t_0, \infty)$, \ $\lim_{t\to\infty} \|\bQ(t)\| = 0$ \ and
 \[
   \bQ(t) \langle \bM \rangle_t \, \bQ(t)^\top \stoch \bfeta \bfeta^\top
   \qquad \text{as \ $t \to \infty$,}
 \]
 where \ $\bfeta$ \ is a \ $d \times d$ random matrix.
Then, for each $\RR^k$-valued random vector \ $\bv$ \ defined on \ $(\Omega, \cF, \PP)$, \ we
 have
 \[
   (\bQ(t) \bM_t, \bv) \distr (\bfeta \bZ, \bv) \qquad \text{as \ $t \to \infty$,}
 \]
 where \ $\bZ$ \ is a \ $d$-dimensional standard normally distributed random vector independent
 of \ $(\bfeta, \bv)$.
\end{Thm}

\section{Some explicit formulae in case of \ $\alpha = \frac{3}{2}$}
\label{App_alpha=3/2}

First, in the special case of \ $\alpha = \frac{3}{2}$, \ we make explicit the Laplace transform
 of the stationary distribution in the subcritical and critical cases given in Theorem \ref{Ergodicity} by
 evaluating the integral in its expression.

\begin{Ex}\label{Ex_Laplace_stat_3/2}
We calculate the Laplace transform of the unique stationary distribution \ $\pi$ \ given in
 Theorem \ref{Ergodicity} in case of \ $\alpha = \frac{3}{2}$.
\ Let \ $u \in \RR_{--}$.
\ By Theorem \ref{Ergodicity},
 \[
   \int_0^\infty \ee^{uy} \, \pi(\dd y)
   = \exp\biggl\{- a \int_0^{-u}
                      \frac{1}
                           {\frac{\sigma^2}{2}x+\frac{2\delta^{\frac{3}{2}}}{3}x^{\frac{1}{2}}
                            +b}
                      \, \dd x\biggr\} .
 \]
By substitution \ $x = y^2$,
 \[
   \int_0^{-u}
    \frac{1}{\frac{\sigma^2}{2}x+\frac{2\delta^{\frac{3}{2}}}{3}x^{\frac{1}{2}}+b} \, \dd x
   = \int_0^{(-u)^{\frac{1}{2}}}
      \frac{2y}{\frac{\sigma^2}{2}y^2+\frac{2\delta^{\frac{3}{2}}}{3}y+b} \, \dd y .
 \]
First we consider the case of \ $b \in \RR_{++}$ \ and \ $\sigma \in \RR_{++}$.
\ Then we can write
 \[
   \frac{2y}{\frac{\sigma^2}{2}y^2+\frac{2\delta^{\frac{3}{2}}}{3}y+b}
   = \frac{2}{\sigma^2}
     \cdot \frac{2y+\frac{4\delta^{\frac{3}{2}}}{3\sigma^2}}
                {y^2+\frac{4\delta^{\frac{3}{2}}}{3\sigma^2}y+\frac{2b}{\sigma^2}}
     - \frac{8\delta^{\frac{3}{2}}}{3\sigma^4}
       \cdot \frac{1}{y^2+\frac{4\delta^{\frac{3}{2}}}{3\sigma^2}y+\frac{2b}{\sigma^2}} .
 \]
We have
 \[
   \int_0^{(-u)^{\frac{1}{2}}}
    \frac{2y+\frac{4\delta^{\frac{3}{2}}}{3\sigma^2}}
         {y^2+\frac{4\delta^{\frac{3}{2}}}{3\sigma^2}y+\frac{2b}{\sigma^2}}
    \, \dd y
   = \left[\log\biggl(y^2 + \frac{4\delta^{\frac{3}{2}}}{3\sigma^2} y
                      + \frac{2b}{\sigma^2}\biggr)\right]_{y=0}^{y=(-u)^{\frac{1}{2}}}
   = \log\biggl(\frac{\sigma^2}{2b}(-u)+\frac{2\delta^{\frac{3}{2}}}{3b}(-u)^{\frac{1}{2}}
                + 1\biggr) .
 \]
Moreover, using
 \[
   y^2 + \frac{4\delta^{\frac{3}{2}}}{3\sigma^2} y + \frac{2b}{\sigma^2}
   = \biggl(y + \frac{2\delta^{\frac{3}{2}}}{3\sigma^2}\biggr)^2
     + \frac{2b}{\sigma^2} - \frac{4\delta^3}{9\sigma^4} ,
 \]
 we obtain
 \[
   \int \frac{1}{y^2+\frac{4\delta^{\frac{3}{2}}}{3\sigma^2}y+\frac{2b}{\sigma^2}} \, \dd y
   = \begin{cases}
      \frac{1}{\sqrt{\frac{2b}{\sigma^2}-\frac{4\delta^3}{9\sigma^4}}}
      \arctan\Biggl(\frac{y+\frac{2\delta^{\frac{3}{2}}}{3\sigma^2}}
                         {\sqrt{\frac{2b}{\sigma^2}-\frac{4\delta^3}{9\sigma^4}}}\Biggr)
      + C
       & \text{if \ $b \in \bigl(\frac{2\delta^3}{9\sigma^2}, \infty\bigr)$,} \\[4mm]
      - \frac{1}{y+\frac{2\delta^{\frac{3}{2}}}{3\sigma^2}} + C
       & \text{if \ $b = \frac{2\delta^3}{9\sigma^2}$,} \\[4mm]
      \frac{1}{\sqrt{\frac{4\delta^3}{9\sigma^4}-\frac{2b}{\sigma^2}}}
      \log\Biggl(\Biggl|\frac{y+\frac{2\delta^{\frac{3}{2}}}{3\sigma^2}
                              -\sqrt{\frac{4\delta^3}{9\sigma^4}-\frac{2b}{\sigma^2}}}
                             {y+\frac{2\delta^{\frac{3}{2}}}{3\sigma^2}
                              +\sqrt{\frac{4\delta^3}{9\sigma^4}-\frac{2b}{\sigma^2}}}\Biggr|
          \Biggr)
      + C
       & \text{if \ $b \in \bigl(0, \frac{2\delta^3}{9\sigma^2}\bigr)$,}
     \end{cases}
 \]
 where \ $C \in \RR$.
\ If \ $b \in \bigl(\frac{2\delta^3}{9\sigma^2}, \infty\bigr)$ \ and \ $\sigma \in \RR_{++}$
 \ then, applying the formula
 \ $\arctan(u) - \arctan(v) = \arctan\bigl(\frac{u-v}{1+uv}\bigr)$, \ $u, v \in \RR_+$, \ we get
 \begin{align*}
  \int_0^{(-u)^{\frac{1}{2}}}
   \frac{1}{y^2+\frac{4\delta^{\frac{3}{2}}}{3\sigma^2}y+\frac{2b}{\sigma^2}} \, \dd y
  &= \frac{1}{\sqrt{\frac{2b}{\sigma^2}-\frac{4\delta^3}{9\sigma^4}}}
     \left(\arctan\Biggl(\frac{(-u)^{\frac{1}{2}}+\frac{2\delta^{\frac{3}{2}}}{3\sigma^2}}
                              {\sqrt{\frac{2b}{\sigma^2}-\frac{4\delta^3}{9\sigma^4}}}\Biggr)
           - \arctan\Biggl(\frac{\frac{2\delta^{\frac{3}{2}}}{3\sigma^2}}
                                {\sqrt{\frac{2b}{\sigma^2}-\frac{4\delta^3}{9\sigma^4}}}\Biggr)
     \right) \\
  &= \frac{1}{\sqrt{\frac{2b}{\sigma^2}-\frac{4\delta^3}{9\sigma^4}}}
     \arctan\left(\frac{\sqrt{\frac{2b}{\sigma^2}-\frac{4\delta^3}{9\sigma^4}}}
                       {\frac{2b}{\sigma^2}(-u)^{-\frac{1}{2}}
                        +\frac{2\delta^{\frac{3}{2}}}{3\sigma^2}}\right) ,
 \end{align*}
 and hence
 \[
   \int_0^\infty \ee^{uy} \, \pi(\dd y)
   = \biggl(\frac{\sigma^2}{2b}(-u)+\frac{2\delta^{\frac{3}{2}}}{3b}(-u)^{\frac{1}{2}}
            + 1\biggr)^{-\frac{2a}{\sigma^2}}
     \exp\left\{\frac{\frac{8a\delta^{\frac{3}{2}}}{3\sigma^4}}
                     {\sqrt{\frac{2b}{\sigma^2}-\frac{4\delta^3}{9\sigma^4}}}
                \arctan\left(\frac{\sqrt{\frac{2b}{\sigma^2}-\frac{4\delta^3}{9\sigma^4}}}
                                  {\frac{2b}{\sigma^2}(-u)^{-\frac{1}{2}}
                                   +\frac{2\delta^{\frac{3}{2}}}{3\sigma^2}}\right)\right\} .
 \]
If \ $b = \frac{2\delta^3}{9\sigma^2}$ \ and \ $\sigma \in \RR_{++}$, \ then
 \[
   \int_0^{(-u)^{\frac{1}{2}}}
    \frac{1}{y^2+\frac{4\delta^{\frac{3}{2}}}{3\sigma^2}y+\frac{2b}{\sigma^2}} \, \dd y
   = - \frac{1}{(-u)^{\frac{1}{2}}+\frac{2\delta^{\frac{3}{2}}}{3\sigma^2}}
     + \frac{1}{\frac{2\delta^{\frac{3}{2}}}{3\sigma^2}}
   = \frac{(-u)^{\frac{1}{2}}}{(-u)^{\frac{1}{2}}+\frac{2\delta^{\frac{3}{2}}}{3\sigma^2}}
   = \frac{1}{1+\frac{2\delta^{\frac{3}{2}}}{3\sigma^2}(-u)^{-\frac{1}{2}}} ,
 \]
 and hence
 \[
   \int_0^\infty \ee^{uy} \, \pi(\dd y)
   = \biggl(\frac{9\sigma^4}{4\delta^3}(-u)
            + \frac{3\sigma^2}{\delta^{\frac{3}{2}}}(-u)^{\frac{1}{2}}
            + 1\biggr)^{-\frac{2a}{\sigma^2}}
     \exp\Biggl\{\frac{\frac{8a\delta^{\frac{3}{2}}}{3\sigma^4}}
                      {1+\frac{2\delta^{\frac{3}{2}}}{3\sigma^2}(-u)^{-\frac{1}{2}}}\Biggr\} .
 \]
If \ $b \in \bigl(0, \frac{2\delta^3}{9\sigma^2}\bigr)$ \ and \ $\sigma \in \RR_{++}$, \ then
 \begin{align*}
  &\int_0^{(-u)^{\frac{1}{2}}}
    \frac{1}{y^2+\frac{4\delta^{\frac{3}{2}}}{3\sigma^2}y+\frac{2b}{\sigma^2}} \, \dd y \\
  &= \frac{1}{\sqrt{\frac{4\delta^3}{9\sigma^4}-\frac{2b}{\sigma^2}}}
     \left(\log\left(\frac{(-u)^{\frac{1}{2}}+\frac{2\delta^{\frac{3}{2}}}{3\sigma^2}
                           -\sqrt{\frac{4\delta^3}{9\sigma^4}-\frac{2b}{\sigma^2}}}
                          {(-u)^{\frac{1}{2}}+\frac{2\delta^{\frac{3}{2}}}{3\sigma^2}
                           +\sqrt{\frac{4\delta^3}{9\sigma^4}-\frac{2b}{\sigma^2}}}\right)
           - \log\left(\frac{\frac{2\delta^{\frac{3}{2}}}{3\sigma^2}
                             -\sqrt{\frac{4\delta^3}{9\sigma^4}-\frac{2b}{\sigma^2}}}
                            {\frac{2\delta^{\frac{3}{2}}}{3\sigma^2}
                             +\sqrt{\frac{4\delta^3}{9\sigma^4}-\frac{2b}{\sigma^2}}}\right)
     \right) \\
  &= \frac{1}{\sqrt{\frac{4\delta^3}{9\sigma^4}-\frac{2b}{\sigma^2}}}
     \log\left(\frac{(-u)^{\frac{1}{2}}
                     \Bigl(\frac{2\delta^{\frac{3}{2}}}{3\sigma^2}
                           +\sqrt{\frac{4\delta^3}{9\sigma^4}-\frac{2b}{\sigma^2}}\Bigr)
                     +\frac{2b}{\sigma^2}}
                    {(-u)^{\frac{1}{2}}
                     \Bigl(\frac{2\delta^{\frac{3}{2}}}{3\sigma^2}
                           -\sqrt{\frac{4\delta^3}{9\sigma^4}-\frac{2b}{\sigma^2}}\Bigr)
                     +\frac{2b}{\sigma^2}}\right) ,
 \end{align*}
 and hence
 \[
   \int_0^\infty \ee^{uy} \, \pi(\dd y)
   = \biggl(\frac{\sigma^2}{2b}(-u)+\frac{2\delta^{\frac{3}{2}}}{3b}(-u)^{\frac{1}{2}}
            + 1\biggr)^{-\frac{2a}{\sigma^2}}
     \left(\frac{(-u)^{\frac{1}{2}}
                 \Bigl(\frac{2\delta^{\frac{3}{2}}}{3\sigma^2}
                       +\sqrt{\frac{4\delta^3}{9\sigma^4}-\frac{2b}{\sigma^2}}\Bigr)
                 +\frac{2b}{\sigma^2}}
                 {(-u)^{\frac{1}{2}}
                 \Bigl(\frac{2\delta^{\frac{3}{2}}}{3\sigma^2}
                       -\sqrt{\frac{4\delta^3}{9\sigma^4}-\frac{2b}{\sigma^2}}\Bigr)
                 +\frac{2b}{\sigma^2}}
     \right)^\frac{\frac{8a\delta^{\frac{3}{2}}}{3\sigma^4}}
                  {\sqrt{\frac{4\delta^3}{9\sigma^4}-\frac{2b}{\sigma^2}}} .
 \]
Next we consider the case of \ {$b \in \RR_{++}$ \ and} \ $\sigma = 0$.
\ Then we can write
 \[
   \frac{2y}{\frac{2\delta^{\frac{3}{2}}}{3}y+b}
   = \frac{3}{\delta^{\frac{3}{2}}}
     - \frac{\frac{9b}{2\delta^3}}{y+\frac{3b}{2\delta^{\frac{3}{2}}}} ,
 \]
 thus
 \begin{align*}
  \int_0^{(-u)^{\frac{1}{2}}} \frac{2y}{\frac{2\delta^{\frac{3}{2}}}{3}y+b} \dd y
  &= \frac{3}{\delta^{\frac{3}{2}}} (-u)^{\frac{1}{2}}
     - \frac{9b}{2\delta^3}
       \biggl(\log\biggl((-u)^{\frac{1}{2}}+\frac{3b}{2\delta^{\frac{3}{2}}}\biggr)
              - \log\biggl(\frac{3b}{2\delta^{\frac{3}{2}}}\biggr)\biggr) \\
  &= \frac{3}{\delta^{\frac{3}{2}}} (-u)^{\frac{1}{2}}
     - \frac{9b}{2\delta^3}
       \log\biggl(\frac{2\delta^{\frac{3}{2}}}{3b}(-u)^{\frac{1}{2}} + 1\biggr) ,
 \end{align*}
 and hence
 \[
   \int_0^\infty \ee^{uy} \, \pi(\dd y)
   = \exp\biggl\{-\frac{3a}{\delta^{\frac{3}{2}}}(-u)^{\frac{1}{2}}\biggr\}
     \biggl(1 + \frac{2\delta^{\frac{3}{2}}}{3b}
                (-u)^{\frac{1}{2}}\biggr)^{\frac{9ba}{2\delta^3}} .
 \]

Next we consider the case of \ $b = 0$ \ and \ $\sigma \in \RR_{++}$.
\ Let \ $u \in \RR_-$.
\ By Theorem \ref{Ergodicity},
 \[
   \int_0^\infty \ee^{uy} \, \pi(\dd y)
   = \exp\biggl\{- a \int_0^{-u}
                      \frac{1}
                           {\frac{\sigma^2}{2}x+\frac{2\delta^{\frac{3}{2}}}{3}x^{\frac{1}{2}}}
                      \, \dd x\biggr\} .
 \]
By substitution \ $x = y^2$,
 \[
   \int_0^{-u}
    \frac{1}{\frac{\sigma^2}{2}x+\frac{2\delta^{\frac{3}{2}}}{3}x^{\frac{1}{2}}} \, \dd x
   = \int_0^{(-u)^{\frac{1}{2}}}
      \frac{2y}{\frac{\sigma^2}{2}y^2+\frac{2\delta^{\frac{3}{2}}}{3}y} \, \dd y
   = \int_0^{(-u)^{\frac{1}{2}}}
      \frac{2}{\frac{\sigma^2}{2}y+\frac{2\delta^{\frac{3}{2}}}{3}} \, \dd y .
 \]
Consequently,
 \begin{align*}
  \int_0^{(-u)^{\frac{1}{2}}}
   \frac{2}{\frac{\sigma^2}{2}y+\frac{2\delta^{\frac{3}{2}}}{3}} \, \dd y
  &= \left[\frac{4}{\sigma^2}
           \log\biggl(y + \frac{4\delta^{\frac{3}{2}}}{3\sigma^2}\biggr)
     \right]_{y=0}^{y=(-u)^{\frac{1}{2}}}
   = \frac{4}{\sigma^2}
     \log\biggl((-u)^{\frac{1}{2}} + \frac{4\delta^{\frac{3}{2}}}{3\sigma^2}\biggr)
     - \frac{4}{\sigma^2} \log\biggl(\frac{4\delta^{\frac{3}{2}}}{3\sigma^2}\biggr) \\
  &= \frac{4}{\sigma^2}
     \log\biggl(\frac{3\sigma^2}{4\delta^{\frac{3}{2}}} (-u)^{\frac{1}{2}} + 1\biggr) ,
 \end{align*}
 hence
 \[
   \int_0^\infty \ee^{uy} \, \pi(\dd y)
   = \biggl(\frac{3\sigma^2}{4\delta^{\frac{3}{2}}} (-u)^{\frac{1}{2}}
            + 1\biggr)^{-\frac{4a}{\sigma^2}} .
 \]
Finally, by the proof of Theorem \ref{Ergodicity}, if \ $b = 0$, \ $\sigma = 0$ \ and
 \ $\alpha = \frac{3}{2}$, \ then
 \ $\int_0^\infty \ee^{uy} \, \pi(\dd y) = \exp\{-\frac{3a}{\delta^{3/2}} (-u)^{1/2}\}$,
 \ $u \in \RR_{-}$.
\proofend
\end{Ex}

\begin{Ex}\label{Ex_Laplace_crit_3/2}
Now we formulate a special case of Theorem \ref{Thm_Laplace} giving the Laplace transform of
 \ $Y_t$ \ in case of \ $\alpha = \frac{3}{2}$.
\ Let \ $(Y_t)_{t\in\RR_+}$ \ be the unique strong solution of the SDE \eqref{stable_CIR}
 satisfying \ $\PP(Y_0 = y_0) = 1$ \ with some \ $y_0 \in \RR_+$, \ with \ $a \in \RR_+$,
 \ $b \in \RR$, \ $\sigma \in \RR_{++}$, \ $\delta \in \RR_{++}$ \ and \ $\alpha = \frac{3}{2}$.
\ Then, by Theorem \ref{Thm_Laplace}, for all \ $u \in \RR_-$,
 \begin{align}\label{Laplace_crit_3/2}
  \EE(\ee^{uY_t}) = \exp\left\{\psi_{u,0}(t) y_0 + a \int_0^t \psi_{u,0}(s)\,\dd s \right\} ,
  \qquad t \in \RR_+ ,
 \end{align}
 where the function \ $\psi_{u,0} : \RR_+ \to \RR_-$ \ is the unique locally bounded solution to
 the differential equation
 \begin{align}\label{psi_3/2}
  \psi'_{u,0}(t)
  = \frac{\sigma^2}{2} \psi_{u,0}(t)^2
    + \frac{2\delta^{\frac{3}{2}}}{3} (-\psi_{u,0}(t))^{\frac{3}{2}}
    - b \psi_{u,0}(t) ,
  \qquad t \in \RR_+ , \qquad
  \psi_{u,0}(0) = u .
 \end{align}
In case of \ $u = 0$, \ the unique locally bounded solution of \eqref{psi_3/2} is
 \ $\psi_{0,0}(t) = 0$, \ $t \in \RR_+$.
\ Let us consider the function
 \ $g_u(t) := (-\psi_{u,0}(t))^{\frac{1}{2}} \in \RR_{++}$, \ $t \in \RR_+$.
\ Then we have \ $\psi_{u,0}(t) = - g_u(t)^2$, \ $\psi_{u,0}(t)^2 = g_u(t)^4$,
 \ $(-\psi_{u,0}(t))^{\frac{3}{2}} = g_u(t)^3$ \ and
 \ $\psi'_{u,0}(t) = - 2 g_u(t) g'_u(t)$ \ for all \ $t \in \RR_+$ \ and \ $u \in \RR_{--}$,
 \ hence \eqref{psi_3/2} yields
 \begin{align}\label{g_3/2}
  g'_u(t)
  = - \frac{\sigma^2}{4} g_u(t)^3 - \frac{\delta^{\frac{3}{2}}}{3} g_u(t)^2
    - \frac{b}{2} g_u(t) ,
  \qquad t \in \RR_+ , \qquad
  g_u(0) = (-u)^{\frac{1}{2}} .
 \end{align}
In case of \ $b \in \RR_+$, \ \eqref{g_3/2} has a constant solution if and only if \ $u = 0$,
 \ and then \ $g_0(t) = \psi_{0,0}(t) = 0$ \ for all \ $t \in \RR_+$.
\ In case of \ $b \in \RR_{--}$, \ \eqref{g_3/2} has a constant solution if and only if
 \ $u = 0$
 \ or
 \[
   u = u_0 := - \biggl(- \frac{2\delta^{\frac{3}{2}}}{3\sigma^2}
                       + \sqrt{\frac{4\delta^3}{9\sigma^4} - \frac{2b}{\sigma^2}}\biggr)^2 ,
 \]
 and then \ $g_0(t) = \psi_{0,0}(t) = 0$ \ for all \ $t \in \RR_+$ \ or
 \ $g_{u_0}(t) = (-u_0)^{{\frac{1}{2}}}$, \ and hence \ $\psi_{u_0,0}(t) = u_0$
 \ for all \ $t \in \RR_+$.
\ In the sequel, we suppose that \ $u \in \RR_{--}$, \ and in case of \ $b \in \RR_{--}$, \ in
 addition, we suppose that \ $u \ne u_0$.
\ Then, by separation of variables, we have
 \[
   \frac{1}
        {\bigl(g_u^2+\frac{4\delta^{\frac{3}{2}}}{3\sigma^2}g_u+\frac{2b}{\sigma^2}\bigr)g_u}
   \dd g_u
   = - \frac{\sigma^2}{4} \dd t .
 \]
If \ $b \ne 0$, \ then
 \begin{align*}
  \frac{1}
       {\bigl(g_u^2+\frac{4\delta^{\frac{3}{2}}}{3\sigma^2}g_u+\frac{2b}{\sigma^2}\bigr)g_u}
  &= \frac{\sigma^2}{2bg_u}
     - \frac{\frac{\sigma^2}{2b}g_u+\frac{2\delta^{\frac{3}{2}}}{3b}}
            {g_u^2+\frac{4\delta^{\frac{3}{2}}}{3\sigma^2}g_u+\frac{2b}{\sigma^2}} \\
  &= \frac{\sigma^2}{2bg_u}
     - \frac{\frac{\sigma^2}{2b}g_u+\frac{\delta^{\frac{3}{2}}}{3b}}
            {g_u^2+\frac{4\delta^{\frac{3}{2}}}{3\sigma^2}g_u+\frac{2b}{\sigma^2}}
     - \frac{\frac{\delta^{\frac{3}{2}}}{3b}}
            {g_u^2+\frac{4\delta^{\frac{3}{2}}}{3\sigma^2}g_u+\frac{2b}{\sigma^2}} ,
 \end{align*}
 and we have
 \begin{align*}
  \int
   \left(\frac{\sigma^2}{2bg_u}
         - \frac{\frac{\sigma^2}{2b}g_u+\frac{\delta^{\frac{3}{2}}}{3b}}
                {g_u^2+\frac{4\delta^{\frac{3}{2}}}{3\sigma^2}g_u+\frac{2b}{\sigma^2}}\right)
   \dd g_u
  &= \frac{\sigma^2}{2b} \log(|g_u|)
     - \frac{\sigma^2}{4b}
       \log\biggl(\biggl|g_u^2 + \frac{4\delta^{\frac{3}{2}}}{3\sigma^2} g_u
                         + \frac{2b}{\sigma^2}\biggr|\biggr)
     + C \\
  &= - \frac{\sigma^2}{4b}
       \log\biggl(\biggl|1 + \frac{4\delta^{\frac{3}{2}}}{3\sigma^2g_u}
                         + \frac{2b}{\sigma^2g_u^2}\biggr|\biggr)
     + C ,
 \end{align*}
 where \ $C \in \RR$.
\ Moreover, using
 \[
   g_u^2 + \frac{4\delta^{\frac{3}{2}}}{3\sigma^2} g_u + \frac{2b}{\sigma^2}
   = \biggl(g_u + \frac{2\delta^{\frac{3}{2}}}{3\sigma^2}\biggr)^2
     + \frac{2b}{\sigma^2} - \frac{4\delta^3}{9\sigma^4} ,
 \]
 we obtain
 \begin{align}\label{help_alpha_3/2_1}
   \int
    \frac{\frac{\delta^{\frac{3}{2}}}{3b}}
         {g_u^2+\frac{4\delta^{\frac{3}{2}}}{3\sigma^2}g_u+\frac{2b}{\sigma^2}}
    \dd g_u
   = \begin{cases}
      \frac{\frac{\delta^{\frac{3}{2}}}{3b}}
           {\sqrt{\frac{2b}{\sigma^2}-\frac{4\delta^3}{9\sigma^4}}}
      \arctan\Biggl(\frac{g_u+\frac{2\delta^{\frac{3}{2}}}{3\sigma^2}}
                         {\sqrt{\frac{2b}{\sigma^2} - \frac{4\delta^3}{9\sigma^4}}}\Biggr)
      + C
       & \text{if \ $b > \frac{2\delta^3}{9\sigma^2}$,} \\[4mm]
      - \frac{\frac{\delta^{\frac{3}{2}}}{3b}}{g_u+\frac{2\delta^{\frac{3}{2}}}{3\sigma^2}} + C
       & \text{if \ $b = \frac{2\delta^3}{9\sigma^2}$,} \\[4mm]
      \frac{\frac{\delta^{\frac{3}{2}}}{3b}}
           {\sqrt{\frac{4\delta^3}{9\sigma^4}-\frac{2b}{\sigma^2}}}
      \log\Biggl(\Biggl|\frac{g_u+\frac{2\delta^{\frac{3}{2}}}{3\sigma^2}
                              -\sqrt{\frac{4\delta^3}{9\sigma^4}-\frac{2b}{\sigma^2}}}
                             {g_u+\frac{2\delta^{\frac{3}{2}}}{3\sigma^2}
                              +\sqrt{\frac{4\delta^3}{9\sigma^4}-\frac{2b}{\sigma^2}}}\Biggr|
          \Biggr)
      + C
       & \text{if \ $b < \frac{2\delta^3}{9\sigma^2}$.}
     \end{cases}
 \end{align}
Consequently, if \ $b > \frac{2\delta^3}{9\sigma^2}$, \ then
 \[
   - \frac{\sigma^2}{4b}
     \log\biggl(1 + \frac{4\delta^{\frac{3}{2}}}{3\sigma^2g_u(t)}
                + \frac{2b}{\sigma^2g_u(t)^2}\biggr)
   - \frac{\frac{\delta^{\frac{3}{2}}}{3b}}
          {\sqrt{\frac{2b}{\sigma^2}-\frac{4\delta^3}{9\sigma^4}}}
     \arctan\left(\frac{g_u(t)+\frac{2\delta^{\frac{3}{2}}}{3\sigma^2}}
                       {\sqrt{\frac{2b}{\sigma^2} - \frac{4\delta^3}{9\sigma^4}}}\right)
   = - \frac{\sigma^2}{4} t + C , \qquad t \in \RR_+ ,
 \]
 with some \ $C \in \RR_+$.
\ Using the initial value \ $g_u(0) = (-u)^{\frac{1}{2}}$, \ we obtain
 \[
   C = - \frac{\sigma^2}{4b}
         \log\biggl(1 + \frac{4\delta^{\frac{3}{2}}}{3\sigma^2(-u)^{\frac{1}{2}}}
                    - \frac{2b}{\sigma^2u}\biggr)
       - \frac{\frac{\delta^{\frac{3}{2}}}{3b}}
              {\sqrt{\frac{2b}{\sigma^2}-\frac{4\delta^3}{9\sigma^4}}}
         \arctan\left(\frac{(-u)^{\frac{1}{2}}+\frac{2\delta^{\frac{3}{2}}}{3\sigma^2}}
                           {\sqrt{\frac{2b}{\sigma^2} - \frac{4\delta^3}{9\sigma^4}}}\right) ,
 \]
 and hence, by \ $g_u(t) = (-\psi_{u,0}(t))^{\frac{1}{2}}$, \ we conclude
 \begin{align*}
  &\frac{\sigma^2}{4b}
   \log\left(\frac{1 + \frac{4\delta^{\frac{3}{2}}}{3\sigma^2(-\psi_{u,0}(t))^{\frac{1}{2}}}
                   - \frac{2b}{\sigma^2\psi_{u,0}(t)}}
                  {1 + \frac{4\delta^{\frac{3}{2}}}{3\sigma^2(-u)^{\frac{1}{2}}}
                   - \frac{2b}{\sigma^2u}}\right) \\
  &+ \frac{\frac{\delta^{\frac{3}{2}}}{3b}}
          {\sqrt{\frac{2b}{\sigma^2}-\frac{4\delta^3}{9\sigma^4}}}
     \left(\arctan\left(\frac{(-\psi_{u,0}(t))^{\frac{1}{2}}
                              +\frac{2\delta^{\frac{3}{2}}}{3\sigma^2}}
                             {\sqrt{\frac{2b}{\sigma^2} - \frac{4\delta^3}{9\sigma^4}}}\right)
           - \arctan\left(\frac{(-u)^{\frac{1}{2}}+\frac{2\delta^{\frac{3}{2}}}{3\sigma^2}}
                               {\sqrt{\frac{2b}{\sigma^2} - \frac{4\delta^3}{9\sigma^4}}}\right)
     \right)
   = \frac{\sigma^2}{4} t .
 \end{align*}
In a similar way, if \ $b = \frac{2\delta^3}{9\sigma^2}$, \ then
 \[
   - \frac{\sigma^2}{4b}
     \log\biggl(1 + \frac{4\delta^{\frac{3}{2}}}{3\sigma^2g_u(t)}
                + \frac{2b}{\sigma^2g_u(t)^2}\biggr)
   + \frac{\frac{\delta^{\frac{3}{2}}}{3b}}{g_u(t)+\frac{2\delta^{\frac{3}{2}}}{3\sigma^2}}
   = - \frac{\sigma^2}{4} t + C ,\qquad t \in \RR_+ ,
 \]
 with some \ $C \in \RR_+$.
\ Using the initial value \ $g_u(0) = (-u)^{\frac{1}{2}}$, \ we obtain
 \[
   C = - \frac{\sigma^2}{4b}
         \log\biggl(1 + \frac{4\delta^{\frac{3}{2}}}{3\sigma^2(-u)^{\frac{1}{2}}}
                    - \frac{2b}{\sigma^2u}\biggr)
       + \frac{\frac{\delta^{\frac{3}{2}}}{3b}}
              {(-u)^{\frac{1}{2}}+\frac{2\delta^{\frac{3}{2}}}{3\sigma^2}} ,
 \]
 \ and hence, by \ $g_u(t) = (-\psi_{u,0}(t))^{\frac{1}{2}}$, \ we conclude
 \[
   \frac{\sigma^2}{4b}
   \log\left(\frac{1 + \frac{4\delta^{\frac{3}{2}}}{3\sigma^2(-\psi_{u,0}(t))^{\frac{1}{2}}}
                   - \frac{2b}{\sigma^2\psi_{u,0}(t)}}
                  {1 + \frac{4\delta^{\frac{3}{2}}}{3\sigma^2(-u)^{\frac{1}{2}}}
                   - \frac{2b}{\sigma^2u}}\right)
   - \frac{\frac{\delta^{\frac{3}{2}}}{3b}}
          {(-\psi_{u,0}(t))^{\frac{1}{2}}+\frac{2\delta^{\frac{3}{2}}}{3\sigma^2}}
   + \frac{\frac{\delta^{\frac{3}{2}}}{3b}}
          {(-u)^{\frac{1}{2}}+\frac{2\delta^{\frac{3}{2}}}{3\sigma^2}}
   = \frac{\sigma^2}{4} t .
 \]
Further, if \ $b \ne 0$ \ and \ $b < \frac{2\delta^3}{9\sigma^2}$, \ then
 \begin{align}\label{help_alpha_3/2_2}
  \begin{split}
  &\frac{\sigma^2}{2b} \log(g_u(t))
   - \frac{\sigma^2}{4b}
     \log\biggl(\biggl|g_u(t)^2 + \frac{4\delta^{\frac{3}{2}}}{3\sigma^2} g_u(t)
                       + \frac{2b}{\sigma^2}\biggr|\biggr) \\
  &- \frac{\frac{\delta^{\frac{3}{2}}}{3b}}
          {\sqrt{\frac{4\delta^3}{9\sigma^4}-\frac{2b}{\sigma^2}}}
     \log\left(\left|\frac{g_u(t)+\frac{2\delta^{\frac{3}{2}}}{3\sigma^2}
                           -\sqrt{\frac{4\delta^3}{9\sigma^4}-\frac{2b}{\sigma^2}}}
                          {g_u(t)+\frac{2\delta^{\frac{3}{2}}}{3\sigma^2}
                           +\sqrt{\frac{4\delta^3}{9\sigma^4}-\frac{2b}{\sigma^2}}}\right|
         \right)
   = - \frac{\sigma^2}{4} t + C ,\qquad t \in \RR_+ ,
  \end{split}
 \end{align}
 with some \ $C \in \RR$.
\ Using the initial value \ $g_u(0) = (-u)^{\frac{1}{2}}$, \ we obtain
 \begin{align*}
  C &= \frac{\sigma^2}{2b} \log((-u)^{\frac{1}{2}})
       - \frac{\sigma^2}{4b}
         \log\biggl(\biggl|-u + \frac{4\delta^{\frac{3}{2}}}{3\sigma^2} (-u)^{\frac{1}{2}}
                           + \frac{2b}{\sigma^2}\biggr|\biggr) \\
    &\quad
       - \frac{\frac{\delta^{\frac{3}{2}}}{3b}}
              {\sqrt{\frac{4\delta^3}{9\sigma^4}-\frac{2b}{\sigma^2}}}
         \log\left(\left|\frac{(-u)^{\frac{1}{2}}+\frac{2\delta^{\frac{3}{2}}}{3\sigma^2}
                               -\sqrt{\frac{4\delta^3}{9\sigma^4}-\frac{2b}{\sigma^2}}}
                              {(-u)^{\frac{1}{2}}+\frac{2\delta^{\frac{3}{2}}}{3\sigma^2}
                               +\sqrt{\frac{4\delta^3}{9\sigma^4}-\frac{2b}{\sigma^2}}}\right|
             \right) ,
 \end{align*}
 and hence, by \ $g_u(t) = (-\psi_{u,0}(t))^{\frac{1}{2}}$, \ we conclude
 \begin{equation}\label{psi_supercrit}
  \begin{aligned}
   &- \frac{\sigma^2}{2b} \log((-\psi_{u,0}(t))^{\frac{1}{2}})
    + \frac{\sigma^2}{4b}
      \log\biggl(\biggl|- \psi_{u,0}(t)
                        + \frac{4\delta^{\frac{3}{2}}}{3\sigma^2} (-\psi_{u,0}(t))^{\frac{1}{2}}
                        + \frac{2b}{\sigma^2}\biggr|\biggr) \\
   &+ \frac{\sigma^2}{2b} \log((-u)^{\frac{1}{2}})
        - \frac{\sigma^2}{4b}
          \log\biggl(\biggl|-u + \frac{4\delta^{\frac{3}{2}}}{3\sigma^2} (-u)^{\frac{1}{2}}
                            + \frac{2b}{\sigma^2}\biggr|\biggr) \\
   &+ \frac{\frac{\delta^{\frac{3}{2}}}{3b}}
           {\sqrt{\frac{4\delta^3}{9\sigma^4}-\frac{2b}{\sigma^2}}}
      \log\left(\left|\frac{(-\psi_{u,0}(t))^{\frac{1}{2}}
                            +\frac{2\delta^{\frac{3}{2}}}{3\sigma^2}
                            -\sqrt{\frac{4\delta^3}{9\sigma^4}-\frac{2b}{\sigma^2}}}
                           {(-\psi_{u,0}(t))^{\frac{1}{2}}
                            +\frac{2\delta^{\frac{3}{2}}}{3\sigma^2}
                            +\sqrt{\frac{4\delta^3}{9\sigma^4}-\frac{2b}{\sigma^2}}}\right|
          \right) \\
   &- \frac{\frac{\delta^{\frac{3}{2}}}{3b}}
           {\sqrt{\frac{4\delta^3}{9\sigma^4}-\frac{2b}{\sigma^2}}}
      \log\left(\left|\frac{(-u)^{\frac{1}{2}}+\frac{2\delta^{\frac{3}{2}}}{3\sigma^2}
                            -\sqrt{\frac{4\delta^3}{9\sigma^4}-\frac{2b}{\sigma^2}}}
                           {(-u)^{\frac{1}{2}}+\frac{2\delta^{\frac{3}{2}}}{3\sigma^2}
                            +\sqrt{\frac{4\delta^3}{9\sigma^4}-\frac{2b}{\sigma^2}}}\right|
          \right)
    = \frac{\sigma^2}{4} t .
  \end{aligned}
 \end{equation}
Finally, if \ $b = 0$, \ then, by separation of variables, we have
 \[
   \frac{1}{\bigl(g_u+\frac{4\delta^{\frac{3}{2}}}{3\sigma^2}\bigr)g_u^2} \dd g_u
   = - \frac{\sigma^2}{4} \dd t ,
 \]
 where
 \[
   \frac{1}{\bigl(g_u+\frac{4\delta^{\frac{3}{2}}}{3\sigma^2}\bigr)g_u^2}
   = \frac{\frac{9\sigma^4}{16\delta^3}}{g_u+\frac{4\delta^{\frac{3}{2}}}{3\sigma^2}}
     - \frac{\frac{9\sigma^4}{16\delta^3}}{g_u} + \frac{3\sigma^2}{4\delta^{\frac{3}{2}}g_u^2} ,
 \]
 hence
 \begin{align*}
  \int \frac{1}{\bigl(g_u+\frac{4\delta^{\frac{3}{2}}}{3\sigma^2}\bigr)g_u^2} \dd g_u
  &= \frac{9\sigma^4}{16\delta^3} \log\biggl(g_u+\frac{4\delta^{\frac{3}{2}}}{3\sigma^2}\biggr)
     - \frac{9\sigma^4}{16\delta^3} \log(g_u)
     - \frac{3\sigma^2}{4\delta^{\frac{3}{2}}g_u}
     + C \\
  &= \frac{9\sigma^4}{16\delta^3}
     \log\biggl(1 + \frac{4\delta^{\frac{3}{2}}}{3\sigma^2g_u}\biggr)
     - \frac{3\sigma^2}{4\delta^{\frac{3}{2}}g_u}
     + C ,
 \end{align*}
 yielding
 \[
   \frac{9\sigma^4}{16\delta^3}
   \log\biggl(1 + \frac{4\delta^{\frac{3}{2}}}{3\sigma^2g_u(t)}\biggr)
   - \frac{3\sigma^2}{4\delta^{\frac{3}{2}}g_u(t)}
   = - \frac{\sigma^2}{4} t + C ,\qquad t \in \RR_+ ,
 \]
 with some \ $C \in \RR$.
\ Using the initial value \ $g_u(0) = (-u)^{\frac{1}{2}}$, \ we obtain
 \[
   C = \frac{9\sigma^4}{16\delta^3}
       \log\biggl(1 + \frac{4\delta^{\frac{3}{2}}}{3\sigma^2(-u)^{\frac{1}{2}}}\biggr)
       - \frac{3\sigma^2}{4\delta^{\frac{3}{2}}(-u)^{\frac{1}{2}}} ,
 \]
 and hence, by \ $g_u(t) = (-\psi_{u,0}(t))^{\frac{1}{2}}$, \ we conclude
 \[
   \frac{9\sigma^4}{16\delta^3}
   \log\left(\frac{1+\frac{4\delta^{\frac{3}{2}}}{3\sigma^2(-\psi_{u,0}(t))^{\frac{1}{2}}}}
                  {1+\frac{4\delta^{\frac{3}{2}}}{3\sigma^2(-u)^{\frac{1}{2}}}}\right)
   - \frac{3\sigma^2}{4\delta^{\frac{3}{2}}}
     \left(\frac{1}{(-\psi_{u,0}(t))^{\frac{1}{2}}} - \frac{1}{(-u)^{\frac{1}{2}}}\right)
   = - \frac{\sigma^2}{4} t .
 \]
\end{Ex}

\begin{Ex}\label{Ex_Laplace_V_3/2}
We derive an explicit formula for the Laplace transform of \ $V$ \ given in Theorem
 \ref{Thm_supercritical_convergence}
 in case of \ $\alpha = \frac{3}{2}$.
\ In fact, we present two detailed calculations, the first one is based on the representation of
 \ $\psi_u^*$ \ given in part (iii) of Theorem \ref{Thm_supercritical_convergence}, and the second
 one is based on part (ii) of Theorem \ref{Thm_supercritical_convergence}.

\noindent {\sl Calculations based on part (iii) of Theorem \ref{Thm_supercritical_convergence}.}
We have
 \[
   \EE(\ee^{uV})
   = \exp\Biggl\{y_0 \psi_u^*
                 + a
                   \int_0^{-\psi_u^*}
                    \frac{1}
                         {\frac{\sigma^2}{2}x+\frac{2\delta^{\frac{3}{2}}}{3}x^{\frac{1}{2}}+b}
                    \, \dd x\Biggr\} , \qquad u \in \RR_- ,
 \]
 with \ $\psi_0^* = 0$ \ and \ $\psi_u^* = - K^{-1}(-u)$ \ for \ $u \in \RR_{--}$, \ where
 \ $K^{-1}$ \ is the inverse of the strictly increasing function
 \ $K : (0, \theta_0) \to \RR_{++}$ \ given by
 \[
   K(\lambda)
   = \lambda
     \exp\Biggl\{-\int_0^\lambda
                   \frac{\frac{\sigma^2}{2}+\frac{2\delta^{\frac{3}{2}}}{3}x^{-\frac{1}{2}}}
                        {\frac{\sigma^2}{2}x+\frac{2\delta^{\frac{3}{2}}}{3}x^{\frac{1}{2}}+b}
                   \, \dd x\Biggr\} , \qquad \lambda \in (0, \theta_0) ,
 \]
 where we used \eqref{b/R(z)-1/z} and
 \[
   \theta_0
   = \inf\biggl\{x \in \RR_{++}
                 : \frac{\sigma^2}{2} x^2 + \frac{2\delta^{\frac{3}{2}}}{3} x^{\frac{3}{2}}
                   + b x
                   \in \RR_+\biggr\}
   = \biggl(- \frac{2\delta^{\frac{3}{2}}}{3\sigma^2}
            + \sqrt{\frac{4\delta^3}{9\sigma^4} - \frac{2b}{\sigma^2}}\biggr)^2 .
 \]
By substitution \ $x = y^2$, \ for all \ $\lambda \in (0, \theta_0)$, \ we have
 \[
   \int_0^\lambda
    \frac{\frac{\sigma^2}{2}+\frac{2\delta^{\frac{3}{2}}}{3}x^{-\frac{1}{2}}}
         {\frac{\sigma^2}{2}x+\frac{2\delta^{\frac{3}{2}}}{3}x^{\frac{1}{2}}+b}
    \, \dd x
   = \int_0^{\lambda^{\frac{1}{2}}}
      \frac{\sigma^2y+\frac{4\delta^{\frac{3}{2}}}{3}}
           {\frac{\sigma^2}{2}y^2+\frac{2\delta^{\frac{3}{2}}}{3}y+b}
      \, \dd y .
 \]
First we consider the case of \ $\sigma \in \RR_{++}$.
\ Then we can write
 \[
   \frac{\sigma^2x+\frac{4\delta^{\frac{3}{2}}}{3}}
        {\frac{\sigma^2}{2}x^2+\frac{2\delta^{\frac{3}{2}}}{3}x+b}
   = \frac{2x+\frac{4\delta^{\frac{3}{2}}}{3\sigma^2}}
          {x^2+\frac{4\delta^{\frac{3}{2}}}{3\sigma^2}x+\frac{2b}{\sigma^2}}
     + \frac{4\delta^{\frac{3}{2}}}{3\sigma^2}
       \cdot \frac{1}{x^2+\frac{4\delta^{\frac{3}{2}}}{3\sigma^2}x+\frac{2b}{\sigma^2}} .
 \]
As in Example \ref{Ex_Laplace_stat_3/2}, we have
 \[
   \int_0^{\lambda^{\frac{1}{2}}}
    \frac{2x+\frac{4\delta^{\frac{3}{2}}}{3\sigma^2}}
         {x^2+\frac{4\delta^{\frac{3}{2}}}{3\sigma^2}x+\frac{2b}{\sigma^2}}
    \, \dd x
   = \log\biggl(\frac{\sigma^2}{2b} \lambda
                + \frac{2\delta^{\frac{3}{2}}}{3b} \lambda^{\frac{1}{2}} + 1\biggr)
 \]
 and
 \[
   \int_0^{\lambda^{\frac{1}{2}}}
    \frac{1}{x^2+\frac{4\delta^{\frac{3}{2}}}{3\sigma^2}x+\frac{2b}{\sigma^2}} \, \dd x
   = \frac{1}{\sqrt{\frac{4\delta^3}{9\sigma^4}-\frac{2b}{\sigma^2}}}
     \log\left(\frac{\lambda^{\frac{1}{2}}
                     \Bigl(\frac{2\delta^{\frac{3}{2}}}{3\sigma^2}
                           +\sqrt{\frac{4\delta^3}{9\sigma^4}-\frac{2b}{\sigma^2}}\Bigr)
                     +\frac{2b}{\sigma^2}}
                    {\lambda^{\frac{1}{2}}
                     \Bigl(\frac{2\delta^{\frac{3}{2}}}{3\sigma^2}
                           -\sqrt{\frac{4\delta^3}{9\sigma^4}-\frac{2b}{\sigma^2}}\Bigr)
                     +\frac{2b}{\sigma^2}}\right) ,
 \]
 and hence
 \begin{align*}
  K(\lambda)
  &= \frac{\lambda}
          {\frac{\sigma^2}{2b}\lambda+\frac{2\delta^{\frac{3}{2}}}{3b}\lambda^{\frac{1}{2}}+1}
     \left(\frac{\lambda^{\frac{1}{2}}
                 \Bigl(\frac{2\delta^{\frac{3}{2}}}{3\sigma^2}
                       +\sqrt{\frac{4\delta^3}{9\sigma^4}-\frac{2b}{\sigma^2}}\Bigr)
                 +\frac{2b}{\sigma^2}}
                {\lambda^{\frac{1}{2}}
                 \Bigl(\frac{2\delta^{\frac{3}{2}}}{3\sigma^2}
                       -\sqrt{\frac{4\delta^3}{9\sigma^4}-\frac{2b}{\sigma^2}}\Bigr)
                 +\frac{2b}{\sigma^2}}
     \right)^{-\frac{\frac{4\delta^{\frac{3}{2}}}{3\sigma^2}}
                    {\sqrt{\frac{4\delta^3}{9\sigma^4}-\frac{2b}{\sigma^2}}}} \\
  &= \biggl(\frac{\delta^{\frac{3}{2}}}{3b}+\sqrt{\frac{\delta^3}{9b^2}-\frac{\sigma^2}{2b}}
            +\lambda^{-\frac{1}{2}}
     \biggr)^{-1-\frac{\frac{4\delta^{\frac{3}{2}}}{3\sigma^2}}
                      {\sqrt{\frac{4\delta^3}{9\sigma^4}-\frac{2b}{\sigma^2}}}}
     \biggl(\frac{\delta^{\frac{3}{2}}}{3b}-\sqrt{\frac{\delta^3}{9b^2}-\frac{\sigma^2}{2b}}
            +\lambda^{-\frac{1}{2}}
     \biggr)^{-1+\frac{\frac{4\delta^{\frac{3}{2}}}{3\sigma^2}}
                      {\sqrt{\frac{4\delta^3}{9\sigma^4}-\frac{2b}{\sigma^2}}}} .
 \end{align*}
Note that an explicit formula for \ $K^{-1}$ \ is not available.
Next we consider the case of \ $\sigma = 0$.
\ By Remark \ref{Thm_supercritical_convergence_sigma=0}, we have
 \ $\theta_0 = \frac{9b^2}{4\delta^3}$,
 \begin{align*}
  K(\lambda)
  = \biggl(\frac{2\delta^{\frac{3}{2}}}{3b} + \lambda^{-\frac{1}{2}}\biggr)^{-2} ,
  \qquad \lambda \in \biggl(0, \frac{9b^2}{4\delta^3}\biggr) ,
 \end{align*}
 thus
 \[
   \psi_u^*
   = - K^{-1}(-u)
   = - \biggl((-u)^{-\frac{1}{2}} - \frac{2\delta^{\frac{3}{2}}}{3b}\biggr)^{-2} ,
 \]
 and hence
 \[
   \EE(\ee^{uV})
   = \exp\Biggl\{- y_0
                   \biggl((-u)^{-\frac{1}{2}} - \frac{2\delta^{\frac{3}{2}}}{3b}\biggr)^{-2}
                 + a \int_0^{\bigl((-u)^{-\frac{1}{2}}
                                   -\frac{2\delta^{\frac{3}{2}}}{3b}\bigr)^{-2}}
                      \frac{1}{\frac{2\delta^{\frac{3}{2}}}{3}x^{\frac{1}{2}}+b}
                      \, \dd x\Biggr\} ,
   \qquad u \in \RR_{--} .
 \]
By substitution \ $x = y^2$, \ for all \ $u \in \RR_{--}$, \ we have
 \[
   \int_0^{\bigl((-u)^{-\frac{1}{2}}-\frac{2\delta^{\frac{3}{2}}}{3b}\bigr)^{-2}}
    \frac{1}{\frac{2\delta^{\frac{3}{2}}}{3}x^{\frac{1}{2}}+b} \, \dd x
   = \int_0^{\bigl((-u)^{-\frac{1}{2}}-\frac{2\delta^{\frac{3}{2}}}{3b}\bigr)^{-1}}
      \frac{2y}{\frac{2\delta^{\frac{3}{2}}}{3}y+b} \, \dd y .
 \]
As in Example \ref{Ex_Laplace_stat_3/2}, we can write
 \[
   \frac{2y}{\frac{2\delta^{\frac{3}{2}}}{3}y+b}
   = \frac{3}{\delta^{\frac{3}{2}}}
     - \frac{\frac{9b}{2\delta^3}}{y+\frac{3b}{2\delta^{\frac{3}{2}}}} ,
 \]
 thus
 \begin{align*}
  \int_0^{\bigl((-u)^{-\frac{1}{2}}-\frac{2\delta^{\frac{3}{2}}}{3b}\bigr)^{-1}}
   \frac{2y}{\frac{2\delta^{\frac{3}{2}}}{3}y+b} \, \dd y
  &= \frac{3}{\delta^{\frac{3}{2}}}
     \biggl((-u)^{-\frac{1}{2}}-\frac{2\delta^{\frac{3}{2}}}{3b}\biggr)^{-1} \\
  &\quad
     - \frac{9b}{2\delta^3}
       \biggl(\log\biggl(\biggl|\biggl((-u)^{-\frac{1}{2}}
                                       -\frac{2\delta^{\frac{3}{2}}}{3b}\biggr)^{-1}
                                +\frac{3b}{2\delta^{\frac{3}{2}}}\biggr|\biggr)
              - \log\biggl(\biggl|\frac{3b}{2\delta^{\frac{3}{2}}}\biggr|\biggr)\biggr) \\
  &= \frac{3}{\delta^{\frac{3}{2}}}
     \biggl((-u)^{-\frac{1}{2}}-\frac{2\delta^{\frac{3}{2}}}{3b}\biggr)^{-1}
     + \frac{9b}{2\delta^3}
       \log\biggl(1 - \frac{2\delta^{\frac{3}{2}}}{3b} (-u)^{\frac{1}{2}}\biggr) ,
 \end{align*}
 and hence
 \[
   \EE(\ee^{uV})
   = \exp\biggl\{- y_0
                   \biggl((-u)^{-\frac{1}{2}}
                          - \frac{2\delta^{\frac{3}{2}}}{3b}\biggr)^{-2}\biggr\}
     \exp\biggl\{\frac{3a}{\delta^{\frac{3}{2}}}
                 \biggl((-u)^{-\frac{1}{2}}-\frac{2\delta^{\frac{3}{2}}}{3b}\biggr)^{-1}\biggr\}
     \biggl(1 - \frac{2\delta^{\frac{3}{2}}}{3b}
                (-u)^{\frac{1}{2}}\biggr)^{\frac{9ba}{2\delta^3}}
 \]
 for \ $u \in \RR_{--}$.

\noindent
{\sl Calculations based on part (ii) of Theorem \ref{Thm_supercritical_convergence}.}
By \eqref{Laplace_supercritical_limit}, it is enough to know \ $\psi_u^*$ \ to have
 an explicit formula for the Laplace transform of \ $V$.
\ We carry out this calculation only in case of \ $\sigma \in \RR_{++}$.
\ By \eqref{psi_supercrit}, for all \ $t \in \RR_+$ \ and \ $u \in \RR_{--}$ \ with
 \ $u \ee^{bt} \ne u_0$, \ we obtain
 \begin{align*}
  &- \frac{\sigma^2}{2b} \log((-\psi_{u\ee^{bt},0}(t))^{\frac{1}{2}})
   + \frac{\sigma^2}{4b}
     \log\biggl(\biggl|- \psi_{u\ee^{bt},0}(t)
                       + \frac{4\delta^{\frac{3}{2}}}{3\sigma^2}
                         (-\psi_{u\ee^{bt},0}(t))^{\frac{1}{2}}
                       + \frac{2b}{\sigma^2}\biggr|\biggr) \\
  &+ \frac{\sigma^2}{2b} \log((-u\ee^{bt})^{\frac{1}{2}})
   - \frac{\sigma^2}{4b}
     \log\biggl(\biggl|-u\ee^{bt}
                       + \frac{4\delta^{\frac{3}{2}}}{3\sigma^2} (-u\ee^{bt})^{\frac{1}{2}}
                       + \frac{2b}{\sigma^2}\biggr|\biggr) \\
  &+ \frac{\frac{\delta^{\frac{3}{2}}}{3b}}
          {\sqrt{\frac{4\delta^3}{9\sigma^4}-\frac{2b}{\sigma^2}}}
     \log\left(\left|\frac{(-\psi_{u\ee^{bt},0}(t))^{\frac{1}{2}}
                           +\frac{2\delta^{\frac{3}{2}}}{3\sigma^2}
                           -\sqrt{\frac{4\delta^3}{9\sigma^4}-\frac{2b}{\sigma^2}}}
                          {(-\psi_{u\ee^{bt},0}(t))^{\frac{1}{2}}
                           +\frac{2\delta^{\frac{3}{2}}}{3\sigma^2}
                           +\sqrt{\frac{4\delta^3}{9\sigma^4}-\frac{2b}{\sigma^2}}}\right|
         \right) \\
  &- \frac{\frac{\delta^{\frac{3}{2}}}{3b}}
          {\sqrt{\frac{4\delta^3}{9\sigma^4}-\frac{2b}{\sigma^2}}}
     \log\left(\left|\frac{(-u\ee^{bt})^{\frac{1}{2}}+\frac{2\delta^{\frac{3}{2}}}{3\sigma^2}
                           -\sqrt{\frac{4\delta^3}{9\sigma^4}-\frac{2b}{\sigma^2}}}
                          {(-u\ee^{bt})^{\frac{1}{2}}+\frac{2\delta^{\frac{3}{2}}}{3\sigma^2}
                           +\sqrt{\frac{4\delta^3}{9\sigma^4}-\frac{2b}{\sigma^2}}}\right|
         \right)
   = \frac{\sigma^2}{4} t .
 \end{align*}
Using
 \ $\frac{\sigma^2}{2b} \log((-u\ee^{bt})^{\frac{1}{2}})
    = \frac{\sigma^2}{4b} \log(-u) + \frac{\sigma^2}{4} t$,
 \ we conclude that \ $\psi_u^*$ \ satisfies the equation
 \begin{equation}\label{|psi_u^*|}
  \begin{aligned}
   &- \frac{\sigma^2}{4} \log(-\psi_u^*)
    + \frac{\sigma^2}{4}
      \log\biggl(\biggl|- \psi_u^*
                        + \frac{4\delta^{\frac{3}{2}}}{3\sigma^2} (-\psi_u^*)^{\frac{1}{2}}
                        + \frac{2b}{\sigma^2}\biggr|\biggr)
    + \frac{\sigma^2}{4} \log(-u)
    - \frac{\sigma^2}{4} \log\biggl(-\frac{2b}{\sigma^2}\biggr) \\
   &+ \frac{\frac{\delta^{\frac{3}{2}}}{3}}
           {\sqrt{\frac{4\delta^3}{9\sigma^4}-\frac{2b}{\sigma^2}}}
      \left(\log\left(\left|\frac{(-\psi_u^*)^{\frac{1}{2}}
                                  +\frac{2\delta^{\frac{3}{2}}}{3\sigma^2}
                                  -\sqrt{\frac{4\delta^3}{9\sigma^4}-\frac{2b}{\sigma^2}}}
                                 {(-\psi_u^*)^{\frac{1}{2}}
                                  +\frac{2\delta^{\frac{3}{2}}}{3\sigma^2}
                                  +\sqrt{\frac{4\delta^3}{9\sigma^4}
                                  -\frac{2b}{\sigma^2}}}\right|\right)
            - \log\left(\frac{\sqrt{\frac{4\delta^3}{9\sigma^4}-\frac{2b}{\sigma^2}}
                              -\frac{2\delta^{\frac{3}{2}}}{3\sigma^2}}
                             {\sqrt{\frac{4\delta^3}{9\sigma^4}-\frac{2b}{\sigma^2}}
                              +\frac{2\delta^{\frac{3}{2}}}{3\sigma^2}}\right)\right)
    = 0 .
  \end{aligned}
 \end{equation}
By Theorem \ref{Thm_supercritical_convergence}, \ $R(-\psi_u^*) < 0$.
\ We have
 \ $R(-\psi_u^*)
    = \frac{\sigma^2}{2}
      \bigl(- \psi_u^* + \frac{4\delta^{\frac{3}{2}}}{3\sigma^2} (-\psi_u^*)^{\frac{1}{2}}
            + \frac{2b}{\sigma^2}\bigr)
      (-\psi_u^*)$,
 \ where \ $-\psi_u^* \in (0,\theta_0)$, \ hence
 \ $- \psi_u^* + \frac{4\delta^{\frac{3}{2}}}{3\sigma^2} (-\psi_u^*)^{\frac{1}{2}}
    + \frac{2b}{\sigma^2} < 0$.
\ Moreover,
 \[
   - \psi_u^* + \frac{4\delta^{\frac{3}{2}}}{3\sigma^2} (-\psi_u^*)^{\frac{1}{2}}
   + \frac{2b}{\sigma^2}
   = \biggl((-\psi_u^*)^{\frac{1}{2}} + \frac{2\delta^{\frac{3}{2}}}{3\sigma^2}
            - \sqrt{\frac{4\delta^3}{9\sigma^4}-\frac{2b}{\sigma^2}}\biggr)
     \biggl((-\psi_u^*)^{\frac{1}{2}} + \frac{2\delta^{\frac{3}{2}}}{3\sigma^2}
            + \sqrt{\frac{4\delta^3}{9\sigma^4}-\frac{2b}{\sigma^2}}\biggr) ,
 \]
 where
 \ $(-\psi_u^*)^{\frac{1}{2}} + \frac{2\delta^{\frac{3}{2}}}{3\sigma^2}
    + \sqrt{\frac{4\delta^3}{9\sigma^4}-\frac{2b}{\sigma^2}} > 0$,
 \ thus
 \ $(-\psi_u^*)^{\frac{1}{2}} + \frac{2\delta^{\frac{3}{2}}}{3\sigma^2}
    - \sqrt{\frac{4\delta^3}{9\sigma^4}-\frac{2b}{\sigma^2}} < 0$.
Consequently, \ $\psi_u^*$ \ satisfies the equation
 \begin{equation}\label{psi_u^*}
  \begin{aligned}
   &- \frac{\sigma^2}{4} \log(-\psi_u^*)
    + \frac{\sigma^2}{4}
      \log\biggl(\psi_u^*
                 - \frac{4\delta^{\frac{3}{2}}}{3\sigma^2} (-\psi_u^*)^{\frac{1}{2}}
                 - \frac{2b}{\sigma^2}\biggr)
    + \frac{\sigma^2}{4} \log(-u)
    - \frac{\sigma^2}{4} \log\biggl(-\frac{2b}{\sigma^2}\biggr) \\
   &+ \frac{\frac{\delta^{\frac{3}{2}}}{3}}
           {\sqrt{\frac{4\delta^3}{9\sigma^4}-\frac{2b}{\sigma^2}}}
      \left(\log\left(\frac{\sqrt{\frac{4\delta^3}{9\sigma^4}-\frac{2b}{\sigma^2}}
                            -(-\psi_u^*)^{\frac{1}{2}}
                            -\frac{2\delta^{\frac{3}{2}}}{3\sigma^2}}
                           {\sqrt{\frac{4\delta^3}{9\sigma^4}-\frac{2b}{\sigma^2}}
                            +(-\psi_u^*)^{\frac{1}{2}}
                            +\frac{2\delta^{\frac{3}{2}}}{3\sigma^2}}\right)
            - \log\left(\frac{\sqrt{\frac{4\delta^3}{9\sigma^4}-\frac{2b}{\sigma^2}}
                              -\frac{2\delta^{\frac{3}{2}}}{3\sigma^2}}
                             {\sqrt{\frac{4\delta^3}{9\sigma^4}-\frac{2b}{\sigma^2}}
                              +\frac{2\delta^{\frac{3}{2}}}{3\sigma^2}}\right)\right)
    = 0 .
  \end{aligned}
 \end{equation}
Note that, by Theorem \ref{Thm_supercritical_convergence}, \eqref{psi_u^*} is equivalent to
 \ $K(-\psi_u^*) = -u$.

We show another way to derive this equation.
By Theorem \ref{Thm_supercritical_convergence}, for all sufficiently small
 \ $\lambda \in \RR_{++}$,
 \begin{equation}\label{psi_u*_K}
  \begin{aligned}
   &\int_\lambda^{-\psi_u^*}
     \frac{1}{\frac{\sigma^2}{2}z^2+\frac{2\delta^{\frac{3}{2}}}{3}z^{\frac{3}{2}}+bz} \, \dd z
    = \int_\lambda^{f\big(\frac{-u}{K(\lambda)},\lambda\big)} \frac{1}{R(z)} \, \dd z
    = \frac{1}{b} \log\biggl(-\frac{u}{K(\lambda)}\biggr)\\
   &= \frac{1}{b} \log(-u) - \frac{1}{b} \log(K(\lambda))
    = \frac{1}{b} \log(-u) - \frac{1}{b} \log(\lambda)
      - \frac{1}{b} \int_0^\lambda \biggl(\frac{b}{R(z)} - \frac{1}{z}\biggr) \dd z .
  \end{aligned}
 \end{equation}
By substitution \ $z = y^2$, \ and using \eqref{help_alpha_3/2_2},
 \begin{align*}
  &\int_\lambda^{-\psi_u^*}
    \frac{1}{\frac{\sigma^2}{2}z^2+\frac{2\delta^{\frac{3}{2}}}{3}z^{\frac{3}{2}}+bz} \, \dd z
   = \int_{\sqrt{\lambda}}^{(-\psi_u^*)^{\frac{1}{2}}}
      \frac{2}{\frac{\sigma^2}{2}y^3+\frac{2\delta^{\frac{3}{2}}}{3}y^2+by} \, \dd y \\
  &= \left[\frac{2}{b} \log(|y|)
           - \frac{1}{b}
             \log\biggl(\biggl|y^2 + \frac{4\delta^{\frac{3}{2}}}{3\sigma^2} y
                               + \frac{2b}{\sigma^2}\biggr|\biggr)
           - \frac{\frac{4\delta^{\frac{3}{2}}}{3\sigma^2b}}
                  {\sqrt{\frac{4\delta^3}{9\sigma^4}-\frac{2b}{\sigma^2}}}
             \log\left(\left|\frac{y+\frac{2\delta^{\frac{3}{2}}}{3\sigma^2}
                                   -\sqrt{\frac{4\delta^3}{9\sigma^4}-\frac{2b}{\sigma^2}}}
                                  {y+\frac{2\delta^{\frac{3}{2}}}{3\sigma^2}
                                   +\sqrt{\frac{4\delta^3}{9\sigma^4}
                                          -\frac{2b}{\sigma^2}}}\right|\right)
     \right]_{y=\sqrt{\lambda}}^{y=(-\psi_u^*)^{\frac{1}{2}}} \\
  &= \frac{1}{b} \log(-\psi_u^*) - \frac{1}{b} \log(\lambda)
     - \frac{1}{b}
       \log\biggl(\biggl|- \psi_u^*
                         + \frac{4\delta^{\frac{3}{2}}}{3\sigma^2} (-\psi_u^*)^{\frac{1}{2}}
                         + \frac{2b}{\sigma^2}\biggr|\biggr)
     + \frac{1}{b}
       \log\biggl(\biggl|\lambda
                         + \frac{4\delta^{\frac{3}{2}}}{3\sigma^2} \sqrt{\lambda}
                         + \frac{2b}{\sigma^2}\biggr|\biggr) \\
  &\quad
     - \frac{\frac{4\delta^{\frac{3}{2}}}{3\sigma^2b}}
                  {\sqrt{\frac{4\delta^3}{9\sigma^4}-\frac{2b}{\sigma^2}}}
       \left(\log\left(\left|\frac{(-\psi_u^*)^{\frac{1}{2}}
                                   +\frac{2\delta^{\frac{3}{2}}}{3\sigma^2}
                                   -\sqrt{\frac{4\delta^3}{9\sigma^4}-\frac{2b}{\sigma^2}}}
                                  {(-\psi_u^*)^{\frac{1}{2}}
                                   +\frac{2\delta^{\frac{3}{2}}}{3\sigma^2}
                                   +\sqrt{\frac{4\delta^3}{9\sigma^4}
                                          -\frac{2b}{\sigma^2}}}\right|\right)
             - \log\left(\left|\frac{\sqrt{\lambda}
                                     +\frac{2\delta^{\frac{3}{2}}}{3\sigma^2}
                                     -\sqrt{\frac{4\delta^3}{9\sigma^4}-\frac{2b}{\sigma^2}}}
                                    {\sqrt{\lambda}
                                     +\frac{2\delta^{\frac{3}{2}}}{3\sigma^2}
                                     +\sqrt{\frac{4\delta^3}{9\sigma^4}
                                          -\frac{2b}{\sigma^2}}}\right|\right)\right) .
 \end{align*}
Moreover,
 \[
   \int_0^\lambda \biggl(\frac{b}{R(z)} - \frac{1}{z}\biggr) \dd z
   = \int_0^\lambda
      \left(\frac{b}{\frac{\sigma^2}{2}z^2+\frac{2\delta^{\frac{3}{2}}}{3}z^{\frac{3}{2}}+bz}
            - \frac{1}{z}\right)
      \dd z
   = - \int_0^\lambda
        \frac{\frac{\sigma^2}{2}+\frac{2\delta^{\frac{3}{2}}}{3}z^{-\frac{1}{2}}}
             {\frac{\sigma^2}{2}z+\frac{2\delta^{\frac{3}{2}}}{3}z^{\frac{1}{2}}+b}
        \, \dd z .
 \]
By substitution \ $z = y^2$,
 \[
   \int_0^\lambda
    \frac{\frac{\sigma^2}{2}+\frac{2\delta^{\frac{3}{2}}}{3}z^{-\frac{1}{2}}}
         {\frac{\sigma^2}{2}z+\frac{2\delta^{\frac{3}{2}}}{3}z^{\frac{1}{2}}+b}
    \, \dd z
   = \int_0^{\sqrt{\lambda}}
      \frac{\sigma^2y+\frac{4\delta^{\frac{3}{2}}}{3}}
           {\frac{\sigma^2}{2}y^2+\frac{2\delta^{\frac{3}{2}}}{3}y+b}
      \, \dd y .
 \]
We can write
 \[
   \frac{\sigma^2y+\frac{4\delta^{\frac{3}{2}}}{3}}
        {\frac{\sigma^2}{2}y^2+\frac{2\delta^{\frac{3}{2}}}{3}y+b}
   = \frac{\sigma^2y+\frac{2\delta^{\frac{3}{2}}}{3}}
          {\frac{\sigma^2}{2}y^2+\frac{2\delta^{\frac{3}{2}}}{3}y+b}
     + \frac{\frac{2\delta^{\frac{3}{2}}}{3}}
            {\frac{\sigma^2}{2}y^2+\frac{2\delta^{\frac{3}{2}}}{3}y+b} ,
 \]
 hence, by \eqref{help_alpha_3/2_1},
 \begin{align*}
  &\int_0^\lambda \biggl(\frac{b}{R(z)} - \frac{1}{z}\biggr) \dd z
   = - \int_0^{\sqrt{\lambda}}
        \frac{\sigma^2y+\frac{4\delta^{\frac{3}{2}}}{3}}
             {\frac{\sigma^2}{2}y^2+\frac{2\delta^{\frac{3}{2}}}{3}y+b}
        \, \dd y \\
  &= - \left[\log\biggl(\biggl|\frac{\sigma^2}{2}y^2
                               + \frac{2\delta^{\frac{3}{2}}}{3}y + b\biggr|\biggr)
             + \frac{\frac{4\delta^{\frac{3}{2}}}{3\sigma^2}}
                    {\sqrt{\frac{4\delta^3}{9\sigma^4}-\frac{2b}{\sigma^2}}}
               \log\left(\left|\frac{y+\frac{2\delta^{\frac{3}{2}}}{3\sigma^2}
                                     -\sqrt{\frac{4\delta^3}{9\sigma^4}-\frac{2b}{\sigma^2}}}
                                    {y+\frac{2\delta^{\frac{3}{2}}}{3\sigma^2}
                                     +\sqrt{\frac{4\delta^3}{9\sigma^4}
                                            -\frac{2b}{\sigma^2}}}\right|\right)
       \right]_{y=0}^{y=\sqrt{\lambda}} \\
  &= - \log\biggl(\biggl|\frac{\sigma^2}{2} \lambda
                         + \frac{2\delta^{\frac{3}{2}}}{3} \sqrt{\lambda} + b\biggr|\biggr)
     + \log(-b) \\
  &\quad
     - \frac{\frac{4\delta^{\frac{3}{2}}}{3\sigma^2}}
            {\sqrt{\frac{4\delta^3}{9\sigma^4}-\frac{2b}{\sigma^2}}}
       \left(\log\left(\left|\frac{\sqrt{\lambda}+\frac{2\delta^{\frac{3}{2}}}{3\sigma^2}
                                   -\sqrt{\frac{4\delta^3}{9\sigma^4}-\frac{2b}{\sigma^2}}}
                                  {\sqrt{\lambda}+\frac{2\delta^{\frac{3}{2}}}{3\sigma^2}
                                   +\sqrt{\frac{4\delta^3}{9\sigma^4}
                                    -\frac{2b}{\sigma^2}}}\right|\right)
             - \log\left(\frac{\sqrt{\frac{4\delta^3}{9\sigma^4}-\frac{2b}{\sigma^2}}
                               -\frac{2\delta^{\frac{3}{2}}}{3\sigma^2}}
                              {\sqrt{\frac{4\delta^3}{9\sigma^4}-\frac{2b}{\sigma^2}}
                               +\frac{2\delta^{\frac{3}{2}}}{3\sigma^2}}\right)\right) .
 \end{align*}
Consequently, \eqref{psi_u*_K} yields again that \ $\psi_u^*$ \ satisfies equation
 \eqref{|psi_u^*|}, and hence, equation \eqref{psi_u^*}.
\end{Ex}

\section*{Acknowledgements}

We are grateful to Cl\'ement Foucart for providing us an idea how to derive
 \eqref{Laplace_supercritical_limit}, a formula for the Laplace transform of \ $V$ \ in Theorem
 \ref{Thm_supercritical_convergence}.
We would like to thank Hatem Zaag for explaining us several methods that may be used for
 describing the asymptotic behavior of the ordinary differential equation \eqref{psi_DE}.
We would like to thank the referees for their comments that helped us to improve the paper.

\end{document}